\documentclass[10pt]{amsart}
\usepackage[english]{babel}
\usepackage{graphicx}
\usepackage{amsmath,amssymb,hyperref}
\newcommand{\diff}[2]{\mbox{{\rm Diff}{${\,}_{#1}({\mathbb
C}^{#2},0)$}}}
\newcommand{\diffh}[2]{\mbox{$\widehat{\rm Diff}{{\,}_{#1}({\mathbb
C}^{#2},0)}$}}

\newcommand{\cn}[1]{\mbox{(${\mathbb C}^{#1},0$)}}
\newcommand{\Xt}{{\mathcal X}_{tp1}\cn{2}}
\newcommand{\Xnt}{{\mathcal X}_{p1} \cn{2}}

\newcommand{\pn}[1]{{\mathbb P}^{#1}({\mathbb C})}

\newcommand{\ex}{\'{e}}

\newcommand{\ox}{\'{o}}

\newtheorem{pro}{Proposition}[section]
\newtheorem{teo}{Theorem}[section]
\newtheorem*{main}{Main Theorem}
\newtheorem{cor}{Corollary}[section]

\newtheorem{lem}{Lemma}[section]
\newtheorem{rem}{Remark}[section]

\begin{document}

\title[Unfoldings of resonant diffeomorphisms]{Modulus of analytic
classification for unfoldings of resonant diffeomorphisms}
\author{Javier Rib\ox n}
\thanks{IMPA, Estrada Dona Castorina 110, Rio de Janeiro, Brasil, 
22460-320}
\thanks{e-mail address: jfribon@impa.br}
\thanks{MSC-class. Primary: 37F45; Secondary: 37G10, 37F75}
\date{\today}
\maketitle

\bibliographystyle{plain}
\section*{Abstract}
We provide a complete system of analytic invariants for unfoldings
of non-linearizable resonant complex analytic diffeomorphisms
as well as its geometrical interpretation. In order to fulfill this goal 
we develop an extension of the Fatou coordinates with controlled
asymptotic behavior in the neighborhood of the fixed points. 
The classical constructions are based on finding regions where
the dynamics of the unfolding is topologically stable. We introduce
a concept of infinitesimal stability leading to Fatou coordinates
reflecting more faithfully the analytic nature of the unfolding.   
These improvements allow us to control the domain of definition of
a conjugating mapping and its power series expansion. 
\section{Introduction}
In this paper we provide a complete analytic classification for unfoldings
of non-linearizable resonant complex analytic diffeomorphisms. The group of
1-dimensional unfoldings of elements of $\diff{}{}$ is
\[ \diff{p}{2} = \{ \varphi(x,y) \in \diff{}{2} \ : \ x \circ \varphi = x \}
. \]
Our main result is stated in the set $\diff{p1}{2}$ composed by the elements 
$\varphi$ of $\diff{p}{2}$ such that 
$\varphi_{|x=0}$ is tangent to the identity 
(i.e. $j^{1} \varphi_{|x=0} = Id$) but $\varphi_{|x=0} \neq Id$.
Given $\varphi_{1}, \varphi_{2}$ we denote
$\varphi_{1} \sim \varphi_{2}$ if they share the same set of fixed points
$Fix \varphi_{1} = Fix \varphi_{2}$ and $\varphi_{1}, \varphi_{2}$ are 
conjugated
by a holomorphic diffeomorphism respecting the fixed points and the fibers
$x=constant$.
\begin{main}
\label{teo:modi}
Consider $\varphi_{1}, \varphi_{2} \in \diff{p1}{2}$ with
$Fix \varphi_{1}= Fix \varphi_{2}$.
Then $\varphi_{1} \sim \varphi_{2}$ if and only if
there exists a real constant $r \in {\mathbb R}^{+}$ such that
for all $x_{0}$ in a pointed neighborhood of $0$ the restrictions
$(\varphi_{1})_{|x=x_{0}}$ and $(\varphi_{2})_{|x=x_{0}}$ are conjugated
by an injective holomorphic mapping defined in $B(0,r)$.
\end{main}
There is no hypothesis on the dependance on
$x_{0}$ of the analytic mappings conjugating
$(\varphi_{1})_{|x=x_{0}}$ and $(\varphi_{2})_{|x=x_{0}}$.
We only require a uniform domain of definition.

The connection between the main theorem and a system of analytic invariants
is obtained by getting an extension of the Fatou coordinates of
$\varphi_{|x=0}$ to the nearby parameters for $\varphi \in \diff{p1}{2}$.
In order to prove the main result we need an improvement of the classical
constructions
(Lavaurs \cite{Lavaurs}-Shishikura \cite{Shishi}-Oudkerk \cite{Oudkerk})
reflecting better the geometry of $\varphi$.

Our system of invariants is a generalization of the 
Mardesic-Roussarie-Rousseau's
system \cite{MRR} for generic unfoldings of codimension $1$ 
tangent to the identity elements of $\diff{}{}$.
We drop here the hypotheses on generic character and codimension.

Let us precise our main statement. Consider the set $\diff{pr}{2}$ of
$\varphi$ in $\diff{p}{2}$ such that $j^{1} \varphi_{|x=0}$ is periodic but
$\varphi_{|x=0}$ is not.
  As a consequence of the Jordan-Chevalley decomposition in linear algebraic
groups the analytic invariants of $\varphi \in \diff{pr}{2}$ coincide with
those of an iterate $\varphi^{\circ (q)} \in \diff{p1}{2}$.

Mardesic, Roussarie and Rousseau apply a refinement of Shishikura's
construction \cite{Shishi}
to get extensions of the Fatou coordinates supported in Lavaurs
sectors $V_{\delta}^{L}$ describing an angle as close to $4 \pi$ as
desired in the $x$-variable. Indeed the extensions
are multi-valuated around $x=0$.
They define analytic invariants a
la Martinet-Ramis.
More precisely they define a classifying space ${\mathcal M}$ and a mapping
$m_{\varphi}:V_{\delta}^{L} \to {\mathcal M}$. They claim that 
$\varphi \sim \zeta$ is equivalent to
$m_{\varphi} \equiv m_{\zeta}$. We skip the details of
the definition of $m_{\varphi}$ but we stress that
$m_{\varphi}(x_{0})$ depends only on $\varphi_{|x=x_{0}}$.
We generalize the definition of $m_{\varphi}$ for all
$\varphi \in \diff{p1}{2}$.
Their result induces to think that the uniform hypothesis in our main
theorem is superfluous. That is not the case, we provide a counterexample
to the main theorem in \cite{MRR}.
\begin{teo}
\label{teoi:couMRR}
There exist $\varphi, \zeta \in \diff{p1}{2}$ such that
\begin{enumerate}
\item $Fix \varphi=Fix \zeta$ and $\varphi_{|x=0} \equiv \zeta_{|x=0}$.
\item $\varphi$, $\zeta$ are conjugated by an injective
analytic $\sigma$ defined in $|y| < C_{0}/\sqrt[\nu]{|\ln x|}$
such that $\sigma(e^{2 \pi i}x,y)= \zeta \circ \sigma(x,y)$
for some $(C_{0}, \nu)  \in {\mathbb R}^{+} \times {\mathbb N}$.
\end{enumerate}
but $\varphi_{1} \not \sim \varphi_{2}$. Given $\eta \in \diff{p1}{2}$ we
can suppose that $(y \circ \varphi -y) = (y \circ \eta -y)$.
\end{teo}
In particular the counterexample is obtained by fixing
$(y \circ \varphi-y)=(x-y^{2})$, the conditions 1 and 2 imply
$m_{\varphi} \equiv m_{\zeta}$.
Their statement can be easily repaired, it is just too optimistic.
In this paper we do it in two different ways: by giving a uniform version
of the analytic system of invariants and also by studying under what
rigidity conditions the uniform hypothesis is no longer necessary.
\begin{teo}
\label{teoi:rig}
Let $\varphi, \zeta$ be formally conjugated elements of $\diff{p1}{2}$
such that $Fix \varphi = Fix \zeta$. Suppose that
$\varphi_{|x=0}$ is not analytically trivial. Then $\varphi \sim \zeta$
if and only if $m_{\varphi} \equiv m_{\zeta}$.
\end{teo}
Denote $\diff{1}{} = \{ \phi \in \diff{}{} : j^{1} \phi = Id\}$.
We say that $\varphi \in \diff{p1}{2}$ (resp. $\phi \in \diff{1}{}$) 
is analytically trivial if it is the exponential
of a germ of nilpotent vector field. A consequence of theorem \ref{teoi:rig}
is that the main theorem in \cite{MRR} is valid in the generic case.

The complete system of analytic invariants is based on building Fatou
coordinates for $\varphi \in \diff{p1}{2}$. Shishikura \cite{Shishi} considers
``transversals'' to the dynamics of $\varphi$. A transversal $T$ and its image
$\varphi(T)$ enclose a strip $S(T)$. The space of orbits of $\varphi_{|S(T)}$
is biholomorphic to ${\mathbb C}^{*}$. Given a biholomorphism
$\rho_{T}$ conjugating $S(T)/\varphi$ and ${\mathbb C}^{*}$ the function
$\psi_{T}^{\varphi} = (1/2 \pi i) \ln \rho_{T}$ is a Fatou coordinate of
$\varphi$ in the fundamental domain $S(T)$.
Shishikura's construction is only valid if $\varphi_{|x=0}$
is of codimension $1$. Since the
Mardesic-Roussarie-Rousseau's system of invariants is obtained by an
improvement of Shishikura's construction then it has the same limitation.

We use Oudkerk's point of view \cite{Oudkerk} based on obtaining transversals
to the dynamics of $\varphi$
by considering trajectories of the real flows of holomorphic vector fields 
$X$ whose exponential ${\rm exp}(X)$ is close to $\varphi$. He obtains Fatou 
coordinates for unfoldings of every $\phi \in \diff{1}{}$ independently 
of the codimension of $\phi$.

Our choice of ${\rm exp}(X)$ is a convergent normal form.
Denote by $\Xnt$ the set of germs of
vector field of the form $f(x,y) \partial/\partial y$ with
$f(0)=(\partial f/\partial y)(0)=0$. Given $\varphi \in \diff{p1}{2}$
there exists $X \in \Xnt$ such that
$y \circ \varphi - y \circ {\rm exp}(X)$ belongs to $(y \circ \varphi -y)^{2}$.
We say that ${\rm exp}(X)$ is a {\it convergent normal form} of $\varphi$
since they are formally conjugated and the infinitesimal generator
$X$ of ${\rm exp}(X)$ is convergent. The transversals are curves of the
form ${\rm exp}(\mu {\mathbb R} X)(x_{0},y_{0})$ for some
$\mu \in {\mathbb S}^{1} \setminus \{-1,1\}$ and
$(x_{0},y_{0}) \in {\mathbb C}^{2}$.
This point of view can be used even if
we do not work with unfoldings and just with discrete deformations
of $\phi \in \diff{1}{} \setminus \{Id\}$ since there exists a universal
theory of unfoldings of germs of vector fields in one variable \cite{Kostov}.

The approach of
Lavaurs-Shishikura-Oudkerk-Mardesic-Roussarie-Rousseau
is of topological type. These constructions imply that the Lavaurs vector
field, i.e.
the unique holomorphic vector field $X_{T}^{\varphi}$ in $S(T)$ such that
$X_{T}^{\varphi}(\psi_{T}^{\varphi}) \equiv 1$, is singular at the fixed
points.
Our approach provides:
\begin{itemize}
 \item Asymptotic developments of $X_{T}^{\varphi}$ until the first non-zero
term in the neighborhood of the fixed points.
\item Accurate estimates for the domains of definition of
${\rm exp}(c X_{T}^{\varphi})$ for $c \in {\mathbb C}$.
\item Canonical normalizing conditions for the Fatou coordinates.
\end{itemize}
These improvements allow us to:
\begin{itemize}
 \item Identify the Taylor's series expansion of the analytic mappings
conjugating $\varphi$, $\zeta \in \diff{p1}{2}$.
\item Study the dependance of the domain of definition of a conjugation
with respect to the parameter.
\item Give a geometrical
interpretation of our complete system of analytic invariants (main theorem).
\end{itemize}
 We use some of the techniques in
\cite{topology} like the dynamical splitting and also others like the
study of polynomial vector fields related to deformations introduced
by Douady-Estrada-Sentenac in \cite{DES}. The polynomial vector fields
that we consider are different. Ours are related
to the infinitesimal properties of the unfolding. They appear after
blow-up transformations. These techniques allow to define fundamental
domains depending on $x$ and representing better the dynamics of
$\varphi$ when $x \to 0$.

%

  Let us remark that the study of germs of diffeomorphism is
useful to classify singular foliations.
For instance consider codimension $1$ complex analytic foliations
defined in a $2$-dimensional manifold. Up to birrational transformation
we can suppose that the singularities are reduced.
Denote by $\Omega_{red} \cn{2}$ the set of germs of 
reduced complex analytic codimension $1$ foliations.
Let $\omega \in \Omega_{red} \cn{2}$;
if the quotient of the eigenvalues $q(\omega)$ is in the domain of
Poincar\ex\ (i.e. $q(\omega) \not \in {\mathbb R}^{-} \cup \{0\}$)
then $\omega$ is conjugated to its linear part.
Anyway, the analytic class of $\omega \in \Omega_{red} \cn{2}$ is determined
by the analytic class of the holonomy of $\omega$ along a ``strong''
integral curve \cite{MaMo:Aen}. Such a holonomy is formally linearizable
if $q(\omega) \in {\mathbb R}^{-} \setminus {\mathbb Q}^{-}$ and
resonant whenever $q(\omega) \in {\mathbb Q}^{-} \cup \{0\}$.
Traditionally a singularity $\omega \in \Omega_{red} \cn{2}$
such that $q(\omega) \in {\mathbb Q}^{-} \setminus \{0\}$ is called resonant 
whereas it is
called a saddle-node if $q(\omega)=0$. The modulus of analytic classification
for both resonant and saddle-node singularities have been described
by Martinet-Ramis \cite{MaRa:aen} \cite{MaRa:ihes}.
Then it is natural to study unfoldings of resonant diffeomorphisms in
order to study unfoldings of resonant singularities and saddle-nodes.
This point of view has been developed by Martinet, Ramis
\cite{Ramis:Tok}, Glutsyuk \cite{Gluglu} and
Mardesic-Roussarie-Rousseau \cite{MRR}. Moreover
Rousseau classifies generic unfoldings of codimension $1$ saddle-nodes
\cite{Rou:sad}. This program can not be carried in higher codimension
without a complete system of analytic invariants for unfoldings of elements
of $\diff{1}{}$ of codimension greater than $1$. We remove such an obstacle
in this paper.

We comment the structure of the paper.
In section \ref{sec:infgen} we introduce the concepts of infinitesimal
generator and convergent normal forms for germs of unipotent diffeomorphism.
We prove that every element of $\diff{p1}{2}$ has a convergent normal form.
Section \ref{sec:onevar} is basically a quick survey about the topological,
formal and analytic classifications of tangent to the identity germs of
diffeomorphism in one variable. We study the formal properties of elements
of $\diff{p1}{2}$ in section \ref{sec:forcon}. We describe the formal
invariants and the structure of the formal centralizer of an
element of $\diff{p1}{2}$. We also reduce the problem of classifying
unfoldings of resonant diffeomorphisms to the tangent to the identity
case via the semisimple-unipotent decomposition.
In section \ref{sec:real} we give a concept of
stability for the real flows of elements of $\Xnt$ and then we describe their
topological behavior in the stable zones.
In section \ref{sec:Fatou} we give a quantitative mesure of how much
$\varphi \in \diff{p1}{2}$ is similar to a convergent normal form.
The estimates are a key ingredient in our refinement of the
Shishikura-Oudkerk-Mardesic-Roussarie-Rousseau's construction. In this way
we obtain Fatou coordinates with controlled asymptotic behavior in the
neighborhood of the fixed points. In section \ref{sec:defanai}
we define the analytic invariants, we describe its nature and compare
with the ones in \cite{MRR}.
Section \ref{sec:trityp} deals with the special case of unfoldings
in which the fixed points set is parameterized by $x$. We can use then
a parameterized version of the Ecalle-Voronin theory.
In section \ref{sec:app} we prove the
main theorem, moreover we provide a complete system of analytic invariants
in both the general and the particular rigid cases.
We prove the optimality of our results in section \ref{sec:otimo}.
\section{Notations and definitions}
Let $\diff{}{n}$ be the group of complex analytic germs of
diffeomorphism at $0 \in {\mathbb C}^{n}$.
Consider coordinates $(x_{1},\hdots,x_{n-1},y) \in {\mathbb C}^{n}$. We
say that $\varphi \in \diff{}{n}$
is a parameterized diffeomorphism if $x_{j} \circ \varphi = x_{j}$ for
all $1 \leq j <n$. We denote by $\diff{p}{n}$ the group of parameterized
diffeomorphisms. Let $\diff{u}{n}$ be the subgroup of $\diff{}{n}$
of unipotent diffeomorphisms, i.e.
$\varphi \in \diff{u}{n}$ if $j^{1} \varphi$ is unipotent. We define
\[ \diff{up}{n} = \diff{u}{n} \cap \diff{p}{n} \]
the group of germs of unipotent parameterized diffeomorphisms.
The formal completions of the previous groups will be denoted with a
hat, for instance $\diffh{}{n}$
is the formal completion of $\diff{}{n}$.

Let $\diff{1}{}$ be the subgroup of $\diff{}{}$ of germs
whose linear part is the identity. We define the set
\[ \diff{p1}{2} = \{ \varphi \in \diff{p}{2} :
\varphi_{|x=0} \in \diff{1}{} \setminus {Id} \} . \]
Then $\diff{p1}{2}$ is the set of one dimensional unfoldings of
one dimensional tangent to the identity germs of diffeomorphism
(excluding the identity).


We define a formal vector field $\hat{X}$ as a derivation of the
maximal ideal of the ring
${\mathbb C}[[x_{1},\hdots,x_{n-1},y]]$. We also express $\hat{X}$ in
the more conventional form
\[ \hat{X} =
\sum_{j=1}^{n-1} \hat{X}(x_{j}) \partial / \partial{x_{j}} +
\hat{X}(y) \partial / \partial{y} . \]
We consider the set $\hat{\mathcal X}_{N} \cn{n}$ of nilpotent formal
vector fields, i.e. the formal vector
fields $\hat{X}$ such that $j^{1} \hat{X}$ is nilpotent. We denote by
${\mathcal X} \cn{n}$ the set of germs
of analytic vector field at $0 \in {\mathbb C}^{n}$.

We denote the rings
${\mathbb C}\{x_{1},\hdots,x_{n-1},y\}$ and
${\mathbb C}[[x_{1},\hdots,x_{n-1},y]]$
by $\vartheta_{n}$ and $\hat{\vartheta}_{n}$ respectively.
We denote $f \sim g$ if $f=O(g)$ and $g=O(f)$.

Let $\varphi \in \diff{up}{n}$.
Denote by $Fix \varphi$ the fixed points set of $\varphi$.
Denote by $Z(\varphi)$ (resp. $\hat{Z}(\varphi)$) the subgroup of
$\diff{p}{n}$ (resp. $\diffh{p}{n}$) whose elements satisfy
$\sigma_{|Fix \varphi} \equiv Id$.
\section{The infinitesimal generator}
\label{sec:infgen}
In this section we associate a formal vector field to every element of
$\diff{u}{n}$.
The properties of this object can be used to provide a complete system
of formal invariants
for the elements of $\diff{up}{n}$ \cite{UPD}. Here, we introduce the
properties that we will use later on.

Let $X \in {\mathcal X} \cn{n}$; suppose that $X$ is singular at $0$.
We denote by ${\rm exp}(tX)$ the flow of the vector field $X$, it is
the unique solution
of the differential equation
\[ \frac{\partial}{\partial{t}} {\rm exp}(tX) = X({\rm exp}(tX)) \]
with initial condition ${\rm exp}(0X)=Id$. We define the exponential
${\rm exp}(X)$ of $X$ as
${\rm exp}(1X)$. We can define the exponential operator for
$\hat{X} \in \hat{\mathcal X}_{N} \cn{n}$.
Moreover the definition coincides with the previous one if $\hat{X}$ is
convergent.
We define
\[
\begin{array}{rccc}
{\rm exp}(\hat{X}): & \hat{\vartheta}_{n} & \to & \hat{\vartheta}_{n}
\\
& g & \to & \sum_{j=0}^{\infty} \frac{\hat{X}^{\circ(j)}}{j!} (g) .
\end{array}
\]
The nilpotent character of $\hat{X}$ implies that the power series
${\rm exp}(\hat{X})(g)$ converges in the Krull
topology for all $g \in \hat{\vartheta}_{n}$. Moreover, since $\hat{X}$
is a derivation
then ${\rm exp}(\hat{X})$ acts like a diffeomorphism, i.e.
${\rm exp}(\hat{X})(g_{1} g_{2})= {\rm exp}(\hat{X})(g_{1})
{\rm exp}(\hat{X})(g_{2})$
for all $g_{1},g_{2} \in \hat{\vartheta}_{n}$.
Moreover $j^{1} {\rm exp}(\hat{X})= {\rm exp}(j^{1} \hat{X})$, thus
$j^{1} {\rm exp}(\hat{X})$ is a unipotent linear
isomorphism. The following proposition is classical.
\begin{pro}
The mapping ${\rm exp}: \hat{\mathcal X}_{N} \cn{n} \to \diffh{u}{n}$
is a bijection.
\end{pro}
Consider the inverse mapping $\log: \diffh{u}{n} \to \hat{\mathcal
H}_{N} \cn{n}$. We can interpret
$\varphi \in \diffh{u}{n}$ as a linear operator $\varphi:
\hat{\mathfrak m} \to \hat{\mathfrak m}$ where
$\hat{\mathfrak m}$ is the maximal ideal of $\hat{\vartheta}_{n}$.
Denote by
$\Theta$ the operator $\varphi - Id$, we have
\[ (\log \varphi)(g) =
\sum_{j=1}^{\infty} {(-1)}^{j+1} j^{-1} \Theta^{\circ (j)} (g) \]
for all $g \in \hat{\vartheta}_{n}$. The power series in the right hand
side converges
in the Krull topology since $\varphi$ is unipotent. Moreover
$j^{1}(\log \varphi)= \log(j^{1} \varphi)$ is
nilpotent and $\log \varphi$ satisfies the Leibnitz rule. We say that
$\log \varphi$ is the {\it infinitesimal generator} of $\varphi$. The
exponential mapping has
a geometrical nature; next proposition claims that $\log \varphi$
preserves the orbits of
$\partial/\partial{y}$ for $\varphi \in \diff{up}{n}$ and also that
$Sing(\log \varphi)=Fix \varphi$.
\begin{pro}
\label{pro:logform}
Let $\varphi \in \diff{up}{n}$. Then $\log \varphi$ is of the form
$\hat{u} (y \circ \varphi -y) \partial/\partial{y}$ for some formal
unit
$\hat{u} \in \hat{\vartheta}_{n}$.
\end{pro}
\begin{proof}
Let $\Theta = \varphi - Id$. We have that $\log \varphi$ is of the form
$\hat{f} \partial / \partial{y}$
since $\Theta(x_{j})=0$ and then $\Theta^{\circ (k)}(x_{j})=0$ for all
$j \in \{1,\hdots,n-1\}$ and all
$k \in {\mathbb N}$. We have $\Theta(y)= y \circ \varphi - y$, moreover
since
\begin{equation}
\label{equ:Taylor}
g \circ \varphi = g + \frac{\partial{g}}{\partial y} (y \circ \varphi
-y) + \sum_{j=2}^{\infty}
\frac{\partial^{j}{g}}{\partial{y}^{j}} \frac{{(y \circ \varphi
-y)}^{j}}{j!}
\end{equation}
we obtain that $\Theta^{\circ (2)}(y) \in (y \circ \varphi -y)
\hat{\mathfrak m}$
where $\hat{\mathfrak m}$ is the maximal ideal of
$\hat{\vartheta}_{n}$.
Again by using the Taylor series
expansion we can prove that $\Theta^{\circ (j)}(y) \in (y \circ \varphi
-y) \hat{\mathfrak m}$
for all $j \geq 2$. Thus  $\log \varphi = (\log \varphi)(y)
\partial/\partial{y}$ is of the form
$\hat{u} (y \circ \varphi -y) \partial/\partial{y}$ for some  $\hat{u}
\in \hat{\vartheta}_{n}$
such that $\hat{u}(0)=1$.
\end{proof}
Let $\varphi = {\rm exp}(\hat{u}(y \circ \varphi -y)
\partial/\partial{y}) \in \diff{up}{n}$. We say that
$\alpha \in \diff{up}{n}$ is a {\it convergent normal form} of
$\varphi$ if
$\log \alpha =  u (y \circ \varphi -y) \partial/\partial{y}$ for some
$u \in {\vartheta}_{n}$
and $y \circ \varphi - y \circ \alpha \in {(y \circ \varphi -y)}^{2}$.
The last condition is equivalent to
$\hat{u} - u \in (y \circ \varphi -y)$.
If $\log \varphi \in {\mathcal X}\cn{n}$ then we say that
$\varphi$ is {\it analytically trivial}.
\begin{pro}
\label{pro:excofor}
Let $\varphi = {\rm exp}(\hat{u} (y \circ \varphi - y) \partial/\partial{y})
\in \diff{up}{n}$. Then $\varphi$ has a convergent normal form.
\end{pro}
\begin{proof}
Let $\Theta = \varphi - Id$. We have
$(\log \varphi)(y) =
\sum_{j=1}^{l} {(-1)}^{j+1} \Theta^{\circ (j)}(y) / j$. Consider the
irreducible decomposition
$f_{1}^{l_{1}} \hdots f_{p}^{l_{p}} g_{1} \hdots g_{q}$
of $y \circ \varphi -y \in {\vartheta}_{n}$
where $l_{j} \geq 2$ for all $j \in \{1, \hdots, p\}$.
Denote $f=y \circ \varphi - y$; we define
$u_{2} =(\ln(1+z)/z) \circ \partial f /\partial y$.
By equation \ref{equ:Taylor} we obtain that
\[ (\log \varphi)(y)/(y \circ \varphi - y) -
\left({ 1 - \frac{\partial{f}/\partial{y}}{2} +
\frac{{(\partial{f}/\partial{y})}^{2}}{3} + \hdots }\right)
\in (f_{1} \hdots f_{p} g_{1} \hdots g_{q}) .\]
We deduce that $\hat{u} - u_{2}$ belongs to $(g_{1} \hdots g_{p})$.

We claim that
$\Theta^{\circ (k)}(y) \in (f_{1}^{l_{1}+k-1} \hdots f_{p}^{l_{p}+k-1})$
for all $k \in {\mathbb N}$.
The result is true for $k=1$ by equation \ref{equ:Taylor}. Since
$f_{j} \circ \varphi - f_{j} \in (f_{j}^{2})$ and
$h \circ \varphi - h \in (y \circ \varphi -y)$
for all $h \in \hat{\vartheta}_{n}$ we deduce that
$\Theta^{\circ (k)}(g) \in (f_{j}^{l_{j}+k-1})$ implies
$\Theta^{\circ (k+1)}(g) \in (f_{j}^{l_{j}+k})$.

Denote $l = \max (l_{1},\hdots,l_{p})$
and
$u_{1}=(\sum_{j=1}^{l} {(-1)}^{j+1} \Theta^{\circ (j)}(y)/j)/f$.
We have that $\hat{u} - u_{1} \in (f_{1}^{l_{1}} \hdots f_{p}^{l_{p}})$.
The function $u_{1}-u_{2}$ belongs to the formal ideal
$(f_{1}^{l_{1}} \hdots f_{p}^{l_{p}} ,g_{1} \hdots g_{q})$; by faithful
flatness there exist $A,B \in \vartheta_{n}$ such that
\[ u_{1} - u_{2} = A f_{1}^{l_{1}} \hdots f_{p}^{l_{p}}
+ B g_{1} \hdots g_{q} .  \]
We define $u = u_{1} - A  f_{1}^{l_{1}} \hdots f_{p}^{l_{p}} =
u_{2} + B  g_{1} \hdots g_{q}$. By construction it is clear that
$\hat{u} - u$ belongs to
$(f_{1}^{l_{1}} \hdots f_{p}^{l_{p}}) \cap (g_{1} \hdots g_{q})$
and then to $(y \circ \varphi -y)$.
\end{proof}
Let $X$ be a holomorphic vector field defined in a connected domain
$U \subset {\mathbb C}$ such that $X \neq 0$. Consider $P \in Sing X$.
There exists a unique meromorphic differential form
$\omega$ in $U$ such that $\omega(X)=0$. We denote by
$Res(X,P)$ the residue of $\omega$ at the point $P$.
Given $Y=f(\underline{x},y) \partial/\partial y$ and
a point $P=(\underline{x}^{0},y^{0}) \in Sing X$ such that
$Sing X$ does not contain $\underline{x}=\underline{x}^{0}$
we define
$Res(X,P)=Res(f(\underline{x}^{0},y) \partial/\partial y,y^{0})$.

Let $\varphi \in \diff{up}{n}$. Consider a convergent normal form
$\alpha$ of $\varphi$. By definition
$Res(\varphi,P)=Res(\log \alpha,P)$ for $P \in Fix \varphi$.
The definition does not depend on the choice of $\alpha$ since given
another convergent normal form $\beta$ of $\varphi$ we have
that $dy/(\log \alpha)(y) - dy/(\log \beta)(y) \in {\vartheta}_{n} dy$.
We denote the function $P \to Res(\varphi, P)$ defined in
$Fix \varphi$ by $Res(\varphi)$.
\section{One variable theory}
\label{sec:onevar}
We introduce here for the sake of completeness some classical results
concerning tangent to
the identity complex analytic germs of diffeomorphism in one variable.
\subsection{Formal theory}
Let $\varphi \in \diff{1}{} = \diff{u}{}$.
We define $\nu(\varphi)$ the order of $\varphi$ as
$\nu(\varphi)= \nu(\varphi(y)-y)-1$.
\begin{pro}
\label{pro:cftg}
Let $\varphi_{1}, \varphi_{2} \in \diff{1}{} \setminus \{ Id \}$. Then
$\varphi_{1}$ is formally
conjugated to $\varphi_{2}$ if and only if
$\nu(\varphi_{1})=\nu(\varphi_{2})$ and
$Res(\varphi_{1})=   Res(\varphi_{2})$. In such a case if $\log
\varphi_{1}$ and $\log \varphi_{2}$
are convergent then $\varphi_{1}$ and $\varphi_{2}$ are analytically
conjugated.
\end{pro}
Supposed that $\varphi_{1}, \varphi_{2}$ are formally conjugated by
$\hat{\sigma} \in \diffh{}{}$.
Then every other formal conjugation can be expressed in the form
$\hat{\tau} \circ \hat{\sigma}$
where $\hat{\tau}$ belongs to the formal centralizer $\hat{Z}(\varphi_{2})$
of $\varphi_{2}$.
As a consequence it is interesting to describe the structure of
$\hat{Z}(\varphi)$ for classification purposes.
\begin{pro}
\label{pro:comonevar}
Let $\varphi \in \diff{1}{} \setminus \{Id\}$. Then there exists
$\hat{\tau}_{0}(\varphi) \in
\diffh{}{}$ satisfying
$(\partial{\hat{\tau}_{0}(\varphi)}/\partial{y})(0)={e}^{2 i \pi/\nu(\varphi)}$
and
$\hat{\tau}_{0}(\varphi)^{\circ (\nu(\varphi))}=Id$ such that
\[ \hat{Z}(\varphi) = \{ \hat{\tau}_{0}(\varphi)^{\circ (r)} \circ {\rm exp}(t
\log \varphi) \
{\rm for} \ r \in {\mathbb Z}/(\nu(\varphi) {\mathbb Z}) \ {\rm and} \
t \in {\mathbb C} \}. \]
Moreover $\hat{Z}(\varphi)$ is a commutative group.
\end{pro}
We say that $\hat{\tau}_{0}(\varphi)$
is the {\it generating symmetry} of $\varphi$.
Let $\kappa_{r}={e}^{2 i  r \pi/\nu(\varphi)}$.
We denote $\hat{\tau}_{0}(\varphi)^{\circ (r)} \circ {\rm
exp}(t \log \varphi)$ by $Z_{\varphi}^{\kappa_{r},t}$. The mapping
$Z_{\varphi}^{\kappa,t} \mapsto (\kappa,t)$ is a bijection from
$\hat{Z}(\varphi)$
to $<{e}^{2 i  \pi/\nu(\varphi)}> \times {\mathbb C}$.
\subsection{Topological behavior}
\label{subsec:topbeh}
Let ${\rm exp}(X)$ be a convergent normal form of
$\varphi \neq Id$ in $\diff{1}{}$. The vector field  $X$ is of the form
$X=(r_{0}{e}^{i \theta_{0}} y^{\nu+1} +
\sum_{j=\nu+2}^{\infty} a_{j} y^{j}) \partial / \partial{y}$
where $\nu = \nu(\varphi)$ and $r_{0} \neq 0$. Consider the blow-up
$\pi: ({\mathbb R}^{+} \cup \{0\}) \times {\mathbb S}^{1} \to {\mathbb R}^{2}$
given by $\pi(r,e^{i \theta}) = r e^{i \theta}$. We denote by
$\tilde{X}$ the strict transform of $Re(X)$, we have
$\tilde{X} = (\pi^{*} Re(X))/r^{\nu}$. We obtain that
\[ \tilde{X} = r \left({
r_{0} Re(e^{i (\nu \theta+ \theta_{0})}) + O(r) }\right)
\frac{\partial}{\partial r} +
\left({ r_{0} Re(-i e^{i (\nu \theta + \theta_{0})}) + O(r)
}\right) \frac{\partial}{\partial \theta}  . \]
We define
$D_{1}(X)=\{ \lambda \in {\mathbb S}^{1} : \lambda^{\nu} e^{i \theta_{0}}
=-1 \}$ and
$D_{-1}(X)=\{ \lambda \in {\mathbb S}^{1} : \lambda^{\nu} e^{i \theta_{0}}
=1 \}$. We have that $\sharp D_{1}(X)= \sharp D_{-1}(X) = \nu$ and
$Sing(\tilde{X}_{|r=0}) = D_{1}(X) \cup D_{-1}(X)$. Moreover, since
\[ \tilde{X}_{|r=0} = (- r_{0} \nu s (\theta-\theta_{1}) +
O({(\theta - \theta_{1})}^{2})) \partial / \partial \theta \]
in the neighborhood of $e^{i \theta_{1}} \in D_{s}(X)$ then
the points in $D_{1}(X)$ are attracting points for $\tilde{X}_{|r=0}$
whereas the points of $D_{-1}(X)$ are repelling.

We define
$\eta = -1/(r_{0} e^{i \theta_{0}} \nu y^{\nu})$, we get
$\tilde{X}(\eta)= r^{- \nu} (1 + O(r))$. Let $\lambda_{1} \in D_{1}(X)$
and consider the set
$S(r_{1},\lambda_{1})= [0 \leq r < r_{1}] \cap
[\lambda \in \lambda_{1} e^{(-i\pi/(4 \nu), i \pi/(4 \nu))}]$. We obtain
$\eta(r,\lambda) \in e^{(-i\pi/4,i\pi/4)}/(\nu r_{0} r^{\nu})$
for all $(r,\lambda) \in S(r_{1},\lambda_{1})$. Since
$\tilde{X}(\eta)=r^{- \nu} (1 + O(r))$ then the points in $S(r_{1},\lambda_{1})$
are attracted to $(0,\lambda_{1})$ by the positive flow
of $\tilde{X}$ for $r_{1}>0$ small enough. Analogously
$(0,\lambda_{1})$ is a repelling point for $\tilde{X}$ if
$\lambda_{1} \in D_{-1}(X)$.

The dynamics of $\varphi$ is a small deformation of the dynamics of
${\rm exp}(X)$. We denote $D_{s}(\varphi)=D_{s}(X)$ for $s \in \{-1,1\}$
and $D(\varphi)=D_{-1}(\varphi) \cup D_{1}(\varphi)$.
These definitions do not depend on the choice of convergent normal form.
Suppose that $\varphi$ and $\varphi^{\circ (-1)}$ are holomorphic in an
small enough open set $U \ni 0$. It is easy to prove that
\[ V_{\varphi}^{\lambda} = \{ P \in  U \setminus \{ 0 \} :
\varphi^{\circ (sn)}(P) \in U \ \forall n \in {\mathbb N} \ {\rm and} \
\lim_{n \to \infty} \varphi^{\circ (sn)}(P) =(0,\lambda) \} \]
is an open set for all $\lambda \in D_{s}(\varphi)$.
A domain $V_{\varphi}^{\lambda}$ for $\lambda \in D_{1}(\varphi)$ is
called an {\it attracting petal}.
A domain $V_{\varphi}^{\lambda}$ for $\lambda \in D_{-1}(\varphi)$ is
called a {\it repelling petal}.

We say that $V(\lambda, \theta)$ is a {\it sector} of direction
$\lambda \in {\mathbb S}^{1}$ and angle
$\theta \in {\mathbb R}^{+}$ if there exists $\mu \in {\mathbb R}^{+}$
such that
$V(\lambda, \theta) = \lambda e^{i[-\theta/2, \theta/2]} (0,\mu]$.
We say that $W(\lambda,\theta)$ is a {\it sectorial domain} of direction
$\lambda \in {\mathbb S}^{1}$ and
angle $\theta \in {\mathbb R}^{+}$ if it contains a sector of direction
$\lambda$ and angle $\theta'$
for all $\theta' \in (0,\theta)$.

The next proposition is a consequence of the previous discussion.
\begin{pro}
Let $\varphi \in \diff{1}{}$. Fix a domain a domain of definition
$0 \in U$. We have
\begin{itemize}
\item $V_{\varphi}^{\lambda}$ is a sectorial domain of direction
$\lambda$ and angle $2\pi/\nu(\varphi)$ for all
$\lambda \in D(\varphi)$.
\item $\{0\} \cup \cup_{\lambda \in D(\varphi)} V_{\varphi}^{\lambda}$
is a neighborhood of $0$.
\item $V_{\varphi}^{\lambda_{0}} \cap V_{\varphi}^{\lambda_{1}} = \emptyset$
if $\lambda_{1} \not \in \{e^{-i \pi/\nu(\varphi)}\lambda_{0},\lambda_{0},
e^{i \pi / \nu(\varphi)}\lambda_{0}  \}$.
\item $V_{\varphi}^{\lambda_{0}} \cap V_{\varphi}^{\lambda_{1}}$ is a
sectorial domain of direction $\lambda_{0} {e}^{i \pi/(2 \nu(\varphi))}$
and angle $\pi/\nu(\varphi)$ for
$\lambda_{1} = e^{i\pi/\nu(\varphi)} \lambda_{0}$.
\end{itemize}
\end{pro}
\subsection{Analytic properties}
\label{subsec:anapro}
Next, we describe the analytic invariants of elements of $\varphi \in \diff{1}{}$.
Choose a convergent normal form $\alpha \in \diff{1}{}$ of $\varphi$.
Consider the equation $(\log \alpha) (\psi_{\alpha})=1$.
A holomorphic solution $\psi_{\alpha}$ is called
a {\it Fatou coordinate} of $\alpha$ or also
an {\it integral of the time form} (or dual form) of
$\log \alpha$. The function $\psi_{\alpha}$ is unique up to an additive constant.
Indeed $\psi_{\alpha}$ is of the form
\[ \psi_{\alpha} = \frac{-1}{\nu(\varphi) a_{\nu(\varphi)+1}}
\frac{1}{{y}^{\nu(\varphi)}}
\left({ 1 + \sum_{j=1}^{\infty} b_{j} y^{j}
}\right) + Res(\varphi) \log y  \]
where $\varphi = y + a_{\nu(\varphi)+1}y^{\nu(\varphi)+1} +
O(y^{\nu(\varphi)+2})$.
Let $\lambda \in  D(\varphi)$; we say that
$\eta \in \vartheta(V_{\varphi}^{\lambda})$
is a Fatou coordinate of $\varphi$ in $V_{\varphi}^{\lambda}$ if
$\eta \circ \varphi = \eta +1$
and $\eta - \psi_{\alpha}$ is bounded.
Clearly the definition does not depend on the choice of $\alpha$.
\begin{pro}
Let $\varphi \in \diff{1}{}$. Consider a convergent normal form
$\alpha$ of $\varphi$ and
a direction $\lambda \in D(\varphi)$. Then there exists a unique
Fatou coordinate
$\psi_{\varphi}^{\lambda}$ of $\varphi$ in $V_{\varphi}^{\lambda}$ such that
$\lim_{y \to 0} (\psi_{\varphi}^{\lambda} - \psi_{\alpha})(y)=0$ in every sector
of direction $\lambda$ and angle
lesser than $2 \pi / \nu(\varphi)$ contained in $V_{\varphi}^{\lambda}$.
Moreover $\psi_{\varphi}^{\lambda}$ is injective.
\end{pro}
We can provide a formula for $\psi_{\varphi}^{\lambda}$.
We define $\Delta = \psi_{\alpha} \circ \varphi - (\psi_{\alpha} + 1)$.
By Taylor's formula we obtain
$\Delta \sim (\varphi(y)-\alpha(y)) \partial \psi_{\alpha} / \partial y$ and then
$\Delta \in {\mathbb C} \{ y \} \cap (y^{\nu(\varphi)+1})$.
Since $(\psi_{\varphi}^{\lambda} - \psi_{\alpha}) -
(\psi_{\varphi}^{\lambda} - \psi_{\alpha}) \circ \varphi = \Delta$
we can obtain $\psi_{\varphi}^{\lambda} - \psi_{\alpha}$ as a telescopic sum.
More precisely
let $\psi_{\alpha}^{\lambda} \in \vartheta(V_{\varphi}^{\lambda})$ be a
Fatou coordinate of $\alpha$. We have
\[ \psi_{\varphi}^{\lambda} = \psi_{\alpha}^{\lambda} +
\sum_{j=0}^{\infty} \Delta \circ \varphi^{\circ (j)}
\ {\rm and} \
\psi_{\varphi}^{\lambda} = \psi_{\alpha}^{\lambda} -
\sum_{j=1}^{\infty} \Delta \circ \varphi^{\circ (-j)} \]
for $\lambda \in D_{1}(\varphi)$ and $\lambda \in D_{-1}(\varphi)$
respectively.

  Let $\varphi \in \diff{1}{}$ with convergent normal form $\alpha$.
Denote $\nu=\nu(\varphi)$.
Consider that  $\psi_{\alpha}^{\lambda} \in \vartheta(V_{\varphi}^{\lambda})$
is chosen for all $\lambda \in D(\varphi)$. We define
\[ \xi_{\varphi}^{\lambda}(z) = \psi_{\varphi}^{\lambda e^{i \pi/ \nu}} \circ
{\left({ \psi_{\varphi}^{\lambda} }\right)}^{\circ (-1)} (z)  \]
for $\lambda \in D(\varphi)$.
The dynamics of $\varphi$ in every $V_{\varphi}^{\lambda}$ is
$z \mapsto z +1$ in the coordinate $\psi_{\varphi}^{\lambda}$. Then
$\xi_{\varphi}^{\lambda}$ is the change of chart which allow to glue
two $z \mapsto z +1$ models corresponding to consecutive petals.
In particular we have
$\xi_{\varphi}^{\lambda} \circ (z+1) \equiv  \xi_{\varphi}^{\lambda}(z)+1$
for all $\lambda \in D(\varphi)$.  Fix $\lambda_{0} \in D(\varphi)$ and
$\psi_{\alpha}^{\lambda_{0}}$. Denote
$\lambda_{j}=\lambda_{0} e^{i \pi j / \nu}$.
There are several possible definitions for $\psi_{\alpha}^{\lambda_{j}}$.
We consider {\it homogeneous coordinates},
supposed $\psi_{\alpha}^{\lambda_{j}}$ is defined we extend it
to
$V_{\varphi}^{\lambda_{j}} \cup V_{\varphi}^{\lambda_{j+1}}$
by analytic continuation. Then we define
$\psi_{\alpha}^{\lambda_{j+1}} = \psi_{\alpha}^{\lambda_{j}} -
\pi i Res(\varphi)/\nu$.
Let us remark that
$\psi_{\varphi}^{\lambda_{0}}=\psi_{\varphi}^{\lambda_{2 \nu}}$.
The definition of $\xi_{\varphi}^{\lambda}$ depends on the choice of
$\psi_{\alpha}^{\lambda_{0}}$.
If we replace $\psi_{\alpha}^{\lambda_{0}}$ with
$\psi_{\alpha}^{\lambda_{0}} + K$
for some $K \in {\mathbb C}$ then $\xi_{\varphi}^{\lambda}$
becomes $(z + K) \circ \xi_{\varphi}^{\lambda} \circ (z-K)$ for all
$\lambda \in D(\varphi)$.
Denote $\zeta_{\varphi} = - \pi i Res(\varphi)/\nu(\varphi)$.
\begin{pro}
Let $\varphi \in \diff{1}{}$ with convergent normal form $\alpha$.
Consider $\lambda \in D_{s}(\varphi)$.
Then there exists $C \in {\mathbb R}^{+}$ such that
\begin{itemize}
\item $\xi_{\varphi}^{\lambda}$ is defined in $s  Imgz < -C$
and $\xi_{\varphi}^{\lambda} \circ (z+1) \equiv (z+1) \circ
\xi_{\varphi}^{\lambda}$.
\item $\lim_{ |Img (z)| \to \infty}
\xi_{\varphi}^{\lambda}(z) - z=\zeta_{\varphi}$.
\item $\xi_{\varphi}^{\lambda} = z + \zeta_{\varphi} +
\sum_{j=1}^{\infty} a_{\lambda, j}^{\varphi} {e}^{-2 \pi i s j z}$
for some  $\sum_{j=1}^{\infty} a_{\lambda, j}^{\varphi} w^{j} \in
{\mathbb C}\{w\}$.
\end{itemize}
\end{pro}
All the possible changes of charts can be realized.
\begin{pro}
\label{pro:reainv}
\cite{V} \cite{mal:ast}
Let $Y \in {\mathcal X}_{N} \cn{}$. Consider
$\sum_{j=1}^{\infty} a_{\lambda, j} w^{j} \in {\mathbb C}\{w\}$
for all $\lambda \in D_{{\rm exp}(Y)}$. There exists
$\varphi \in \diff{1}{}$ with convergent normal form ${\rm exp}(Y)$ such that
$\xi_{\varphi}^{\lambda} =
z + \zeta_{\varphi} + \sum_{j=1}^{\infty} a_{\lambda, j}
{e}^{-2 \pi i s j z}$
for all $\lambda \in D_{s}({\rm exp}(Y))$ and $s \in \{-1,1\}$.
\end{pro}
\subsection{Analytic classification}
\label{subsec:anacla}
Suppose that $\varphi_{1}, \varphi_{2}$ are formally conjugated.
Let $\alpha_{j}$ be a convergent normal form
of $\varphi_{j}$. Then $\alpha_{1}$ and $\alpha_{2}$ are analytically
conjugated by some $h \in \diff{}{}$ by proposition \ref{pro:cftg}.
Up to replace $\varphi_{2}$ with
$h^{\circ (-1)} \circ \varphi_{2} \circ h$ we can suppose that
$\varphi_{1}$ and $\varphi_{2}$ have common normal form
$\alpha_{1}=\alpha_{2}$ and
in particular $\varphi_{1}(y) - \varphi_{2}(y) \in
({y}^{2(\nu(\varphi_{1})+1)})$.
Indeed $\varphi_{1}$ and $\varphi_{2}$ have common convergent normal
form if and only if
$\nu(\varphi_{1})=\nu(\varphi_{2})$ and
$\varphi_{1}(y) - \varphi_{2}(y) \in ({y}^{2(\nu(\varphi_{1})+1)})$.

Let $\varphi_{1},\varphi_{2} \in \diff{1}{}$ with common convergent
normal form $\alpha$.
There exists $\hat{\sigma}(\varphi_{1},\varphi_{2}) \in
\diffh{}{}$ conjugating
$\varphi_{1}$ and $\varphi_{2}$
such that $\hat{\sigma}(\varphi_{1},\varphi_{2})(y) -y \in
(y^{\nu(\varphi)+2})$. Moreover $\hat{\sigma}(\varphi_{1},\varphi_{2})$
is unique. We say that $\hat{\sigma}(\varphi_{1},\varphi_{2})$ is the
{\it privileged formal conjugation} . Choose $\lambda_{0} \in
D(\varphi_{1})=D(\varphi_{2})$ and
$\psi_{\alpha}^{\lambda_{0}}$. The next couple of propositions
are a consequence of Ecalle's theory. We always use homogeneous coordinates.
\begin{pro}
\label{pro:pricon}
Let $\varphi_{1},\varphi_{2} \in  \diff{1}{}$ with common convergent
normal form $\alpha$.
Then for all $\lambda \in D(\varphi_{1})$ there exists a unique holomorphic
$\sigma_{\lambda}: V_{\varphi_{1}}^{\lambda} \to  V_{\varphi_{2}}^{\lambda}$
conjugating $\varphi_{1}$ and $\varphi_{2}$ and such that
$\hat{\sigma}(\varphi_{1},\varphi_{2})$
is a $\nu(\varphi_{1})$-Gevrey asymptotic development of $\sigma_{\lambda}$
in $V_{\varphi_{1}}^{\lambda}$.
Moreover we have
$\sigma_{\lambda}= {(\psi_{\varphi_{2}}^{\lambda})}^{\circ (-1)}
\circ \psi_{\varphi_{1}}^{\lambda}$.
\end{pro}
The expression
$\sigma_{\lambda}: V_{\varphi_{1}}^{\lambda} \to  V_{\varphi_{2}}^{\lambda}$
implies an abuse of notation.
Rigorously $V_{\varphi_{1}}^{\lambda}$ and $V_{\varphi_{2}}^{\lambda}$ can be
replaced by sectorial
domains $W_{\varphi_{1}}^{\lambda}$ and $W_{\varphi_{2}}^{\lambda}$
of direction $\lambda$ and angle $2 \pi/ \nu(\varphi_{1})$
and such that
$\sigma_{\lambda}: W_{\varphi_{1}}^{\lambda} \to  W_{\varphi_{2}}^{\lambda}$
is a biholomorphism. For simplicity we keep this kind of notation
throughout this section.

The elements of the centralizer $\hat{Z}(\varphi)$ of
$\varphi \in \diff{1}{}$ can be realized in the sectorial domains
$V_{\varphi}^{\lambda}$ for every $\lambda \in D_{\varphi}$.
\begin{pro}
\label{pro:comgev}
Let $\varphi  \in  \diff{1}{}$ with convergent normal form $\alpha$.
Consider an element $Z_{\varphi}^{\kappa, t}$ of $\hat{Z}(\varphi)$.
Then for all $\lambda \in D_{\varphi}$ there exists a unique holomorphic
$\tau_{\lambda}: V_{\varphi}^{\lambda} \to  V_{\varphi}^{\lambda \kappa}$
such that $\varphi \circ \tau_{\lambda} = \tau_{\lambda} \circ \varphi$ and
$Z_{\varphi}^{\kappa, t}$ is a
$\nu(\varphi)$-Gevrey asymptotic development of $\tau_{\lambda}$ in
$V_{\varphi}^{\lambda}$.
Moreover we have $\tau_{\lambda}=
{(\psi_{\varphi}^{\lambda \kappa})}^{\circ (-1)} \circ
(\psi_{\varphi}^{\lambda} + t)$.
\end{pro}
We can combine propositions \ref{pro:pricon} and \ref{pro:comgev} to
obtain:
\begin{pro}
\label{pro:gevforcon}
Let $\varphi_{1},\varphi_{2} \in  \diff{1}{}$ with common convergent
normal form $\alpha$.
Consider
$(\kappa,t) \in <{e}^{2 i  \pi/\nu(\varphi_{1})}> \times {\mathbb C}$.
Then for all $\lambda \in D_{\varphi_{1}}$ there exists a unique holomorphic
$\sigma_{\lambda}^{\kappa,t}: V_{\varphi_{1}}^{\lambda} \to
V_{\varphi_{2}}^{\lambda \kappa}$
conjugating $\varphi_{1}$ and $\varphi_{2}$ and such that
$Z_{\varphi_{2}}^{\kappa,t} \circ
\hat{\sigma}(\varphi_{1},\varphi_{2})$
is a $\nu(\varphi_{1})$-Gevrey asymptotic development of
$\sigma_{\lambda}^{\kappa,t}$ in $V_{\varphi_{1}}^{\lambda}$.
Moreover $\sigma_{\lambda}^{\kappa,t}=
{(\psi_{\varphi_{2}}^{\lambda \kappa})}^{\circ (-1)} \circ
(\psi_{\varphi_{1}}^{\lambda} + t)$ in homogeneous coordinates.
\end{pro}
 By uniqueness of the $\nu(\varphi_{1})$-Gevrey sum in sectors of angle
greater than $\pi/\nu(\varphi_{1})$
we deduce that
$Z_{\varphi_{2}}^{\kappa,t} \circ \hat{\sigma}(\varphi_{1},\varphi_{2})$ is
analytic if and only
if $\sigma_{\lambda}^{\kappa,t}=
\sigma_{\lambda {e}^{i \pi/\nu(\varphi_{1})}}^{\kappa,t}$
in $V_{\varphi_{1}}^{\lambda} \cap
V_{\varphi_{1}}^{\lambda {e}^{i \pi/\nu(\varphi_{1})}}$ for all
$\lambda \in D_{\varphi_{1}}$. These conditions can be expressed in
terms of the changes of charts.
\begin{pro}
\label{pro:proca1}
Let $\varphi_{1},\varphi_{2} \in  \diff{1}{}$ with common convergent
normal form.
Then $\varphi_{1}$, $\varphi_{2}$ are analytically conjugated if and
only if there exists
$\kappa \in <{e}^{2 i  \pi/\nu(\varphi_{1})}>$ and $t \in {\mathbb C}$
such that
\begin{equation}
\label{equ:eqca1}
\xi_{\varphi_{2}}^{\lambda \kappa} \circ (z+t) \equiv (z+t) \circ
\xi_{\varphi_{1}}^{\lambda} \ \ \forall \lambda \in D(\varphi_{1}) .
\end{equation}
Indeed the equation \ref{equ:eqca1} is equivalent to
$Z_{\varphi_{2}}^{\kappa,t} \circ \hat{\sigma}(\varphi_{1},\varphi_{2})
\in \diff{}{}$.
\end{pro}
We can find references where
it is claimed that given $\varphi_{1},\varphi_{2} \in \diff{1}{}$
analytically conjugated and with common normal form then
the conjugation can be chosen of the form $y + O(y^{\nu(\varphi_{1})+2})$.
That is equivalent to $\hat{\sigma}(\varphi_{1},\varphi_{2}) \in \diff{}{}$. This
false statement is obtained by neglecting the role of the centralizer
in the analytic conjugation. A clarifying reference can be found in \cite{rey}.
\begin{rem}
Let $\lambda \in D_{s}(\varphi_{1})$. The condition
$\xi_{\varphi_{2}}^{\lambda \kappa} (z+t) =
(z+t) \circ \xi_{\varphi_{1}}^{\lambda}$
is equivalent to
$a_{\lambda \kappa,j}^{\varphi_{2}} e^{-2 \pi i s j t} =
a_{\lambda,j}^{\varphi_{1}}$ for all $j \in {\mathbb N}$.
\end{rem}
\begin{rem}
\label{rem:logcon}
Let $\varphi \in \diff{1}{}$ with convergent normal form $\alpha$.
Then $\log \varphi$ belongs to ${\mathcal X}\cn{}$ if and only if
$\varphi \sim \alpha$ (prop. \ref{pro:cftg}).
Therefore $\log \varphi \in {\mathcal X} \cn{}$ if and only if
$a_{\lambda, j}^{\varphi}=0$ for all $\lambda \in D(\varphi)$ and all
$j \in {\mathbb N}$.
\end{rem}
\section{Formal conjugation}
\label{sec:forcon}
Part of this paper is devoted to explain the relations among
formal conjugations, analytic conjugations and the centralizer
when dealing with elements of $\diff{p1}{2}$. In this section we
study the formal properties of the diffeomorphisms.

\subsection{Formal invariants}
Let $\varphi_{1}, \varphi_{2} \in \diff{p1}{2}$. Suppose that there
exists $\sigma \in \diff{}{2}$
such that $\sigma \circ \varphi_{1} = \varphi_{2} \circ \sigma$. We want
to express $\sigma$ as a composition $\sigma_{1} \circ \sigma_{2}$
such that the action of $\sigma$ on the formal invariants of $\varphi_{1}$
is the same action induced by $\sigma_{2}$. Moreover identifying a possible
$\sigma_{2}$ is much simpler than finding $\sigma$.

The property
$\sigma \circ \varphi_{1}=\varphi_{2} \circ \sigma$ implies that
$\sigma$ conjugates convergent normal forms of $\varphi_{1}$ and $\varphi_{2}$.
We obtain:
\begin{pro}
\label{pro:posfor}
Let $\varphi_{1}, \varphi_{2} \in \diff{p1}{2}$.
Suppose that $\varphi_{1}$ and $\varphi_{2}$ are analytically
conjugated by $\sigma \in \diff{}{2}$.
Then
\begin{itemize}
\item $[(y \circ \varphi_{2} -y) \circ \sigma]/(y \circ \varphi_{1}-y)$
is a unit.
\item $Res(\varphi_{1},P) = Res(\varphi_{2},\sigma(P))$ for all
$P \in Fix \varphi_{1}$.
\end{itemize}
\end{pro}
\begin{rem}
The residue functions are formal invariants \cite{UPD} but for us it is
enough to know that
they are analytic invariants.
\end{rem}
We denote $\tau(Fix \varphi_{1}) = Fix \varphi_{2}$ if
$[(y \circ \varphi_{2} -y) \circ \tau]/(y \circ \varphi_{1}-y)$
is a unit for some $\tau \in \diffh{}{2}$.
In particular $Fix \varphi_{1}=Fix \varphi_{2}$ means that
$Id(Fix \varphi_{1})=Fix \varphi_{2}$.

Consider $\tau \in \diff{}{2}$ holding the two
conditions in prop. \ref{pro:posfor}.
By replacing $\varphi_{2}$ with
$\tau^{\circ (-1)} \circ \varphi_{2} \circ \tau$ we
can suppose $Fix \varphi_{1}=Fix \varphi_{2}$
and $Res(\varphi_{1}) \equiv Res(\varphi_{2})$.
Thus we consider from now on analytic
(resp. formal) conjugations $\sigma$ satisfying
the natural normalizing conditions
$x \circ \sigma \equiv x$ and
$y \circ \sigma -y \in I(Fix \varphi_{1})$  where $I(Fix \varphi_{1})$
is the ideal of $Fix \varphi_{1}$; if such a conjugation
exists we denote $\varphi_{1} \sim \varphi_{2}$
(resp. $\varphi_{1} \stackrel{*}{\sim} \varphi_{2}$).

 We denote ${\mathcal X}_{p1} \cn{2} = \{ X \in {\mathcal X} \cn{2} :
{\rm exp}(X) \in \diff{p1}{2} \}$
In particular ${\rm exp}({\mathcal X}_{p1} \cn{2})$
is the subset of $\diff{p1}{2}$ of convergent normal forms.
\begin{pro}
\label{pro:anconfl}
Let $\alpha_{1},\alpha_{2} \in \diff{p1}{2}$ such that
$\log \alpha_{j} \in \Xnt$ for
$j \in \{1,2\}$. Suppose that $Fix \alpha_{1}=Fix \alpha_{2}$ and
$Res(\alpha_{1}) \equiv Res(\alpha_{2})$.
Then $\alpha_{1} \sim \alpha_{2}$.
\end{pro}
\begin{lem}
\label{lem:eqho}
Let $f \in {\mathbb C}\{x,y\}$ such that $f(0,y) \not \equiv 0$.
Consider $A \in {\mathbb C}\{x,y\}$ such that $(A(x_{0},y)/f(x_{0},y))dy$
has vanishing residues for all $x_{0}$ in a neighborhood of $0$.
Then there exists a germ of meromorphic function $\beta$ such that
$\partial{\beta}/\partial{y}=A/f$ and $\beta f \in \sqrt{f} \subset
{\mathbb C}\{x,y\}$.
\end{lem}
\begin{proof}
Let $P=(0,y_{0}) \neq (0,0)$ be a point close to the origin. Since
$f(P) \neq 0$ there
exists a unique holomorphic solution $\beta$ defined in the
neighborhood of $P$ such that
$\partial{\beta}/\partial{y}=A/f$ and $\beta(x,y_{0}) \equiv 0$.
The residues vanish, then we extend $\beta$ by analytic continuation
to obtain $\beta \in \vartheta(U \setminus (f=0))$ for some neighborhood
$U$ of $(0,0)$.

Consider $Q \in (U \setminus \{(0,0)\} ) \cap (f=0)$. Up to a change of
coordinates
$(x,y + h(x))$ we can suppose that $f=v(x,y) y^{r}$ in the neighborhood
of $Q$
where $y(Q)=0 \neq v(Q)$ and $r \in {\mathbb N}$. The form $(A/f)dy$ is of the
form
$(\sum_{-1 \neq j \geq -r} c_{j}(x) y^{j}) dy$. Then $\beta$ is of the
form $\sum_{0 \neq j \geq -(r-1)} c_{j-1}(x) y^{j} / j +
\beta_{Q}(x)$
for some $\beta_{Q}$ holomorphic in a neighborhood of $Q$.
As a consequence $\beta f$ is holomorphic and vanishes at $f=0$ in a
neighborhood of $Q$.
Hence $\beta f$ belongs to $\vartheta(U \setminus \{(0,0)\})$ and then
to $\vartheta(U)$
since we can remove codimension $2$ singularities.
Clearly we have $\beta f \in I(f=0) = \sqrt{f}$.
\end{proof}
\begin{proof}[Proof of proposition \ref{pro:anconfl}]
There exists $f \in {\mathbb C}\{x,y\}$ such that $\log
\alpha_{j}=u_{j} f \partial/\partial{y}$
for some unit $u_{j} \in {\mathbb C}\{x,y\}$ and all $j \in \{1,2\}$.
  Let us use the path method (see \cite{Rou:ast} and \cite{Mar:ast}).
We define
\[ X_{1+z} = u_{1+z} f \frac{\partial}{\partial{y}} =
\frac{u_{1}u_{2}f}{zu_{1} + (1-z)u_{2}}  \frac{\partial}{\partial{y}}. \]
We have that $X_{1+z} \in \Xnt$ for all
$z \in {\mathbb C} \setminus \{ c \}$ where
$c =u_{2}(0)/(u_{2}(0)-u_{1}(0))$. Moreover
$Sing X_{1+z}$ and $Res(X_{1+z})$ do not depend on $z$.
It is enough to prove $\log \alpha_{1} \sim \log \alpha_{2}$
for $c \not \in [0,1]$. If $c \in [0,1]$ we define
\[ Y^{1}_{1+z} = \frac{u_{1}u_{1+i}f}{zu_{1} + (1-z)u_{1+i}}
\frac{\partial}{\partial{y}}
\ \ {\rm and} \ \
Y^{2}_{1+z} = \frac{u_{1+i}u_{2}f}{zu_{1+i} + (1-z)u_{2}}
\frac{\partial}{\partial{y}}. \]
Since $u_{1+i}(0)/(u_{1+i}(0)-u_{1}(0))$ and
$u_{2}(0)/(u_{2}(0)-u_{1+i}(0))$ do not belong to
$[0,1]$ then we obtain by composition a diffeomorphism conjugating
$\alpha_{1}$ and $\alpha_{2}$.

  Suppose $c \not \in [0,1]$.
We look for $W \in {\mathcal X} \cn{3}$ of the form
$h(x,y,z) f \partial / \partial{y} + \partial / \partial{z}$
such that $[W,X_{1+z}]=0$. We ask $h f$ to be holomorphic in a
connected domain
$V \times V' \subset {\mathbb C}^{2} \times {\mathbb C}$ containing $\{
(0,0) \} \times [0,1]$.
We also require $hf$ to vanish at $(f=0) \times V'$.
Supposed that such a $W$ exists then
${{\rm exp}(W)}_{|z=0}$ is a diffeomorphism conjugating
$\log \alpha_{1}$ and $\log \alpha_{2}$.
The equation $[W,X_{1+z}]=0$ is equivalent to
\[ u_{1+z} f \frac{\partial{(hf)}}{\partial{y}} - hf
\frac{\partial{(u_{1+z}f)}}{\partial{y}}
= \frac{\partial{(u_{1+z} f)}}{\partial{z}} . \]
By simplifying we obtain
\[ u_{1+z} f \frac{\partial{h}}{\partial{y}} - h f
\frac{\partial{u_{1+z}}}{\partial{y}} =
\frac{\partial{u_{1+z}}}{\partial{z}}
\Rightarrow \frac{\partial{(h/u_{1+z})}}{\partial{y}} =
 \frac{1}{u_{1} f} - \frac{1}{u_{2} f}. \]
Let $\beta$ be a solution of
$\partial{\beta}/\partial{y} = 1/(u_{1}f) - 1/(u_{2}f)$ such that
$\beta f \in \sqrt{f}$.
Since $(1/(u_{1}f) - 1/(u_{2}f))dy$ has vanishing residues by
hypothesis then
such a solution exists by lemma \ref{lem:eqho}. We are done by defining
$h=u_{1+z} \beta$.
\end{proof}
Suppose
$Fix \varphi_{1}=Fix \varphi_{2}$ and
$Res(\varphi_{1}) \equiv Res(\varphi_{2})$ for some
$\varphi_{1},\varphi_{2} \in \diff{p1}{2}$.
Proposition \ref{pro:anconfl} implies that up to replace
$\varphi_{2}$ with $\tau^{\circ (-1)} \circ \varphi_{2} \circ \tau$
for some $\tau$ in $\diff{p}{2}$ we can suppose that $\varphi_{1}$
and $\varphi_{2}$ have common convergent normal form.
\begin{pro}
\label{pro:despcon}
Let $\varphi_{1},\varphi_{2} \in  \diff{p1}{2}$ with common convergent
normal form $\alpha$.
Let $f \in {\mathbb C}\{x,y\}$ such that $(y \circ \varphi_{1}-y)/f$ is
a unit and
denote $\hat{u}_{j}= (log \varphi_{j})(y)/f$ for $j \in \{1,2\}$.
Then $\varphi_{1} \stackrel{*}{\sim} \varphi_{2}$ by
\[ \tau(\hat{\beta},\hat{u}_{1},\hat{u}_{2}) \stackrel{def}{=}
{{\rm exp} \left({ \hat{\beta} \frac{\hat{u}_{1}
\hat{u}_{2}}{z\hat{u}_{1} + (1-z)\hat{u}_{2}}
f \frac{\partial}{\partial{y}} + \frac{\partial}{\partial{z}}
}\right)}_{|z=0}  \]
where $\hat{\beta}$ can be any solution of
$\partial{\hat{\beta}}/\partial{y}=1/(\hat{u}_{1}f) - 1/(\hat{u}_{2}
f)$ in ${\mathbb C}[[x,y]]$.
\end{pro}
\begin{proof}
We have that $1/(\hat{u}_{1}f) - 1/(\hat{u}_{2}f) \in {\mathbb
C}[[x,y]]$ since $\varphi_{1}$
and $\varphi_{2}$ have convergent common normal form. Let $\beta_{k}
\in {\mathbb C}\{x,y\}$
such that $\hat{\beta} - \beta_{k} \in {(x,y)}^{k}$. We choose $u_{1,k}
\in {\mathbb C}\{x,y\}$
such that $\hat{u}_{1} - u_{1,k} \in (f){(x,y)}^{k}$; this is possible by
proposition \ref{pro:excofor}.
We define $u_{2,k} \in {\mathbb C}\{x,y\} \setminus (x,y)$ such that
$\partial{\beta_{k}}/\partial{y}=1/(u_{1,k}f) - 1/(u_{2,k}f)$. Now
${\rm exp}(u_{1,k} f \partial/\partial{y})$ and
${\rm exp}(u_{2,k} f \partial/\partial{y})$ are
formally conjugated by $\tau(\beta_{k},u_{1,k},u_{2,k})$
(prop. \ref{pro:anconfl}). We obtain
$u_{j,k} \to \hat{u}_{j}$ and
$\tau(\beta_{k},u_{1,k},u_{2,k}) \to \hat{\tau}$, the limits
considered in the Krull topology. Thus
$\tau(\hat{\beta},\hat{u}_{1},\hat{u}_{2})$
conjugates $\varphi_{1}$ and $\varphi_{2}$.
\end{proof}
\subsection{Formal centralizer}
Let $\varphi \in \diff{p1}{2}$. Next, we study the groups
$\hat{Z}(\varphi)$ and
$\hat{Z}_{up}(\varphi)= \{ \hat{\sigma} \in \diffh{up}{2} :
\hat{\sigma} \circ \varphi = \varphi \circ \hat{\sigma} \}$.
We say that $Fix \varphi$ is of {\it trivial type} if
$I(Fix \varphi)$ is of the form $(f)$ for some $f \in {\mathbb C}\{x,y\}$
such that $(\partial{f}/\partial{y})(0,0) \neq 0$.
\begin{lem}
\label{lem:comtwodif}
Let $\varphi \in \diff{p1}{2}$.  Then
$\hat{Z}_{up}(\varphi)$ is a commutative group given by
\[ \hat{Z}_{up}(\varphi) = \{ {\exp}(\hat{c}(x) \log \varphi) \ {\rm
for \ some } \ \hat{c}(x) \in {\mathbb C}[[x]] \} . \]
Moreover  we have $\hat{Z}(\varphi)= \hat{Z}_{up}(\varphi)$
if $Fix \varphi$ is not of trivial type.
\end{lem}
\begin{proof}
We have that $\hat{\tau} \in \hat{Z}_{up}(\varphi)$ is equivalent
to $[\log \varphi, \log \hat{\tau}]=0$.
%
%
Thus $\log \hat{\tau}$ is of the form
$(\log \hat{\tau})(y) \partial/\partial{y}$ by the same arguments than in
the proof of proposition \ref{pro:logform}. Hence
$[\log \varphi, \log \hat{\tau}]=0$ is equivalent to
$\partial ((\log \varphi)(y) / (\log \hat{\tau})(y)) / \partial y= 0$.
Since $(\log \varphi)(0,y) \not \equiv 0$ then
$\log \hat{\tau} = \hat{c}(x) \log \varphi$ for some
$\hat{c}(x) \in {\mathbb C}[[x]]$.
This  implies
$\hat{Z}_{up}(\varphi) \subset \hat{Z}(\varphi)$, we always
have  $\hat{Z}(\varphi) \subset \hat{Z}_{up}(\varphi)$
in the non-trivial type case.
\end{proof}
We define the order $\nu(\varphi)$ of $\varphi \in \diff{p1}{2}$
as the order of $\varphi_{|x=0} \in \diff{1}{}$.
We define $\nu(X)=\nu({\rm exp}(X))=\nu(X(y)(0,y))-1$ for $X \in \Xnt$.
\begin{lem}
\label{lem:cenntt}
Let $\varphi \in \diff{p1}{2}$. Suppose that $Fix \varphi$ is of trivial
type. Then
\[ \hat{Z}(\varphi) = \{ \hat{\tau}_{0}(\varphi)^{\circ (r)} \circ
{\exp}(\hat{c}(x) \log \varphi) \
{\rm for \ some} \ r \in {\mathbb Z}/(\nu(\varphi) {\mathbb Z})
\ {\rm and} \ \hat{c}(x) \in {\mathbb C}[[x]] \}  \]
where $\hat{\tau}_{0}(\varphi) \in \diffh{p}{2}$ is periodic and
$(\partial (y \circ \hat{\tau}_{0}(\varphi))/\partial{y})(0,0)
={e}^{2 \pi i/\nu(\varphi)}$.
Moreover $\hat{Z}(\varphi)$ is a commutative group.
\end{lem}
We say that $\hat{\tau}_{0}(\varphi)$ is the
{\it generating symmetry} of $\varphi$.
We denote ${\rm exp}(c(x) \log \varphi)$ by $Z_{\varphi}^{1,c}$
whereas we denote
$\hat{\tau}_{0}(\varphi)^{\circ (r)} \circ {\exp}(c(x) \log \varphi)$
by $Z_{\varphi}^{\kappa,c}$ where
$\kappa={e}^{2 \pi i r/\nu(\varphi)}$.
\begin{proof}
Let $\nu=\nu(\varphi_{1})$.
Up to a change of coordinates $(x,h(x,y))$ we can suppose
$Fix \varphi=[y^{\nu+1}=0]$.
By propositions \ref{pro:despcon} and \ref{pro:anconfl} we obtain
$\varphi \stackrel{*}{\sim} {\rm exp}(X)$ where
\[ X = \frac{y^{\nu +1}}{1 + y^{\nu} Res(\varphi,(x,0))}
\frac{\partial}{\partial{y}} . \]
We can suppose $\varphi = {\rm exp}(X)$.
Denote $\tau_{0}=(x,{e}^{2 \pi i/\nu}y)$.
We remark that $\tau_{0}^{*} X =X$.
Given $\hat{\tau} \in \hat{Z}({\rm exp}(X))$ there exists
$r \in {\mathbb Z}$
such that $({\tau}_{0}^{\circ (-r)} \hat{\tau})_{|x=0}$ is tangent
to the identity (prop. \ref{pro:comonevar}). Hence we
get $\hat{\tau} = {\tau}_{0}^{\circ (r)} \circ
{\exp}(\hat{c}(x) X)$ for some $\hat{c}(x) \in {\mathbb C}[[x]]$
by lemma \ref{lem:comtwodif}. Moreover $\hat{Z}({\rm exp}(X))$
is commutative since $\tau_{0}^{*} X =X$.
\end{proof}
Let $X \in \Xnt$. We denote by $Sing_{V} X$ the set of irreducible
components of $Sing X$ which are parameterized by $x$.
Consider $\gamma \in Sing_{V} X$; we denote by $\nu_{X}(\gamma)$
the only element of ${\mathbb N} \cup \{0\}$ such that
$X(y) \in I(\gamma)^{\nu_{X}(\gamma)+1} \setminus
I(\gamma)^{\nu_{X}(\gamma)+2}$.
\begin{pro}
Let $\varphi_{1}, \varphi_{2} \in \diff{p1}{2}$ with common normal
form ${\rm exp}(X)$. Consider $\gamma \in Sing_{V} X$. Then
$\varphi_{1} \stackrel{*}{\sim} \varphi_{2}$ by a unique
$\hat{\sigma}(\varphi_{1}, \varphi_{2},\gamma) \in \diffh{}{2}$
such that
$y \circ \hat{\sigma}(\varphi_{1}, \varphi_{2},\gamma) -y
\in I(\gamma)^{\nu_{X}(\gamma)+2}$.
\end{pro}
By definition the transformation
$\hat{\sigma}(\varphi_{1}, \varphi_{2},\gamma)$
is the {\it privileged formal conjugation} between $\varphi_{1}$
and $\varphi_{2}$ with respect to $\gamma$.
\begin{proof}
There exists a unique solution $\hat{\beta}$ of
$\partial{\hat{\beta}}/\partial{y}=1/(\log \varphi_{1})(y) -
1/(\log \varphi_{2})(y)$ such that $\hat{\beta}_{|\gamma} \equiv 0$.
The formula in proposition \ref{pro:despcon} provides
$\hat{\sigma}(\varphi_{1},\varphi_{2},\gamma)=\hat{\tau}$ conjugating
$\varphi_{1}$ and $\varphi_{2}$ and such that
$y \circ \hat{\sigma}(\varphi_{1}, \varphi_{2},\gamma) -y
\in I(\gamma)^{\nu_{X}(\gamma)+2}$.

Suppose $\hat{\sigma}(\varphi_{1}, \varphi_{2},\gamma)$ is not unique.
Thus we have $y \circ \hat{h} - y \in I(\gamma)^{\nu_{X}(\gamma)+2}$ for some
$\hat{h} \in \hat{Z}_{up}(\varphi_{1}) \setminus \{Id\}$. By lemma
\ref{lem:comtwodif} the transformation $\hat{h}$ is of the form
$Z_{\varphi_{1}}^{1,c}$ for some $c \in {\mathbb C}[[x]]$.
Since $(\log \hat{h})(y)$ belongs to $I(\gamma)^{\nu_{X}(\gamma)+2}$
then $c \equiv 0$ and $\hat{h} \equiv Id$. We obtain a contradiction.
\end{proof}
\subsection{Unfolding of diffeomorphisms $y \to {e}^{2 \pi i p/q}y +
O(y^{2})$}
Consider the sets
\[ \diff{prs}{2} =  \{ \varphi \in \diff{p}{2} :
j^{1} \varphi_{|x=0} \ {\rm is \ periodic} \} \]
and
\[ \diff{pr}{2} = \{ \varphi \in \diff{p}{2} :
j^{1} \varphi_{|x=0} \ {\rm is \ periodic \ but} \
\varphi_{|x=0} \ {\rm is \ not \ periodic} \}. \]
Given $\varphi \in \diff{prs}{2}$ we denote by $q(\varphi)$
the smallest element of ${\mathbb N}$ such that
$(\partial \varphi/\partial y)(0,0)^{q(\varphi)}=1$.
Clearly $\varphi \in \diff{pr}{2}$ implies
$\varphi^{\circ (q(\varphi))} \in \diff{p1}{2}$.
In this paper we classify analytically the elements of $\diff{p1}{2}$.
We obtain for free a complete system of analytic invariants for
the elements of $\diff{pr}{2}$.
In this subsection conjugations are not supposed to be normalized.
\begin{pro}
Let $\varphi_{1},\varphi_{2} \in \diff{prs}{2}$. Then
$\varphi_{1}$, $\varphi_{2}$ are analytically conjugated if and only if
$(\partial \varphi_{1}/\partial y)(0,0) =
(\partial \varphi_{2}/\partial y)(0,0)$
and
$\varphi_{1}^{\circ (q(\varphi_{1}))}, \varphi_{2}^{\circ (q(\varphi_{1}))}$
are analytically conjugated.
\end{pro}
\begin{proof}
The sufficient condition is obvious.
Every $\varphi \in \diff{u}{n}$ admits a unique formal Jordan decomposition
$\varphi = \varphi_{s} \circ \varphi_{u} = \varphi_{u} \circ \varphi_{s}$
in semisimple $\varphi_{s} \in \diffh{}{n}$ and unipotent
$\varphi_{u} \in \diffh{u}{n}$ parts. Semisimple is equivalent to formally
linearizable. The decomposition is compatible with the filtration in the
space of jets, i.e. $j^{k} \varphi = j^{k} \zeta$ implies
$j^{k} \varphi_{s} = j^{k} \zeta_{s}$ and
$j^{k} \varphi_{u} = j^{k} \zeta_{u}$. Moreover we have
$\varphi_{s}, \varphi_{u} \in \diffh{p}{n}$ for all $\varphi \in \diff{p}{n}$.

Denote $q=q(\varphi_{1})$ and
$\upsilon = (\partial \varphi_{1}/\partial y)(0,0)$, we can
suppose $\upsilon \neq 1$.
Suppose $\varphi_{1}^{\circ (q)} \equiv Id$. This implies
$\varphi_{2}^{\circ (q)} \equiv Id$. Denote by $\eta_{k}$ the
unipotent diffeomorphism
$q^{-1} \sum_{j=0}^{q-1}
(x,\upsilon y)^{\circ (-j)} \circ \varphi_{k}^{\circ (j)}$.
By construction $\eta_{k} \circ \varphi_{k} = (x,\upsilon y) \circ \eta_{k}$
for $k \in \{1,2\}$. The diffeomorphism
$\eta_{2}^{\circ (-1)} \circ \eta_{1}$ conjugates $\varphi_{1}$ and
$\varphi_{2}$.

Suppose $\varphi_{1}^{\circ (q)} \not \equiv Id$.
We have that $j^{1} \varphi_{k}$ is conjugated to
$(x,\upsilon y)$ by a linear isomorphism and then semisimple
for $k \in \{1,2\}$. Thus
we obtain $j^{1} \varphi_{k,s} = j^{1} \varphi_{k}$, moreover since
$\varphi_{k,s}$ is formally linearizable then
$\varphi_{k,s}^{\circ (q)} \equiv Id$ for $k \in \{1,2\}$.
We deduce that $\varphi_{k}^{\circ (q)} = \varphi_{k,u}^{\circ (q)}$
for $k \in \{1,2\}$. We obtain  $\log \varphi_{k,u} \not \equiv 0$
for all $k \in \{1,2\}$.

Let $\sigma \in \diff{}{2}$
conjugating $\varphi_{1}^{\circ (q)}$ and $\varphi_{2}^{\circ (q)}$;
it also conjugates
$q \log \varphi_{1,u}$ and $q \log \varphi_{2,u}$ by uniqueness of the
infinitesimal generator and then
$\sigma \circ \varphi_{1,u} = \varphi_{2,u} \circ \sigma$.
Denote $\eta = \sigma^{\circ (-1)} \circ \varphi_{2,s} \circ \sigma$.
We claim that $\varphi_{1,s}=\eta$, this implies that
$\sigma$ conjugates $\varphi_{1}$ and $\varphi_{2}$.
Denote $\rho = \eta^{\circ (-1)} \circ \varphi_{1,s}$.
We have $x \circ \rho \equiv x$ and
$(\partial \rho/\partial y)(0,0) = 1$.
As a consequence $\rho$ is unipotent.
Since both $\eta$ and $\varphi_{1,s}$ commute with $\varphi_{1,u}$
then  $\rho \circ \varphi_{1,u} = \varphi_{1,u} \circ \rho$.
We deduce that
$[\log \rho, \log \varphi_{1,u}] = 0$.
Since $(\log \rho)(x)=0$ then we obtain
$\log \rho = (\hat{c}(x) / x^{m}) \log \varphi_{1,u}$
for some $\hat{c} \in {\mathbb C}[[x]]$ and $m \in {\mathbb Z}_{\geq 0}$.
The equations
$x \circ \eta =x$ and
$\eta_{*} \log \varphi_{1,u} = \log \varphi_{1,u}$ imply that
$\eta$ commutes with $\rho$. This leads us to
$\rho^{\circ (q)} \equiv Id$. In particular
$\hat{c}$ is identically $0$, we obtain $\eta = \varphi_{1,s}$.
\end{proof}
\section{Dynamics of the real flow of a normal form}
\label{sec:real}
Let $\varphi \in \diff{p1}{2}$ with convergent normal form ${\rm exp}(X)$.
Our goal is splitting a domain $|y|<\epsilon$ in several sets in
which the dynamics of $\varphi$ is simpler to analyze.
Afterwards we intend
to analyze the sectors in the parameter space in which $Re(\lambda X)$
($\lambda \in {\mathbb S}^{1} \setminus \{-1,1\}$) has a stable
behavior. The stability will provide well-behaved transversals to $Re(X)$.
Such transversals are the base to construct Fatou coordinates
of $\varphi$ for all $x$ in a neighborhood of $0$.
%
%

Consider the function
\[
\begin{array}{cccc}
ag_{X}^{\epsilon}: & B(0,\delta) \times \partial B(0,\epsilon) &
\to & {\mathbb S}^{1} \\
& (x,y) & \mapsto & (X(y)/y)/|X(y)/y| .
\end{array}
\]
By lifting $ag_{X}^{\epsilon}$ to
${\mathbb R} = \widetilde{{\mathbb S}^{1}}$ we obtain a mapping
${arg}_{X}^{\epsilon}:
B(0,\delta) \times {\mathbb R} \to {\mathbb R}$ such that
$e^{2 \pi i \theta} \circ {arg}_{X}^{\epsilon}(x,\theta)
= ag_{X}^{\epsilon}(x,\epsilon {e}^{2 \pi i \theta})$.
It is easy to prove that
$(\partial{{arg}_{X}^{\epsilon}}/\partial{\theta})(0,\theta)$
tends uniformly to $\nu(X)$ when $\epsilon \to 0$. By continuity we
obtain that
$\partial{{arg}_{X}^{\epsilon}}/\partial{\theta}$ is very
close to $\nu(X)$ for $0 < \epsilon <<1$ and
$0 < \delta(\epsilon) << 1$.

Let $X \in \Xnt$ and fix $0 < \epsilon <<1$. We define the set
$T_{X}^{\epsilon}(x_{0})$ of tangent points between
$Re(X)_{|x=x_{0}}$ and $\partial B(0,\epsilon)$ for
$x_{0} \in B(0,\delta(\epsilon))$.
Denote
the set $\cup_{x \in B(0,\delta)} \{x\} \times T_{X}^{\epsilon}(x)$
by $T_{X}^{\epsilon}$. We say that a point $y_{0} \in T_{X}^{\epsilon}(x_{0})$
is {\it convex} if the germ of trajectory of
$Re(X)_{|x=x_{0}}$ passing through $y_{0}$
is contained in $\overline{B}(0,\epsilon)$. Next lemma is a consequence of
$\partial{{arg}_{X}^{\epsilon}}/\partial{\theta} \sim \nu(X)$
and $T_{\lambda X}^{\epsilon}(x_{0}) =
ag_{X}^{\epsilon}(x_{0},y)^{\circ (-1)}\{-i/\lambda ,i /\lambda\}$.
\begin{lem}
\label{lem:tgpt}
Let $X \in \Xnt$. There exist $\epsilon_{0}>0$ and
$\delta_{0}:(0,\epsilon_{0}) \to {\mathbb R}^{+}$ such that
$T_{\lambda X}^{\epsilon}(x_{0})$ is composed of $2 \nu(X)$ convex points
for all $\lambda \in {\mathbb S}^{1}$,
$0 < \epsilon < \epsilon_{0}$ and $x_{0} \in B(0,\delta_{0}(\epsilon))$.
Moreover, each connected component of
$\partial{B(0,\epsilon)} \setminus T_{\lambda X}^{\epsilon}(x_{0})$
contains a unique point of $T_{\mu X}^{\epsilon}(x_{0})$ for all
$\mu \in {\mathbb S}^{1} \setminus \{ -\lambda, \lambda \}$.
\end{lem}
\begin{rem}
Fix $\lambda \in {\mathbb S}^{1}$. We have
$ T_{\lambda X}^{\epsilon}(x) =
\{ T_{\lambda X}^{\epsilon, 1}(x), \hdots,
T_{\lambda X}^{\epsilon, 2 \nu(X)}(x)   \}$
for all $x \in B(0,\delta_{0}(\epsilon))$ where
$T_{\lambda X}^{\epsilon, j}: B(0,\delta_{0}(\epsilon)) \to
T_{X}^{\epsilon}$ is continuous for all $1 \leq j \leq 2\nu(X)$.
\end{rem}
\subsection{Splitting the dynamics}
For simplicity we consider
${\mathcal X}_{tp1} \cn{2} \subset {\mathcal X}_{p1} \cn{2}$
and $\diff{tp1}{2} \subset \diff{p1}{2}$
whose elements satisfy that their singular or fixed points sets respectively
are union of smooth curves transversal to $\partial/\partial{y}$.
For all $\varphi \in \diff{p1}{2}$ there exists
$k \in {\mathbb N}$ such that
$(x^{1/k},y) \circ \varphi \circ (x^{k},y) \in \diff{tp1}{2}$.

Let $X \in \Xt$. We define $T_{0}=(|y| \leq \epsilon)$.
Suppose that we have a sequence $\beta=\beta_{0}\hdots\beta_{k}$ where
$\beta \in \{0\} \times {\mathbb C}^{k}$ and $k \geq 0$ and a set
$T_{\beta} = (|t| \leq \eta)$ in coordinates
$(x,t)$ canonically associated to $T_{\beta}$. The coordinates
$(x,y)$ are canonically associated to $T_{0}$. Suppose also that
\[ X = x^{d_{\beta}}
v(x,t) (t-\gamma_{1}(x))^{s_{1}} \hdots (t-\gamma_{p}(x))^{s_{p}}
\partial / \partial{t} \]
where $\gamma_{1}(0)=\hdots=\gamma_{p}(0)=0$
and $(v=0) \cap T_{\beta}=\emptyset$.
Denote $\nu(\beta)=s_{1}+\hdots+s_{p}-1$ and $N(\beta)=p$.
Define
$X_{\beta, E}=(X(t)/x^{d_{\beta}})\partial /\partial{t}$.
Denote by $TE_{\mu X}^{\beta,\eta}(r, \lambda)$ the set of
tangent points between
$Re(\lambda^{d_{\beta}} \mu X_{\beta, E})_{|x=r \lambda}$ and $|t|=\eta$
for $(r,\lambda,\mu) \in {\mathbb R}_{\geq 0} \times {\mathbb S}^{1}
\times {\mathbb S}^{1}$.
If $N(\beta)=1$ then we define $E_{\beta}=T_{\beta}$,
we do not split $T_{\beta}$. We denote $\dot{E}_{\beta} = [|t| < \eta]$.

Suppose $N(\beta)>1$. Denote
$S_{\beta}=\{(\partial \gamma_{1}/\partial x)(0),\hdots,
(\partial \gamma_{p}/\partial x)(0)\}$.
We define $t=xw$ and the sets
$E_{\beta} = T_{\beta} \cap [|t| \geq |x|\rho]$
and $M_{\beta}= (|w| \leq \rho)$
for some $\rho>>0$.
We denote $\dot{E}_{\beta} =[\rho |x|  < |t| < \eta]$. We have
\[ X = x^{d_{\beta}+s_{1}+\hdots+s_{p}-1}v(x,xw)
{\left({ w -  \gamma_{1}(x) / x }\right)}^{s_{1}} \hdots
{\left({ w - \gamma_{p}(x) / x }\right)}^{s_{p}}
\partial / \partial{w}  \]
in $M_{\beta}$, we define $m_{\beta}=d_{\beta}+\nu(\beta)$
and the polynomial vector field
\[ X_{\beta}(\lambda) = \lambda^{m_{\beta}} v(0,0)
{\left({ w - (\partial \gamma_{1} / \partial x)(0) }\right)}^{s_{1}} \hdots
{\left({ w - (\partial \gamma_{p} / \partial x)(0) }\right)}^{s_{p}}
\partial / \partial{w}  \]
for $\lambda \in {\mathbb S}^{1}$.
We define $I_{\beta}=(|w| \leq \rho) \setminus
\cup_{\zeta \in S_{\beta}} (|w - \zeta| < r(\zeta))$
where $r(\zeta)>0$ is small enough for all $\zeta \in S_{\beta}$.
We define $X_{\beta,M}=(X(w)/x^{m_{\beta}})\partial/\partial{w}$;
we denote by $TI_{\mu X}^{\beta,\rho}(r, \lambda)$ the set of
tangent points between
$Re(\lambda^{m_{\beta}} \mu X_{\beta,M})_{|x=r \lambda}$ and $|w|=\rho$.
Finally we define
$\dot{I}_{\beta}=(|w| < \rho) \setminus
\cup_{\zeta \in S_{\beta}} (|w - \zeta| \leq r(\zeta))$.

Fix $\zeta \in S_{\beta}$. We define $d_{\beta \zeta}=m_{\beta}$.
Consider the coordinate $t' = w - \zeta$.
We denote $T_{\beta \zeta}= (|t'| \leq r(\zeta))$.
We have
\[ X = x^{d_{\beta \zeta}} h(x,t')
\prod_{(\partial \gamma_{j}/\partial x)(0) = \zeta}
{\left({ t' - \left({
\gamma_{j}(x) / x - \zeta }\right)
}\right)}^{s_{j}}  \partial / \partial{w} . \]

Every set $M_{\beta}$ with $\beta \neq \emptyset$ is called
a {\it magnifying glass set}. The sets $E_{\beta}$ are called
{\it exterior sets} whereas the sets $I_{\beta}$ are called
{\it intermediate sets}.

In the previous paragraph we introduced
a method to divide $|y| \leq \epsilon$ in a union of exterior
and intermediate sets.

\underline{Example}: Consider $X = y(y-{x}^{2})(y-x) \partial/\partial{y}$.
We have
\[ (|y| \leq \epsilon) = E_{0} \cup I_{0} \cup E_{01} \cup E_{00} \cup I_{00}
\cup E_{000} \cup E_{001} . \]
We have
$X_{0}(1)=w_{1}^{2}(w_{1}-1)\partial/\partial{w_{1}}$ and
$X_{00}(1)=-w_{2}(w_{2}-1) \partial/\partial{w_{2}}$
where $y=xw_{1}$ and $y=x^{2}w_{2}$. We also get $m_{0}=2$ and $m_{00}=3$.
\begin{lem}
\label{lem:tgpt2}
Let $X \in \Xt$ and an exterior set
$E_{\beta}=[\eta \geq |t| \geq \rho|x|]$ associated to $X$
with $0 < \eta <<1$. Then
$TE_{\mu X}^{\beta, \eta}(r, \lambda)$ is composed of
$2 \nu(\beta)$ convex points for all
$(\lambda,\mu) \in {\mathbb S}^{1} \times {\mathbb S}^{1}$
and $r$ close to $0$. Each connected component of
$[|t|=\eta] \setminus TE_{\mu X}^{\beta, \eta}(r, \lambda)$
contains a unique point of
$TE_{\mu' X}^{\beta, \eta}(r, \lambda) \ \forall
\mu' \in {\mathbb S}^{1} \setminus \{ -\mu, \mu \}$.
\end{lem}
\begin{lem}
\label{lem:tgpt3}
Let $X \in \Xt$ and a magnifying glass set
$M_{\beta}=[|w| \leq \rho]$ associated to $X$ with $\rho >>0$.
Then $TI_{\mu X}^{\beta, \rho}(r, \lambda)$ is composed of
$2 \nu(\beta)$ convex points for all
$(\lambda,\mu) \in {\mathbb S}^{1} \times {\mathbb S}^{1}$
and $r$ close to $0$. Moreover each connected component of
$[|w|=\rho] \setminus TI_{\mu X}^{\beta, \rho}(r, \lambda)$
contains a unique point of
$TI_{\mu' X}^{\beta, \rho}(r, \lambda)$ for all
$\mu' \in {\mathbb S}^{1} \setminus \{ -\mu, \mu \}$.
\end{lem}
Lemma \ref{lem:tgpt2} is the analogue of lemma \ref{lem:tgpt}
for exterior sets.
Lemma \ref{lem:tgpt3} is deduced from the polynomial character
of $X_{\beta}(1)$ since
$\partial {\arg}_{X_{\beta}(1)}^{\rho}/\partial \theta \sim \nu(\beta)$
when $\rho \to \infty$.

Let $X \in {\mathcal X} \cn{n}$. Consider a set
$F \subset {\mathbb C}^{n}$ contained in the domain of definition of $X$.
Denote by $\dot{F}$ the interior of $F$.
We define $It(X, P,F)$ the maximal interval where
${\rm exp}(zX)(P)$ is well-defined and belongs to
$F$ for all $z \in It(X, P,F)$ whereas
${\rm exp}(zX)(P)$ belongs to $\dot{F}$ for
all $z \neq 0$ in the interior of $It(X,P,F)$. We define
\[ \partial It(X,P,F) = \{ \inf(It(X,P,F)), \sup(It(X,P,F)) \}
\subset {\mathbb R} \cup \{-\infty, \infty\}. \]
We denote $\Gamma(X,P,F)={\rm exp}(It(X,P,F)X)(P)$.

We will consider coordinates $(x,y) \in {\mathbb C} \times {\mathbb C}$ or
$(r,\lambda,y) \in
{\mathbb R}_{\geq 0} \times {\mathbb S}^{1} \times {\mathbb C}$
in ${\mathbb C}^{2}$. Given a set $F \subset {\mathbb C}^{2}$ we denote by
$F(x_{0})$ the set $F \cap [x=x_{0}]$ and by $F(r_{0},\lambda_{0})$
the set $F \cap [(r,\lambda)=(r_{0},\lambda_{0})]$. In the next subsections
we analyze the dynamics in the exterior and intermediate sets.
\subsection{Parabolic exterior sets}
Let $X \in \Xt$. Suppose we have
\[ X = x^{d_{\beta}}
v(x,t) (t-\gamma_{1}(x))^{s_{1}} \hdots (t-\gamma_{p}(x))^{s_{p}}
\partial / \partial{t} \]
in some exterior set $E_{\beta}=[\eta \geq |t| \geq |x|\rho]$ for some
$\rho \geq 0$. We say that $E_{\beta}$ is parabolic
if $s_{1}+\hdots+s_{p} \geq 2$. In particular $E_{0}$ is always parabolic
since $\nu(X) \geq 1$.
\begin{lem}
\label{lem:bhvtg}
Let $X \in \Xt$ and a parabolic exterior set
$E_{\beta} = [|t|] \leq \eta]$ associated to $X$ with $0<\eta<<1$.
Consider $\mu \in {\mathbb S}^{1}$ and
$t_{0} \in TE_{\mu X}^{\beta, \eta}(r, \lambda)$. Then we have
$It(\mu \lambda^{d_{\beta}} X_{\beta,E}, (r, \lambda,t_{0}),
[|t| \leq \eta])={\mathbb R}$ and
$\lim_{z \in {\mathbb R}, |z| \to \infty}
{\rm exp}(z \mu \lambda^{d_{\beta}} X_{\beta, E})(r \lambda ,t_{0})$
is the point in
$E_{\beta}(r,\lambda) \cap Sing X_{\beta, E}$.
\end{lem}
\begin{proof}
Consider $\eta_{0}>0$  such that
$TE_{\mu X}^{\beta, \eta}(r, \lambda)$
is composed of $2 \nu(\beta)$ convex points for all
$0< \eta < \eta_{0}$,
$(r,\lambda) \in [0,\delta(\eta)) \times {\mathbb S}^{1}$
and $\mu \in {\mathbb S}^{1}$.

Fix $0 < \eta < \eta_{0}$, $\mu \in {\mathbb S}^{1}$
and $(r,\lambda) \in [0,\delta(\eta)) \times {\mathbb S}^{1}$.
Denote $Y=(\mu \lambda^{d_{\beta}} X_{\beta, E})_{x=r \lambda}$.
We have that $Sing Y$ is a point $t=\gamma_{0}$.
Let $\tilde{Y}$ be the strict transform of $Re(Y)$ by the
blow-up $\pi:({\mathbb R}^{+} \cup \{0\}) \times {\mathbb S}^{1}
\to {\mathbb C}$ of
$t=\gamma_{0}$ given by $\pi(s,\gamma)=s \gamma + \gamma_{0}$.
We consider the set
$S \subset \pi^{-1}[|t| < \eta] \setminus Sing \tilde{Y}$
of points $(s,\gamma)$  such that
$It(\tilde{Y},(s,\gamma),[|t| < \eta])={\mathbb R}$ and
$\lim_{z \to \pm \infty}  {\rm exp}(z \tilde{Y})
(s,\gamma) \in Sing \tilde{Y}$.
By the discussion in subsection \ref{subsec:topbeh}
the set $S$ has exactly $2 \nu(\beta)$ connected components. More precisely
every connected component of $S$ contains exactly an arc
$\{0\} \times e^{(i\theta,i(\theta+\pi/\nu(\beta)))}$ for
$e^{i \theta} \in Sing \tilde{Y}$.

Consider a connected component $C$ of $S$.
We have $It(\tilde{Y},t_{0},[|t| < \eta + c] )={\mathbb R}$ for all
$t_{0} \in \partial C$ and $c \in {\mathbb R}^{+}$.
By Poincar\ex -Bendixon's theorem
the $\alpha$ and $\omega$ limits of $t_{0}$ by $Re(Y)$ are
either $\gamma_{0}$ or a cycle enclosing $\gamma_{0}$
since the points in $Sing \tilde{Y}$ are either attracting or repelling.
The second possibility is excluded by Cartan's lemma. We deduce that there
exists $t_{C} \in [|t| = \eta]  \cap \partial C$.
Clearly $t_{C} \in \overline{C}$ implies that
$t_{C} \in TE_{\mu X}^{\beta,\eta}(r, \lambda)$
and that  ${\rm exp}(zY)(t_{C})$ belongs to
$[|t| \leq \eta] $ for all $z \in {\mathbb R}$. Moreover we obtain
$\lim_{|z| \to \infty} {\rm exp}(z Y)(t_{C})= (x_{0},\gamma_{0})$.
The number of connected components of $S$ coincides with
$\sharp TE_{\mu X}^{\beta, \eta}(r, \lambda)$.
We deduce that ${\rm exp}(zY)(t_{C}) \in [|t| < \eta]$  for all
$z \in {\mathbb R} \setminus \{0\}$
since $\overline{C}_{1} \cap \overline{C}_{2} = \emptyset$
for different connected components of $S$.
\end{proof}
\begin{pro}
\label{pro:bhvtg}
Let $X \in \Xt$ and let
$E_{\beta} = [\eta \geq |t| \geq \rho |x|]$ be a
parabolic exterior set associated to $X$.
Consider $t_{0} \in TE_{\mu X}^{\beta, \eta}(r, \lambda)$
and $\mu \in {\mathbb S}^{1}$. Then we have
\[ \lim_{z \in {\mathbb R}, z \to c}
{\rm exp}(z \mu \lambda^{d_{\beta}} X_{\beta, E})(r, \lambda,t_{0})
\in (\partial{E_{\beta}} \cup Sing X_{\beta, E}) \setminus [|t|=\eta] \]
for $c \in \partial
It(\mu \lambda^{d_{\beta}} X_{\beta,E},(r, \lambda,t_{0}),E_{\beta})$.
\end{pro}
\begin{proof}
If $N(\beta)=1$ the result is true
by lemma \ref{lem:bhvtg}. Suppose $N(\beta)>1$.
Consider $\eta_{0}>0$  and $\rho_{0}>0$ such that
$TE_{\mu X}^{\beta, \eta}(r, \lambda)$ and
$TI_{\mu X}^{\beta,\rho}(r, \lambda)$
are both composed of $2 \nu(\beta)$ convex points for all
$0< \eta < \eta_{0}$, $\rho > \rho_{0}$,
$(r,\lambda) \in [0,\delta(\eta,\rho)) \times {\mathbb S}^{1}$
and $\mu \in {\mathbb S}^{1}$.

Fix $0<\eta<\eta_{0}$ and $\rho>\rho_{0}$.
We can suppose that $r \lambda \neq 0$ since otherwise
the proof is analogous to the proof  in lemma \ref{lem:bhvtg}.

Fix $(r,\lambda) \in (0,\delta(\eta,\rho)) \times {\mathbb S}^{1}$
and $\mu \in {\mathbb S}^{1}$.
Consider a point $t_{1} \in TI_{\mu X}^{\beta, \rho}(r, \lambda)$.
There exists exactly one connected component
$H_{s}$ of
$[|w| = \rho] \setminus  TI_{\mu X}^{\beta, \rho}(r, \lambda)$
such that $t_{1} \in \overline{H_{s}}$ and
$Re(s \mu X)$ points towards $|w| < \rho$ for
$s \in \{-1,1\}$.
We define $S(t_{1})$ as the set of points $t$ in
$\dot{E}_{\beta}(r \lambda)$ such that
there exists $c_{-1}(t), c_{1}(t) \in {\mathbb R}^{+}$ satisfying
that
${\rm exp}((-c_{-1},c_{1})\mu X)(r \lambda,t)$
is well-defined and contained in $\dot{E}_{\beta}$
whereas
${\rm exp}(s c_{s} \mu X)(r \lambda,t) \in H_{s}$
for $s \in \{-1,1\}$. Clearly $S(t_{1}) \neq \emptyset$ since
$t_{1} \in \overline{S(t_{1})}$.

Like in lemma \ref{lem:bhvtg}
there exists a unique $t_{0} \in \overline{S(t_{1})} \cap
TE_{\mu X}^{\beta, \eta}(r, \lambda)$.
We deduce that
$It=It(\mu X,(r \lambda,t_{0}),E_{\beta})$
is compact. Moreover we have
\[ {\rm exp}(h_{I} \mu \lambda^{d_{\beta}} X_{\beta, E})
(r \lambda,t_{0}) \in H_{-1} \ {\rm and} \
 {\rm exp}(h_{S} \mu \lambda^{d_{\beta}} X_{\beta, E})
(r \lambda,t_{0}) \in H_{1} \]
where $It=[h_{I}/r^{d_{\beta}},h_{S}/r^{d_{\beta}}]$.
Since $\sharp TE_{\mu X}^{\beta, \eta}(r, \lambda) =
\sharp TI_{\mu X}^{\beta, \rho}(r, \lambda)$ we are done.
\end{proof}
Let $X \in \Xt$.
We define $SC_{\mu X}^{\beta,\eta}(r, \lambda)$ the set
of connected components of
\[ (\dot{E}_{\beta}(r, \lambda) \setminus Sing X_{\beta, E}) \setminus
\cup_{t \in TE_{\mu X}^{\beta, \eta}(r, \lambda)}
\Gamma(\mu \lambda^{d_{\beta}} X_{\beta,E},(r, \lambda,t),E_{\beta}).  \]

The behavior of the trajectories passing through tangent points
characterizes the dynamics of $Re(\mu X)$ in a parabolic exterior set.
It is a topological product. The next results are a consequence of this fact.
\begin{pro}
\label{pro:bhvtg2}
Let $X \in \Xt$ and let
$E_{\beta} = [\eta \geq |t| \geq \rho |x|]$ be a parabolic exterior set
associated to $X$.
Consider $t_{0} \in \dot{E}_{\beta}(r,\lambda)$ and $\mu \in {\mathbb S}^{1}$.
Then we have
\[ \lim_{z \in {\mathbb R}, z \to c}
{\rm exp}(z \mu \lambda^{d_{\beta}} X_{\beta, E})(r, \lambda,t_{0})
\in \partial{E_{\beta}} \cup Sing X_{\beta, E} \]
for $c \in \partial
It(\mu \lambda^{d_{\beta}} X_{\beta,E},(r, \lambda,t_{0}),E_{\beta})$.
\end{pro}
\begin{proof}
Let $C \in SC_{\mu X}^{\beta,\eta}(r, \lambda)$.
Consider the set $L_{C}$ of points in $C$ satisfying the result in
the proposition. It is enough to prove that $C=L_{C}$ for all
$C \in SC_{\mu X}^{\beta,\eta}(r, \lambda)$.

The points in $C$ in the neighborhood of points in
$TE_{\mu X}^{\beta,\eta}(r, \lambda)$ are contained in $L_{C}$ by
proposition \ref{pro:bhvtg} and continuity of the flow.
We have that $C$ is a simply connected open set such that
$C \cap Sing X_{\beta, E}=\emptyset$.
Moreover every trajectory of $Re(\mu \lambda^{d_{\beta}} X_{\beta, E})$
contained in
$E_{\beta}$ and intersecting the set
$TE_{\mu X}^{\beta,\eta}(r, \lambda) \cup TI_{\mu X}^{\beta, \rho}(r, \lambda)$
is disjoint from $C$. Thus the set $L_{C}$ is open and closed in
$C$ and then $L_{C}=C$.
\end{proof}
The next result can be proved like proposition \ref{pro:bhvtg2}, it is true
in the neighborhood of the tangent points by lemma \ref{lem:bhvtg} and
it defines an open and closed
property in connected sets. We skip the proof.
\begin{cor}
\label{cor:npee}
Let $X \in \Xt$ and let $E_{\beta} = [\eta \geq |t| \geq \rho |x|]$ be
a parabolic exterior set associated to $X$. Let
$(\mu_{0},r, \lambda,t_{0}) \in {\mathbb S}^{1} \times [0,\delta)
\times {\mathbb S}^{1}
\times \partial B(0,\eta)$ such that
$Re(\mu_{0} \lambda^{d_{\beta}} X_{\beta, E})(r \lambda,t_{0})$ does not
point towards ${\mathbb C} \setminus \overline{B}(0,\eta)$.
Then we obtain
\[ \lim_{z \to c(\mu_{0},r, \lambda,t_{0})} {\rm exp}(z \mu_{0} \lambda^{d_{\beta}} X_{\beta, E})
(r \lambda,t_{0}) \in
(\partial{E_{\beta}} \cup Sing X_{\beta, E}) \setminus [|t|=\eta] \]
for
$c(\mu_{0},r, \lambda,t_{0}) = \sup
It(\mu_{0} \lambda^{d_{\beta}} X_{\beta,E},(r,\lambda,t_{0}), {E}_{\beta})
\in {\mathbb R}^{+} \cup \{\infty\}$.
%
%
\end{cor}
Let $X=x^{d_{\beta}} v(x,t) \prod_{j=1}^{N(\beta)}
{(t-\gamma_{j}(x))}^{s_{j}} \partial/\partial{t} \in \Xt$.
We define
\[ X_{\beta}^{0}= v(0,t-\gamma_{1}(x))
{(t-\gamma_{1}(x))}^{\nu(\beta)+1} \partial / \partial{t} . \]
Let $\psi_{\beta,E}^{0}$ be a holomorphic integral of the time form of
$X_{\beta}^{0}$ in the neighborhood of $E_{\beta} \setminus Sing X_{\beta, E}$.
We have
$\psi_{\beta, E}^{0}(x,e^{2 \pi i}y) -\psi_{\beta, E}^{0}(x,y) \equiv
2 \pi i Res(X_{\beta}^{0},(0,0))$, in general $\psi_{\beta, E}^{0}$ is
multivaluated.
Consider  a holomorphic integral $\psi_{\beta, E}$ of the time form of
$X_{\beta, E}$
in the neighborhood of $E_{\beta} \setminus Sing X$
such that $\psi_{\beta, E}(0,y) \equiv \psi_{\beta, E}^{0}(0,y)$.
Clearly $\psi_{X, \beta}^{0} =\psi_{\beta, E}^{0}/x^{d_{\beta}}$ and
$\psi_{X,\beta}=\psi_{\beta, E}/x^{d_{\beta}}$ are integrals of the time forms
of $x^{d_{\beta}} X_{\beta}^{0}$ and $X$ respectively. We want to
provide accurate estimates for $\psi_{X,\beta}$.
\begin{lem}
\label{lem:itf}
Let $X \in \Xt$ and let $E_{\beta} = [\eta \geq |t| \geq \rho |x|]$ be
a parabolic exterior set associated to $X$.
Consider $\zeta>0$ and $\theta>0$.
Then $|\psi_{X,\beta}/\psi_{X, \beta}^{0} -1| \leq \zeta$ in
$E_{\beta} \cap [t - \gamma_{1}(x) \in {\mathbb R}^{+} e^{i[-\theta, \theta]}]
\cap [x \in B(0,\delta(\zeta,\theta))]$ for $N(\beta)=1$.
The same inequality is true for $N(\beta) \geq 2$ if $\rho >0$ is big enough.
\end{lem}
\begin{proof}
Consider the change of coordinates $(x,z)=(x, t-\gamma_{1}(x))$. We have
\[ \psi_{\beta, E}^{0} = \frac{-1}{\nu(\beta) v(0,0)} \frac{1}{z^{\nu(\beta)}}
+ Res(X_{\beta}^{0},(0,0)) \ln z + h(z) + b(x) \]
where $h$ is a $O(1/z^{\nu(\beta)-1})$ meromorphic function
and $b(x)$ is a holomorphic function in the neighborhood of $0$.
In a sector of bounded angle in the variable $z$ we have that
$\psi_{\beta, E}^{0}z^{\nu(\beta)}$ is bounded both by above and by below.

We define $K(x,z) = \psi_{\beta, E}(x,z)-\psi_{\beta, E}^{0}(x,z)$.
Consider the function $J = x$ if $N(\beta)=1$ and
$J = x/z$ if $N(\beta)>1$. We have
\[ v(0,z) z^{\nu(\beta)+1} \frac{\partial{K}}{\partial{z}} =
\frac{v(0,z) z^{\nu(\beta)+1-s_{1}}}{v(x,z+\gamma_{1}(x))
\prod_{j=2}^{p} (z+\gamma_{1}(x)-\gamma_{j}(x))^{s_{j}}} -1
= O(J) . \]
Thus $\partial{K}/\partial{z}$ is a
$O(J/z^{\nu(\beta)+1})$.
Let $(x,r e^{i \omega}) \in E_{\beta} \cap (|arg z| \leq \theta)$.
We obtain
\[ |K(x,\eta e^{i \omega})| \leq |K(x,\eta)|  +
\left|{
\int_{\eta}^{\eta e^{i \omega}} \frac{\partial K}{\partial z} dz }\right|
= O(x) + O(x) = O(x) \ \forall  \omega \in [-\theta, \theta]. \]
Consider $\gamma:[0,1] \to {\mathbb C}^{2}$ defined by
$\gamma(\upsilon) = (x, e^{i\omega} [(1-\upsilon) \eta + \upsilon r])$. We
obtain
\[ |K(x,r e^{i \omega}) - K(x,\eta e^{i \omega})|
\leq \left|{ \int_{\gamma} \frac{\partial K}{\partial z} dz
}\right|  \leq \left|{ \int_{0}^{1}  \frac{\partial K}{\partial z}
(\gamma(\upsilon)) \gamma' (\upsilon) d \upsilon }\right| \]
We define $C_{0} \equiv |x|$ if $N(\beta)=1$ and $C_{0} \equiv 1/\rho$ if
$N(\beta)>1$. We get
\[ |K(x,r e^{i \omega}) - K(x,\eta e^{i \omega})|
\leq A C_{0}(x)
\int_{0}^{1} \frac{\eta -r}{{[(1-\upsilon) \eta + \upsilon r]}^{\nu(\beta)+1}}
d \upsilon
\leq B \left|{ \frac{C_{0}(x)}{z^{\nu(\beta)}} }\right|  \]
for some $A,B>0$.
We obtain $|K(x,z)| = O(x)+O(C_{0}(x)/z^{\nu(\beta)})$ and then
\[  \left|{ \frac{\psi_{X, \beta}}{\psi_{X, \beta}^{0}} -1 }\right| =
\left|{ \frac{K}{{\psi}_{\beta, E}^{0}} }\right| \leq D
\left|{ C_{0}(x) }\right|  \]
in $E_{\beta} \cap [|\arg z| \leq \theta] \cap
[x \in B(0,\delta(\zeta,\theta))]$ for some $D>0$ depending on
$\theta$.
\end{proof}
\begin{rem}
The previous lemma implies that
$\psi_{X, \beta}
\sim 1/(x^{d_{\beta}} (t-\gamma_{1}(x))^{\nu(\beta)})$ in a
parabolic exterior set $E_{\beta}$ for
$|\arg(t-\gamma_{1}(x))|$ bounded.
\end{rem}
\begin{pro}
Let $X \in \Xt$ and let
$E_{\beta} = [\eta \geq |t| \geq \rho |x|]$ be a parabolic exterior set
associated to $X$. Consider
$C \in SC_{\mu X}^{\beta,\eta}(r, \lambda)$ for $r \lambda$ in a
neighborhood of $0$ and $\mu \in {\mathbb S}^{1}$. Then $C$ is
contained in a sector
centered at $t=\gamma_{1}(r \lambda)$ of angle lesser than $\theta$
for some $\theta>0$ independent of $r, \lambda$, $C$ and $\mu$.
\end{pro}
\begin{proof}
We use the notations in lemma \ref{lem:itf}.
We have that the extrema of a connected component
of $[|t|=\eta] \setminus TE_{\mu X}^{\beta,\eta}(r, \lambda)$
lie in an angle centered at $z=0$ of angle similar to $\pi/\nu(\beta)$.
Then it is enough to prove that
$\Gamma=
\Gamma(\mu \lambda^{d_{\beta}} X_{\beta,E},(r, \lambda,t_{0}),E_{\beta})$
lies in a sector of bounded angle for
$t_{0} \in TE_{\mu X}^{\beta,\eta}(r, \lambda)$.

Denote
$\psi^{0} = -1/(\nu(\beta)v(0,0) z^{\nu(\beta)})$. We have
$\lim_{z \to 0} \psi_{\beta, E}/\psi^{0}=1$ in big sectors;
we can suppose that
$|\psi_{\beta, E}/\psi^{0}-1| < \zeta$ for arbitrary $\zeta>0$
by taking $0< \eta <<1$. Since the set
$(\psi_{\beta, E}/(\mu \lambda^{d_{\beta}}))(\Gamma)$
is contained in
$(\psi_{\beta, E}/(\mu \lambda^{d_{\beta}}))(r, \lambda,t_{0})+{\mathbb R}$
then it lies in a sector of angle similar to $\pi$.
Since $\psi_{\beta, E}/\psi^{0} \sim 1$ then $\Gamma$ lies in a sector of
center $t=\gamma_{1}(r, \lambda)$ and angle close to $\pi/\nu(\beta)$.
\end{proof}
\begin{rem}
\label{rem:bditfpet}
We have that
$\psi_{X, \beta} \sim 1/(x^{d_{\beta}} (t-\gamma_{1}(x))^{\nu(\beta)})$
in $E_{\beta} \cap \overline{C}$ for a parabolic exterior set $E_{\beta}$
and all $C \in SC_{\mu X}^{\beta,\eta}$
\end{rem}
\subsection{Nature of the polynomial vector fields}
The study of polynomial vector fields related to stability properties
of unfoldings of elements $h \in \diff{1}{}$
has been introduced in \cite{DES}. Their choices
are associated with the elements in the deformation
whereas ours depend on the infinitesimal properties of the unfolding.
\subsubsection{Directions of instability}
\label{subsubsec:diruns}
Let $M_{\beta}$ be a magnifying glass set associated to a vector
field $X \in \Xt$. We consider
\[ X_{\beta}(\lambda) = \lambda^{m_{\beta}}
C (w-w_{1})^{s_{1}} \hdots (w-w_{p})^{s_{p}} \partial / \partial{w} \]
where $C \in {\mathbb C}^{*}$ and $w_{j} \in {\mathbb C}$ for all
$j \in \{1,\hdots,p\}$. Denote
$r_{\beta}^{j}(X) = Res(X_{\beta}(1),w_{j})$ for $1 \leq j \leq p$.
Consider the set $sum_{\beta}(X)$ whose elements are the non-vanishing
sums of the form $\sum_{j \in E} r_{\beta}^{j}$ for any
$E \subset \{1,\hdots,p\}$. We define
\[ B_{\beta}(X) = \{ (\lambda,\mu) \in {\mathbb S}^{1} \times  {\mathbb S}^{1}:
sum_{\beta} \cap \lambda^{m_{\beta}} \mu i {\mathbb R} \neq \emptyset \} . \]
We denote ${\mathbb S}^{1}/\sim$ the quotient of ${\mathbb S}^{1}$
by the equivalence relation identifying $\mu$ and $-\mu$. We denote
by $\tilde{B}_{\beta}(X) \subset {\mathbb S}^{1}/\sim \times
{\mathbb S}^{1}/\sim$
the quotient of $B_{\beta}(X)$. Now we define
\[ B_{\beta, \lambda}(X) =
\{ \mu \in {\mathbb S}^{1} : (\lambda,\mu) \in B_{\beta}(X) \}
\ \ {\rm and} \ \  B_{\beta}^{\mu}(X) =
\{ \lambda \in {\mathbb S}^{1} : (\lambda,\mu) \in B_{\beta}(X) \} . \]
In an analogous way we can define
$\tilde{B}_{\beta, \lambda}(X) \subset {\mathbb S}^{1}/\sim$
and $\tilde{B}_{\beta}^{\mu}(X) \subset {\mathbb S}^{1}/\sim$
for $\lambda, \mu \in {\mathbb S}^{1}/\sim$.
Roughly speaking we claim that $Re(\mu X)$ has a stable
behavior in $I_{\beta}$ at the direction
$x \in {\mathbb R}^{+} \lambda$ for
$(\lambda,\mu) \not \in B_{\beta}(X)$.
We define $B_{X}$ as the union of $B_{\beta}(X)$ for every
magnifying glass set $M_{\beta}$ associated to $X$.
Analogously we can define $B_{X, \lambda}$, $B_{X}^{\mu}$,
$\tilde{B}_{X, \lambda}$ and $\tilde{B}_{X}^{\mu}$.
The sets $B_{X, \lambda}$ and $B_{X}^{\mu}$ are finite for
all $\lambda,\mu \in {\mathbb S}^{1}$. Moreover we have
$B_{X, \lambda'} \cap B_{X, \lambda} = \emptyset$ and
$B_{X}^{\lambda'} \cap B_{X}^{\lambda} = \emptyset$ for all
$\lambda' \in {\mathbb S}^{1}$ in a pointed neighborhood of $\lambda$.
\begin{rem}
\label{rem:spndepnf}
Let $X \in \Xt$. The sets in the dynamical splitting
depend only on $Sing X$ whereas $B_{\beta}(X)$ depends on
$Sing X$ and $Res(X)$ for all magnifying glass set $M_{\beta}$.
\end{rem}
\subsubsection{Non-parabolic exterior sets}
\label{subsub:npes}
Let $E_{\beta w_{1}}$ be a non-parabolic exterior set where
$w_{1} \in {\mathbb C}$. Thus we have
\[ X=x^{m_{\beta}} h(x,w) (w-w_{1}(x)) (w-w_{2}(x))^{s_{2}} \hdots
(w-w_{p}(x))^{s_{p}} \partial / \partial w \]
in $M_{\beta}$ where $w_{1}(0)=w_{1}$ and $h(x,w)-h(0,0) \in (x)$.
This expression implies
\[ X = x^{m_{\beta}} (r_{\beta}^{1})^{-1} (w-w_{1}(x))
(1+H(x,w)) \partial / \partial w \]
in $E_{\beta w_{1}}$ for some $H \in (x,w-w_{1}) = (x,w-w_{1}(x))$.

Fix $\mu \in {\mathbb S}^{1}$ and a compact set
$K_{X}^{\mu} \subset {\mathbb S}^{1} \setminus B_{X}^{\mu}$. By definition
of $B_{X}^{\mu}$ we obtain that
$\lambda^{m_{\beta}} \mu/r_{\beta}^{1} \not \in i {\mathbb R}$ for all
$\lambda \in K_{X}^{\mu}$. This implies
$\lambda^{m_{\beta}} \mu (r_{\beta}^{1})^{-1}
(1+H(r\lambda,w_{1}(r\lambda))) \not \in i {\mathbb R}$
for $(r,\lambda) \in [0,r_{0}) \times K_{X}^{\mu}$ for some $r_{0}>0$
since $K_{X}^{\mu}$ is compact and $H(x,w_{1}(x)) \in (x)$.
We deduce that the singular point $w=w_{1}(x_{0})$ of
$Re(\mu X)_{|x=x_{0}}$ is not a center  for $x_{0} \in (0,r_{0})K_{X}^{\mu}$.
Hence, it is either an attracting or a repelling point.

The set $E_{\beta w_{1}}$ is of the form $|w-w_{1}|<c$ for some
$0<c<<1$. The vector field $Re(\mu X)_{|x=r \lambda}$ and the set
$\partial E_{\beta w_{1}}$ are tangent at the set
\[ TE_{\mu X}^{\beta w_{1},c}(r, \lambda) =
\left[{
\frac{\lambda^{m_{\beta}} \mu}{r_{\beta}^{1}}
\frac{w-w_{1}(r \lambda)}{w-w_{1}} (1+H(r\lambda,w)) \in i {\mathbb R}
}\right]
\cap [|w-w_{1}(0)|=c] . \]
The function $(w-w_{1}(r\lambda))/(w-w_{1})$ tends to $1$ when $r \to 0$
in $|w-w_{1}|=c$. Moreover since $H \in (x,w-w_{1})$ we obtain that
$TE_{\mu X}^{\beta w_{1},c}(r, \lambda) =\emptyset$
for $r \in [0,r_{0}(c))$ and $\lambda \in K_{X}^{\mu}$.
Then  $Re(s \mu X)$ points towards
$\dot{E}_{\beta w_{1}}(x)$ at $\partial {E}_{\beta w_{1}}(x)$
for all $x \in (0,r_{0}) K_{X}^{\mu}$
and either $s=-1$ or $s=1$.
As a consequence $E_{\beta w_{1}} \cap [x=x_{0}]$
is in the basin of attraction of $(x_{0},w_{1}(x_{0}))$ by
$Re(s \mu X)$ for $x_{0} \in (0,r_{0})K_{X}^{\mu}$.
\subsubsection{Connexions at $\infty$}
We already described the dynamics of $Re(\mu X)$ in the exterior
sets for $\mu \in {\mathbb S}^{1}$ and $X \in \Xt$. Next we analyze the
dynamics of $Re(\mu X)$ in the intermediate sets.

Let $Y=C(w-w_{1})^{s_{1}} \hdots (w-w_{p})^{s_{p}} \partial / \partial{w}$
be a polynomial vector field such that
$\nu(Y)=s_{1}+\hdots+s_{p}-1 \geq 1$. Every vector field
$X_{\beta}(\lambda)$ associated to a magnifying glass set is of this form.
We want to characterize the behavior of $Y$ in the neighborhood of
$\infty$. We define the set $Tr_{\to \infty}(Y)$
of trajectories $\gamma:(c,d) \to {\mathbb C}$ of $Re(Y)$
such that $c \in {\mathbb R} \cup \{-\infty\}$, $d \in {\mathbb R}$
and $\lim_{\zeta \to d} \gamma(\zeta) = \infty$. In an analogous way
we define $Tr_{\leftarrow \infty}(Y)=Tr_{\to \infty}(-Y)$.
We define $Tr_{\infty}(Y) = Tr_{\leftarrow \infty}(Y) \cup Tr_{\to \infty}(Y)$.

We consider a change of coordinates $z=1/w$. The meromorphic vector field
\[ Y = \frac{-C (1-w_{1}z)^{s_{1}} \hdots (1-w_{p}z)^{s_{p}}}{z^{\nu(Y)-1}}
\frac{\partial}{\partial z} \]
is analytically conjugated to
$1/(\nu(Y) z^{\nu(Y)-1}) \partial / \partial{z} =
(z^{\nu(Y)})^{*}(\partial/\partial{z})$ in a neighborhood of $\infty$.
We have $Tr_{\to \infty}(\partial/\partial{z})= {\mathbb R}^{-}$
and $Tr_{\leftarrow \infty}(\partial/\partial{z})= {\mathbb R}^{+}$.
Hence the set $Tr_{\infty}(Y)$ has $2 \nu(Y)$ trajectories
in the neighborhood of $\infty$, there is exactly one of them which is
tangent to
$\arg (w) = -\arg(C)/\nu(Y) + k \pi/\nu(Y)$ for
$0 \leq k < 2 \nu(Y)$. The even values of $k$ correspond to
$Tr_{\to \infty}(Y)$.

The complementary of the set $Tr_{\infty}(Y)  \cup \{ \infty \}$
has $2 \nu(Y)$ connected components in the neighborhood of $w=\infty$.
Each of these components is called an {\it angle}, the boundary of an angle
contains exactly one $\to \infty$-trajectory and one
$\leftarrow \infty$-trajectory.

We say that $Re(Y)$ has $\infty$-connections if
$Tr_{\to \infty}(Y) \cap Tr_{\leftarrow \infty}(Y) \neq \emptyset$.
In other words there exists a trajectory
$\gamma:(c_{-1},c_{1}) \to {\mathbb C}$
of $Re(Y)$ such that $c_{-1},c_{1} \in {\mathbb R}$ and
$\lim_{\zeta \to c_{s}} \gamma(\zeta)=\infty$ for all $s \in \{-1,1\}$.
The notion of connexion at $\infty$ has been introduced in \cite{DES}
for the study of deformations of elements of $\diff{1}{}$.
%
%

We define the $\alpha$ and $\omega$ limits $\alpha^{Y}(P)$ and
$\omega^{Y}(P)$ respectively of a point $P \in {\mathbb C}$
by the vector field $Re(Y)$. If $P \in Tr_{\to \infty}(Y)$
we denote $\omega^{Y}(P) = \{\infty\}$ whereas if
$P \in Tr_{\leftarrow \infty}(Y)$ we denote $\alpha^{Y}(P)=\{\infty\}$.
\begin{lem}
Let $Y \in {\mathcal X} \cn{}$ be a polynomial vector field
such that $\nu(Y) \geq 1$. Then $\omega^{Y}(w_{0})=\{ \infty \}$
is equivalent to $w_{0} \in Tr_{\to \infty}(Y)$. Analogously
$\alpha^{Y}(w_{0})=\{ \infty \}$
is equivalent to $w_{0} \in Tr_{\leftarrow \infty}(Y)$
\end{lem}
\begin{proof}
The vector field $Y$ is a ramification of a regular vector field
in a neighborhood of $\infty$. Thus there exists an open neighborhood $V$
of $\infty$ and $c \in {\mathbb R}^{+}$ such that
\[ {\rm exp}(c Y)(V \setminus Tr_{\to \infty}(Y)) \cap V = \emptyset
\ \ {\rm and} \ \
{\rm exp}(-c Y)(V \setminus Tr_{\leftarrow \infty}(Y)) \cap V = \emptyset. \]
We are done since $w_{0} \not \in Tr_{\to \infty}(Y)$ implies
$\omega^{Y}(w_{0}) \cap (\pn{1} \setminus V) \neq \emptyset$.
\end{proof}
We denote by ${\mathcal X}_{\infty}\cn{}$ the set of polynomial vector
fields in ${\mathcal X} \cn{}$ such that $\nu(Y) \geq 1$ and
$2 \pi i \sum_{P \in S} Res(Y,P) \not \in {\mathbb R} \setminus \{0\}$
for all subset $S$ of $Sing Y$.
\begin{lem}
\label{lem:aopol}
Let $Y \in {\mathcal X}_{\infty} \cn{}$. Then
\begin{itemize}
\item $Re(Y)$ has no $\infty$-connections.
\item $\omega^{Y}(w_{0}) \neq \{ \infty \}$ implies that
$\sharp \omega^{Y}(w_{0}) =1$ and
$\omega^{Y}(w_{0}) \cap Sing Y \neq \emptyset$.
\end{itemize}
\end{lem}
\begin{proof}
Let $\Omega$ the unique
meromorphic 1-form defined by $\Omega(Y)=1$.
Suppose that $\gamma:(c_{-1},c_{1}) \to {\mathbb C}$ is
an $\infty$-connexion of $Re(Y)$. Consider the connected
component $U$ of $\pn{1} \setminus (\gamma(c_{-1},c_{1}) \cup \{\infty\})$
such that $Re(iY)$ points towards $U$ at $\gamma$.

There exists a holomorphic integral $\psi$ of the time form of $Y$
in a neighborhood of $w=\infty$ such that $\psi \sim 1/w^{\nu(Y)}$.
Since $\psi(\infty)=0$ then
theorem of the residues
implies that
$2 \pi i \sum_{P \in Sing Y \cap U} Res(Y,P) =
c_{1} - c_{-1} \in {\mathbb R}^{+}$.
This is a contradiction.

It is enough to prove that
$\omega^{Y}(w_{0}) \cap ({\mathbb C} \setminus Sing Y) =\emptyset$
since $\omega^{Y}(w_{0})$ is connected.
Suppose $P \in \omega^{Y}(w_{0}) \cap ({\mathbb C} \setminus Sing Y)$.
Denote $\gamma:[0,\infty) \to {\mathbb C}$ the
trajectory of $Re(Y)$ passing through $w_{0}$.
Consider a germ of transversal $h$ to the vector field $Re(Y)$
passing through $P$. There exists some $\eta>0$ such that
${\rm exp}((0,\eta]Y)(h) \cap h = \emptyset$. There also exists an
increasing sequence of positive real numbers $j_{n} \to \infty $ such that
$\gamma(j_{n}) \in h$, $\lim_{n \to \infty} \gamma(j_{n})=P$ and
$\gamma(j_{n},j_{n+1}) \cap h = \emptyset$ for all $n \in {\mathbb N}$.

Consider a holomorphic integral $\psi_{0}$ of the time form of $Y$ defined in
the neighborhood of $P$. Let $L_{n}$ be the segment of $h$ whose boundary
is $\{\gamma(j_{n}),\gamma(j_{n+1}) \}$. Denote by $V_{n}$ the bounded
component of ${\mathbb C} \setminus (\gamma[j_{n},j_{n+1}] \cup L_{n})$.
By the theorem of the residues we obtain
\[ \int_{\gamma[j_{n},j_{n+1}]} \Omega +
(\psi_{0}(\gamma(j_{n}))-\psi_{0}(\gamma(j_{n+1}))) =
\pm 2 \pi i \sum_{P \in V_{n} \cap Sing Y} Res(Y,P)  \]
By making $n$ to tend to $\infty$ we deduce that
there exists a subset $S$ of $Sing Y$ such that
$\pm 2 \pi i \sum_{P \in S} Res(Y,P) \in [\eta, \infty)$.
That is a contradiction.
\end{proof}
\begin{cor}
Let $X \in \Xt$. Consider a magnifying glass set $M_{\beta}$ associated to
$X$. Then
\begin{itemize}
\item $Re(\mu X_{\beta}(\lambda))$ has no $\infty$-connections.
\item  $\omega^{\mu X_{\beta}(\lambda)}(w_{0}) \neq \infty
\Rightarrow \sharp \omega^{\mu X_{\beta}(\lambda)}(w_{0}) =1$ and
$\omega^{\mu X_{\beta}(\lambda)}(w_{0}) \cap Sing X_{\beta}(\lambda) \neq
\emptyset$.
\end{itemize}
for all $(\lambda,\mu) \not \in B_{\beta}(X)$.
\end{cor}
\subsubsection{The graph}
In this subsection we associate an oriented graph to every vector field
$\mu X_{\beta}(\lambda)$ for $(\lambda,\mu) \not \in B_{\beta}(X)$.
\begin{lem}
Let $Y \in {\mathcal X}_{\infty} \cn{}$. Then
$\omega^{Y}: {\mathbb C} \setminus (Tr_{\to \infty}(Y) \cup Sing Y)
\to Sing Y$ is locally constant.
\end{lem}
\begin{proof}
Let $P \in  {\mathbb C} \setminus (Tr_{\to \infty}(Y) \cup Sing Y)$.
Denote $Q=\omega^{Y}(P)$. The singular point $Q$ is not a center
since then $Re(Y)$ would support cycles (lemma \ref{lem:aopol}).
If $Q$ is attracting there is nothing to prove. If $Q$ is parabolic
then $P \in \cup_{\lambda \in D_{1}(Y)} V_{{\rm exp}(Y)}^{\lambda}$.
We are done since
$\cup_{\lambda \in D_{1}(Y)} V_{{\rm exp}(Y)}^{\lambda}$ is open and
$\omega^{Y}(\cup_{\lambda \in D_{1}(Y)} V_{{\rm exp}(Y)}^{\lambda})=Q$.
\end{proof}
We denote by $Reg(Y)$ the set of connected components of
$\mathbb C \setminus (Tr_{\infty}(Y) \cup Sing Y)$.
Its elements are called {\it regions} of $Re(Y)$.
Every $H \in Reg(Y)$ satisfies that $\alpha^{Y}(H)$ and
$\omega^{Y}(H)$ are points.
We denote by $Reg_{j}(Y)$ the set of regions $H$ of $Re(Y)$
such that $\sharp \{\alpha^{Y}(H),\omega^{Y}(H)\}=j$ for $j \in \{1,2\}$.
We associate an oriented graph to
$Re(Y)$ for $Y \in {\mathcal X}_{\infty} \cn{}$.
The set of vertexes is $Sing Y$, the
edges are the regions of $Re(Y)$.
We say that $H \in Reg(Y)$ joins the points
$\alpha^{Y}(H)$ and $\omega^{Y}(H)$. We denote
$\alpha^{Y}(H) \stackrel{H}{\rightarrow} \omega^{Y}(H)$.
The graph obtained in this way is denoted by ${\mathcal G}_{Y}$.
We denote by ${\mathcal NG}_{Y}$ the unoriented graph obtained
from ${\mathcal G}_{Y}$ by removing the reflexive edges and the
orientations of the edges.

An angle is always contained in a region of $Re(Y)$.
Such a region is characterized by the angles that it contains.
Let $A$ be an angle of the polynomial vector field $Y$.
We denote by $\gamma_{\to \infty}^{A}$ the trajectory of $Tr_{\to \infty}$
contained in the closure of $A$. The definition of
$\gamma_{\leftarrow \infty}^{A}$ is analogous.
\begin{lem}
\label{lem:atleona}
Let $Y \in {\mathcal X}_{\infty} \cn{}$. Consider $H \in Reg(Y)$.
Then $H$ contains an angle $A$. Moreover
$\alpha^{Y}(\gamma_{\to \infty}^{A})=\alpha^{Y}(H)$ and
$\omega^{Y}(\gamma_{\leftarrow \infty}^{A})=\omega^{Y}(H)$.
\end{lem}
\begin{proof}
Let $P \in ({\mathbb C} \setminus Sing Y) \cap \partial H$; such a
point exists since $Tr_{\infty}(Y)$ is contained in the complementary
of $H$. Since $\alpha^{Y}$ and $\omega^{Y}$ are locally constant then
either $\alpha^{Y}(P)=\infty$ or $\omega^{Y}(P)=\infty$.
We have that $P \in \overline{H}$, thus there are points of $H$ in
every neighborhood of $\infty$. As a consequence $H$ contains at least
an angle $A$. The relations
$\alpha^{Y}(\gamma_{\to \infty}^{A})=\alpha^{Y}(H)$ and
$\omega^{Y}(\gamma_{\leftarrow \infty}^{A})=\omega^{Y}(H)$
can be deduced of the locally constant character of
$\alpha^{Y}$ and $\omega^{Y}$.
\end{proof}
\begin{lem}
\label{lem:ssictr}
Let $Y \in {\mathcal X}_{\infty} \cn{}$. Then we have
$Sing Y \subset \overline{Tr_{\infty}(Y)}$.
\end{lem}
\begin{proof}
Let $P \in Sing Y$. Suppose that $V \cap  Tr_{\infty}(Y) = \emptyset$
for some connected neighborhood $V$ of $P$. Let $H$ be the region of $Re(Y)$
containing $V \setminus \{P\}$. Since $P$ is attracting, repelling or
parabolic then  either $\alpha^{Y}(H)=P$ or
$\omega^{Y}(H)=P$. Consider an angle $A \subset H$. We obtain
$P \in \overline{\gamma_{\leftarrow \infty}^{A} \cup \gamma_{\to \infty}^{A}}
\subset \overline{Tr_{\infty}(Y)}$.
\end{proof}
\begin{lem}
Let $Y \in {\mathcal X}_{\infty} \cn{}$. Consider $H \in Reg_{1}(Y)$.
Then $H$ contains exactly one angle.
\end{lem}
\begin{proof}
Let $A$ be an angle contained in $H$. Denote $P=\alpha^{Y}(H)=\omega^{Y}(H)$.
By lemma \ref{lem:atleona} we have that
$\gamma= \{\infty\} \cup \gamma_{\to \infty}^{A} \cup
\gamma_{\leftarrow \infty}^{A} \cup \{P\}$ is a closed simple curve.
Let $V$ the connected component of
$\pn{1} \setminus \gamma$ containing $A$. The set
$Tr_{\infty}(Y) \cap V$ is empty since $A$ is the only angle contained
in $V$. By lemma \ref{lem:ssictr} we have that $V \cap Sing Y = \emptyset$.
Hence $H$ is equal to $V$ and contains only one angle.
\end{proof}
\begin{lem}
\label{lem:reg2}
Let $Y \in {\mathcal X}_{\infty} \cn{}$. Consider $H \in Reg_{2}(Y)$.
Then $H$ contains exactly two
angles. Moreover ${\mathbb C} \setminus H$ has two connected components
$H_{1}$ and $H_{2}$ such that $\alpha^{Y}(H) \in H_{1}$ and
$\omega^{Y}(H) \in H_{2}$.
\end{lem}
\begin{proof}
Let $A_{1}$ be an angle contained in $H$. Fix a trajectory $\gamma_{0}$
of $Re(Y)$ contained in $H$.
Denote
$\gamma_{1}=\gamma_{0} \cup \gamma_{\to \infty}^{A_{1}} \cup
\gamma_{\leftarrow \infty}^{A_{1}} \cup \{\alpha^{Y}(H),\omega^{Y}(H)\}$.
Consider the connected component $V_{1}$ of ${\mathbb C} \setminus \gamma_{1}$
containing $A_{1}$. Since $V_{1}$ contains only one angle then
$V_{1} \subset H$. By proceeding like in lemma \ref{lem:atleona}
we can prove that there exists an angle $A_{2}$ contained in
$H \setminus (V_{1} \cup \gamma_{0})$. Let $V_{2}$ be the connected
component of the set ${\mathbb C} \setminus (\gamma_{0} \cup
\gamma_{\to \infty}^{A_{2}} \cup
\gamma_{\leftarrow \infty}^{A_{2}} \cup \{\alpha^{Y}(H),\omega^{Y}(H)\})$
such that $A_{2} \subset V_{2}$. Clearly we have $A_{2} \neq A_{1}$
and $H=V_{1} \cup \gamma_{0} \cup V_{2}$. Now
\[ {\mathbb C} \setminus (
\gamma_{\to \infty}^{A_{1}} \cup \gamma_{\leftarrow \infty}^{A_{1}} \cup
\gamma_{\to \infty}^{A_{2}} \cup \gamma_{\leftarrow \infty}^{A_{2}} \cup
\{\alpha^{Y}(H),\omega^{Y}(H)\}) \]
has three connected components $H$, $J_{1}$ and $J_{2}$
such that
\[ \partial{J_{1}} = \gamma_{\to \infty}^{A_{1}} \cup
\gamma_{\to \infty}^{A_{2}} \cup  \{\alpha^{Y}(H)\} \ \ {\rm and} \ \
\partial{J_{2}} = \gamma_{\leftarrow \infty}^{A_{1}} \cup
\gamma_{\leftarrow \infty}^{A_{2}} \cup  \{ \omega^{Y}(H)\} . \]
Then $H_{1}=J_{1} \cup \partial{J_{1}}$ and
$H_{2}=J_{2} \cup \partial J_{2}$ are the connected components of
${\mathbb C} \setminus H$.
\end{proof}
\begin{cor}
\label{cor:nocycles}
Let $Y \in {\mathcal X}_{\infty} \cn{}$. Then
${\mathcal NG}_{Y}$ has no cycles.
\end{cor}
\begin{proof}
Consider an edge $P \stackrel{H}{\to} Q$ of ${\mathcal G}_{Y}$
with $P \neq Q$. Consider the notations in the previous lemma.
The fixed points are divided in two sets $H_{1} \cap Sing Y$ and
$H_{2} \cap Sing Y$. The only edge of ${G}_{Y}$ joining a
vertex in the former set with a vertex in the latter set
(or vice-versa) is $P \stackrel{H}{\to}  Q$.
Clearly ${\mathcal NG}_{Y}$ has no cycles.
\end{proof}
\begin{pro}
\label{pro:connected}
Let $Y \in {\mathcal X}_{\infty} \cn{}$. Then the graph
${\mathcal NG}_{Y}$ is connected.
\end{pro}
\begin{proof}
Let $G_{1}, \hdots, G_{l}$ be the set of vertexes of the $l$
connected components of ${\mathcal NG}_{Y}$.
We define the open set
$V_{j}=(\alpha^{Y})^{-1}(G_{j}) \cup (\omega^{Y})^{-1}(G_{j})$
for all $j \in \{1,\hdots,l\}$.
The lack of $\infty$-connexions implies
$\cup_{j=1}^{l} V_{j} = {\mathbb C}$. Moreover
$V_{j} \cap V_{k} = \emptyset$ if $j \neq k$ since otherwise
$G_{j} = G_{k}$. Clearly $l=1$ since ${\mathbb C}$ is connected.
\end{proof}
\begin{cor}
\label{cor:edge}
Let  $Y \in {\mathcal X}_{\infty} \cn{}$. Then
$\sharp Reg_{2}(Y) = \sharp Sing Y -1$.
\end{cor}
Let $Y \in {\mathcal X} \cn{}$. Consider $y_{0} \in Sing Y$.
We define $\nu_{Y}(y_{0})$ as the only
element of ${\mathbb N} \cup \{0\}$ such that
$Y(y) \in (y-y_{0})^{\nu_{Y}(y_{0})+1} \setminus (y-y_{0})^{\nu_{Y}(y_{0})+2}$.
\begin{pro}
\label{pro:redge}
Let $Y \in {\mathcal X}_{\infty} \cn{}$. Consider $y_{0} \in Sing Y$.
Then there exist exactly $2 \nu_{Y}(y_{0})$ regions of $Re(Y)$
contained in $(\alpha^{Y},\omega^{Y})^{-1}(y_{0},y_{0})$.
\end{pro}
\begin{proof}
If $y_{0}$ is not parabolic the result is obvious since
on the one hand $\nu_{Y}(y_{0})=0$ and on the other hand
$(\alpha^{Y},\omega^{Y})^{-1}(y_{0},y_{0})=\{ y_{0} \}$.

Let ${Y}_{0}$ be the germ of $Y$ in the neighborhood of $y_{0}$,
we have $\nu({Y}_{0})=\nu_{Y}(y_{0})$. Consider the strict transform
$\tilde{Y}$ of $Re(Y)$ by the real blow-up
$\pi(r, \lambda)= y_{0}+r \lambda$. By the discussion in
section \ref{subsec:topbeh} there exists a unique region of $Re(Y)$
adhering to
$[(r,\lambda) \in \{0\} \times [\lambda_{0},\lambda_{0}e^{i \pi/\nu(Y_{0})}]]$
for all $\lambda_{0} \in D(Y_{0})$. In this way we find $2 \nu_{Y}(y_{0})$
regions in $(\alpha^{Y},\omega^{Y})^{-1}(y_{0},y_{0})$.
Any other region would adhere to a single point in $D(Y_{0})$. It
would be both attracting and repelling for $\tilde{Y}$; that is
impossible.
\end{proof}
\begin{cor}
Let  $Y \in {\mathcal X}_{\infty} \cn{}$. Then
$\sharp Reg(Y) = 2 \nu(Y) - \sharp(Sing Y) +1$.
\end{cor}
Let  $Y \in {\mathcal X}_{\infty} \cn{}$. Consider a trajectory
$\gamma_{H} \subset H$ for every region $H \in Reg_{2}(Y)$.
There exists $\rho_{0}>0$ such that
\begin{equation}
\label{equ:size}
\left\{{
\begin{array}{l}
Sing Y \subset B(0,\rho_{0}) \  {\rm and} \
\sharp T_{Y}^{\rho} = 2 \nu(Y) \
{\rm for} \  {\rm all} \ \rho \geq \rho_{0}. \\
\gamma_{H} \subset B(0,\rho_{0}) \ {\rm for \ all} \  H \in Reg_{2}(Y).
\end{array} }\right.
\end{equation}
Let $P \in \overline{B}(0,\rho)$. We define $\omega^{Y, \rho}(P)=\infty$
if $It(Y,P,\overline{B}(0,\rho))$ does not contain $(0,\infty)$.
Otherwise we define
$\omega^{Y, \rho}(P)=\omega^{Y}(P)$. We define $\alpha^{Y, \rho}$ in
an analogous way. Denote by $Reg(Y,\rho)$ the
set of connected components of
\[ B(0,\rho) \setminus
((\alpha^{Y, \rho})^{-1}(\infty) \cup (\omega^{Y, \rho})^{-1}(\infty)
\cup Sing Y). \]
Denote
\[ Reg_{j}(Y,\rho) = \{ H \in Reg(Y,\rho) :
\sharp \{\alpha^{Y, \rho}(H),\omega^{Y, \rho}(H)\}=j \} \] for $j \in
\{1,2\}$. The set of connected components of $B(0,\rho) \setminus
(Sing Y \cup  \cup_{H \in Reg(Y,\rho)} \overline{H})$ will be called
$Reg_{\infty}(Y,\rho)$.
The dynamics of $Re(Y)$ in ${\mathbb C}$ and
$B(0,\rho_{0})$ is analogous.
\begin{pro}
\label{pro:gofidi}
Let $Y \in {\mathcal X}_{\infty} \cn{}$. Consider $\rho >>0$.
There exist bijections $F:Reg(Y,\rho) \to Reg(Y)$ and
$G:Reg_{\infty}(Y,\rho) \to Tr_{\infty}(Y)$ such that
\begin{itemize}
\item $H \subset F(H)$ for all $H \in Reg(Y,\rho)$ \item
$\sharp(\partial H \cap T_{Y}^{\rho})=j$ for all $H \in Reg_{j}(Y,\rho)$
and $j \in \{1,2\}$. \item $\sharp(\partial J \cap T_{Y}^{\rho})=1$
for each connected component $J$ of $H \setminus \gamma_{H}$ and
$H \in Reg_{2}(Y,\rho)$. \item $G(K) \cap B(0,\rho) \subset K$ for
all $K \in Reg_{\infty}(Y,\rho)$.
\end{itemize}
\end{pro}
\begin{proof}
We define $F_{1}(H)$ as the element of $Reg(Y,\rho)$
containing $\gamma_{H}$ for $H \in Reg(Y)$. Every $H \in Reg(Y,\rho)$
is contained in a unique $F(H) \in Reg(Y)$. It is clear that
$F \circ F_{1} \equiv Id$. This implies
$\sharp(Reg_{j}(Y,\rho)) \geq \sharp (Reg_{j}(Y))$ for $j \in \{1,2\}$.

Let $H \in Reg(Y,\rho)$.
We have $\partial H \cap \partial B(0,\rho)=\partial H \cap T_{Y}^{\rho}$.
Thus we obtain $\sharp(\partial H \cap T_{Y}^{\rho}) \geq 1$.
Let $H \in Reg_{2}(Y,\rho)$.
Every connected component of $F(H) \setminus \gamma_{H}$ contains
at least a point in $\partial H \cap T_{Y}^{\rho}$
and then $\sharp(\partial H \cap T_{Y}^{\rho}) \geq 2$. We have
\[ 2 \nu(Y) = \sharp T_{Y}^{\rho} \geq
\sharp Reg_{1}(Y,\rho) + 2 \sharp Reg_{2}(Y,\rho) \geq
\sharp Reg_{1}(Y) + 2 \sharp Reg_{2}(Y) = 2 \nu(Y) . \]
Hence all the inequalities are indeed equalities.
We obtain $\sharp Reg_{j}(Y,\rho) = \sharp Reg_{j}(Y)$
and  $\sharp(\partial H \cap T_{Y}^{\rho})=j$ for all $j \in \{1,2\}$
and $H \in Reg_{j}(Y,\rho)$. We deduce that $F_{1}=F^{\circ (-1)}$
and that $\{ \alpha^{Y, \rho}(Q), \omega^{Y, \rho}(Q) \} \subset Sing Y$
for all $Q \in T_{Y}^{\rho}$.

Let $l$ be a connected component of $\partial B(0,\rho) \setminus T_{Y}^{\rho}$
such that $Re(sY)$ points towards $B(0,\rho)$ for some $s \in \{-1,1\}$.
We claim that ${\rm exp}(s(0,\infty)Y)(l)$ is a connected component of
$Reg_{\infty}(Y,\rho)$. Suppose $s=1$ without lack of generality.
Since $\omega^{Y, \rho}(\partial l) \subset Sing Y$ and
${\mathcal NG}_{Y}$ is connected then
$\omega^{Y, \rho}(l) = \omega^{Y, \rho}(\partial l)$
is a singleton contained in $Sing Y$. The claim is proved,
it implies $\sharp Reg_{\infty}(Y,\rho)=2 \nu(Y)$.
There exists a unique $\gamma(l) \in Tr_{\infty}(Y)$ such that
$\gamma(l) \cap l \neq \emptyset$. The mapping
$G(K) = \gamma(\partial K \cap (\partial B(0,\rho) \setminus T_{Y}^{\rho}))$
is the one we are looking for.
\end{proof}
\subsection{Assembling the dynamics of the polynomial vector fields}
Let $X$ in $\Xt$. Throughout this section $K_{X}^{\mu}$ is some compact connected
set contained in ${\mathbb S}^{1} \setminus B_{X}^{\mu}$.
Fix a magnifying glass set $M_{\beta} = [|w| \leq \rho]$ and $\mu \in {\mathbb S}^{1}$.
Given $P$ in $\dot{I}_{\beta} \cap [x \in [0,\delta_{0}) K_{X}^{\mu}]$ we have that
either
${\rm exp}(c(P) \mu \lambda^{m_{\beta}} X_{\beta,M})(P)$ belongs to $[|w| = \rho]$ or
to $[|w-\zeta| = r(\zeta)]$
for some $\zeta \in S_{\beta}$ where
$c(P) = \sup \Gamma(\mu \lambda^{m_{\beta}} X_{\beta,M}, P, I_{\beta})$.
We denote $\omega_{\beta}^{\mu X}(P) = \infty$ and
$\omega_{\beta}^{\mu X}(P) = \zeta$ respectively. The definition of
$\alpha_{\beta}^{\mu X}(P)$ is analogous.
We define $Reg^{*}(I_{\beta}, \mu X, K_{X}^{\mu})$ the set of connected
components of
\[  [I_{\beta} \cap [x \in [0,\delta_{0}) K_{X}^{\mu}]]
\setminus (Sing X
\cup_{x \in [0,\delta_{0}) K_{X}^{\mu}}
\cup_{P \in  TI_{\mu X}^{\beta, \rho}(x)} \Gamma(\mu \lambda^{m_{\beta}} X_{\beta,M}, P, I_{\beta})) . \]
Since the elements of $Reg(\mu X_{\beta}(\lambda),\rho)$ and
$Reg_{\infty}(\mu X_{\beta}(\lambda),\rho)$ depend continuously on
$\lambda \in K_{X}^{\mu} \subset {\mathbb S}^{1} \setminus B_{\beta}^{\mu}(X)$
so they do the elements of $Reg^{*}(I_{\beta}, \mu X, K_{X}^{\mu})$ by continuity of the flow.
Thus $\alpha_{\beta}^{\mu X}$ and $\omega_{\beta}^{\mu X}$ are constant by restriction
to $H \in Reg^{*}(I_{\beta}, \mu X, K_{X}^{\mu})$.
%
%
The dynamics of $Re(\mu X)$
is a topological product in the intermediate sets when we avoid
the directions in $B_{X}^{\mu}$. Such a property is also true in the exterior
sets. We want to assemble the information
attached to the exterior and intermediate sets to describe
the behavior of $Re(\mu X)$ in $|y| \leq \epsilon$.

\begin{lem}
\label{lem:goins}
Let $X \in \Xt$. Fix $\mu \in {\mathbb S}^{1}$.
Let $P_{0}
\in [0,\delta_{0}) \times K_{X}^{\mu} \times \partial{B(0,\epsilon)}$
such that $Re(\mu X)$ does not point towards
${\mathbb C} \setminus \overline{B}(0,\epsilon)$ at $P_{0}$.
Then $[0,\infty)$ is contained in
$It(\mu X, P_{0},\overline{B}(0,\epsilon))$
and $\lim_{\zeta \to \infty} {\rm exp}(\zeta \mu X)(P_{0}) \in Sing X$.
Moreover the intersection of ${\rm exp}((0,\infty) \mu X)(P_{0})$ with
every intermediate of exterior set is connected.
\end{lem}
The last property is important. Once the trajectories of $Re(\mu X)$ enter
a set $M_{\beta}$ or $T_{\beta}$ they never go out.
\begin{proof}
Denote $P_{0}=(r_{0},\lambda_{0},y_{0})$.
We can suppose $r_{0} \neq 0$ and $N(X)>1$.
Otherwise the result is a consequence of
corollary \ref{cor:npee} since
$\{r_{0} \lambda_{0} \} \times \overline{B}(0,\epsilon) \subset E_{0}$.
Since $E_{0}(r_{0},\lambda_{0}) \cap Sing X = \emptyset$
then
$c_{0} = \sup It(\mu \lambda_{0}^{d_{0}} X_{0,E},P_{0}, {E}_{0})$
belongs to ${\mathbb R}^{+}$. Denote
$Q_{0}={\rm exp}(c_{0} \mu \lambda_{0}^{d_{0}} X_{0,E})(P_{0})$.
We have that $Q_{0} \in \partial M_{0}$ and $Re(\mu X)$
points towards $\dot{I}_{0}$ at $Q_{0}$.
The point $Q_{0}$ is contained in the closure of a unique
$H$ in $Reg^{*}(I_{0}, \mu X,K_{X}^{\mu})$. Moreover
$(\alpha_{0}^{\mu X}, \omega_{0}^{\mu X})(H) = (\infty, \zeta)$ for some $\zeta \in {\mathbb C}$
since ${\mathcal G}_{\mu X_{0}(\lambda_{0})}$ is connected.
We have that $d_{1} =\sup(It(\mu X,Q,I_{0}))$ belongs to ${\mathbb R}^{+}$.
Denote $P_{1}={\rm exp}(d_{1} \mu X)(Q_{0})$; we obtain
$P_{1} \in \partial I_{0} \cap \partial E_{0\zeta}$.
Moreover $Re(\mu X)$ points towards $\dot{E}_{0\zeta}$ at
$P_{1}$. Denote $\beta(0)=0$ and $\beta(1)=0 \zeta$.
Analogously there exist sequences $\beta(0)$, $\hdots$, $\beta(k)$ and
\[ (P_{0},0) = (P_{0},d_{0}) , \ (Q_{0},c_{0}), \ (P_{1},d_{1}), \ (Q_{1},c_{1}),
\ \hdots, \ (P_{k},d_{k}) \  k \geq 1 \]
such that we have
${\rm exp}((d_{l},c_{l}) \mu X)(P_{0}) \subset \dot{E}_{\beta(l)}$,
${\rm exp}((c_{j-1},d_{j}) \mu X)(P_{0}) \subset \dot{I}_{\beta(j-1)}$,
\[ Q_{l} ={\rm exp}(c_{l} \mu X)(P_{0}) \in \partial E_{\beta(l)} \cap \partial
M_{\beta(l)}, \
P_{j} = {\rm exp}(d_{j} \mu X)(P_{0}) \in
\partial I_{\beta(j-1)} \cap \partial E_{\beta(j)} \]
for all $1 \leq j, l+1 \leq k$
and $E_{\beta(k)}=T_{\beta(k)}$. By corollary \ref{cor:npee} and
the discussion in subsection \ref{subsub:npes} then
$Re(\mu X)$ points towards $\dot{E}_{\beta(k)}$
at $\partial E_{\beta(k)}$
and $P_{k}$ is in the basin of attraction of
$Sing X \cap E_{\beta(k)}$. Thus we get
${\rm exp}((0,\infty) \mu X)(P_{k}) \subset \dot{E}_{\beta(k)}$ and
$\lim_{z \to \infty}{\rm exp}(z \mu X)(P_{0}) \in Sing X \cap E_{\beta(k)}$.
\end{proof}
We define $Reg^{*}(\epsilon, \mu X, K_{X}^{\mu})$ the set of connected
components of
\[  [[|y| < \epsilon] \cap [x \in [0,\delta_{0}) K_{X}^{\mu}]]
\setminus (Sing X
\cup_{x \in [0,\delta_{0}) K_{X}^{\mu}}
\cup_{P \in  T_{\mu X}^{\epsilon}(x)} \Gamma(\mu X, P, |y| \leq \epsilon)) . \]
We define $\alpha^{\mu X}(P)=\lim_{z \to - \infty} {\rm exp}(z \mu X)(P)$
for all $P \in [|y| \leq \epsilon]$ such that $It(\mu X,P,|y| \leq \epsilon)$
contains $(- \infty,0)$. Otherwise we define $\alpha^{\mu X}(P)=\infty$.
We define $\omega^{\mu X}(P)$ in an analogous way.

Given $H \in Reg^{*}(\epsilon, \mu X, K_{X}^{\mu})$ the functions
$(\alpha^{\mu X})_{|H}$ and $(\omega^{\mu X})_{|H}$ satisfy
that either they are identically $\infty$ or their value is never
$\infty$. Since the basins of attraction and repulsion of the curves
in $Sing_{V} \mu X$ in $x \in [0,\delta_{0}) K_{X}^{\mu}$ are open sets
then $(\alpha^{\mu X})_{|H}$ and $(\omega^{\mu X})_{|H}$
are continuous. Thus
$(\alpha^{\mu X})_{|H(x)}$ and $(\omega^{\mu X})_{|H(x)}$
are constant for all $x \in [0,\delta_{0}) K_{X}^{\mu}$.
Indeed we can interpret $\alpha^{\mu X}(H)$ either as $\infty$
if $(\alpha^{\mu X})_{|H} \equiv \infty$ or as the element
of $Sing_{V} X$ that contains $\alpha^{\mu X}(H)$ otherwise. Denote
\[ Reg_{\infty}(\epsilon, \mu X, K_{X}^{\mu})=
Reg^{*}(\epsilon,\mu X, K_{X}^{\mu})
\cap ((\alpha^{\mu X})^{-1}(\infty) \cup (\omega^{\mu X})^{-1}(\infty)) \]
and $Reg(\epsilon, \mu X, K_{X}^{\mu}) =
Reg^{*}(\epsilon, \mu X, K_{X}^{\mu}) \setminus
Reg_{\infty}(\epsilon, \mu X, K_{X}^{\mu})$.
We define
\[ Reg_{j}(\epsilon, \mu X, K_{X}^{\mu}) =
\{ H \in Reg(\epsilon, \mu X, K_{X}^{\mu}) :
\sharp \{\alpha^{\mu X}(H), \omega^{\mu X}(H) \} =j \} \]
for $j \in \{1,2\}$. We have that the set $H(x)$ is connected for
$H \in Reg(\epsilon, \mu X,K_{X}^{\mu})$ and $x \in (0,\delta_{0})K_{X}^{\mu}$.
The set $H(0)$ is connected for
$H \not \in Reg_{2}(\epsilon, \mu X,K_{X}^{\mu})$
whereas otherwise $H(0)$ has two connected components.

We define an oriented graph ${\mathcal G}(\mu X, K_{X}^{\mu})$.
The set of vertexes is $Sing_{V} X$ whereas the edges
are the elements of $Reg(\epsilon, \mu X, K_{X}^{\mu})$.
The edge $H \in Reg(\epsilon, \mu X, K_{X}^{\mu})$ joins the vertexes
$\alpha^{\mu X}(H)$ and $\omega^{\mu X}(H)$, we denote
$\alpha^{\mu X}(H) \stackrel{H}{\to} \omega^{\mu X}(H)$.
The graph  ${\mathcal NG}(\mu X, K_{X}^{\mu})$ is obtained from
${\mathcal G}(\mu X, K_{X}^{\mu})$ by removing the reflexive edges
and the orientation of edges.
\begin{pro}
\label{pro:acyconn}
Let $X \in \Xt$. Fix $\mu \in {\mathbb S}^{1}$ and a compact connected
set $K_{X}^{\mu} \subset {\mathbb S}^{1} \setminus B_{X}^{\mu}$.
Then the graph ${\mathcal NG}(\mu X, K_{X}^{\mu})$ is acyclic and connected.
\end{pro}
We say that a exterior set $E_{\beta}$ has depth $0$ if $N(\beta)=1$.
In general given $E_{\beta}$ such that $N(\beta)>1$ we define
$depth(E_{\beta})=1 + \sup_{\zeta \in S_{\beta}} depth(E_{\beta \zeta})$.
\begin{proof}
An exterior set $E_{\beta}=[\eta \geq |t| \geq \rho |x|]$ is contained
in $T_{\beta}=[\eta \geq |t|]$. We can associate graphs
${\mathcal G}_{\beta}(\mu X, K_{X}^{\mu})$
and ${\mathcal NG}_{\beta}(\mu X, K_{X}^{\mu})$ to
the vector field $Re(\mu \lambda^{d_{\beta}} X_{\beta, E})$
defined in $T_{\beta}$.

Consider an exterior set $E_{\beta}$ such that $depth(E_{\beta})=0$.
The graph  ${\mathcal NG}_{\beta}(\mu X, K_{X}^{\mu})$ has only
one vertex and no edges, therefore it is connected and acyclic.

Suppose that ${\mathcal NG}_{\beta}(\mu X, K_{X}^{\mu})$ is connected
and acyclic for all exterior set $E_{\beta}$ such that
$depth(E_{\beta}) \leq k$. It is enough to prove that the result is true
for every exterior set $E_{\beta}$ such that $depth(E_{\beta})=k+1$.

Fix $\lambda_{0} \in K_{X}^{\mu}$. The graph
${\mathcal NG}_{\mu X_{\beta}(\lambda_{0})}$ is connected and acyclic
by corollary \ref{cor:nocycles} and proposition \ref{pro:connected}.
Consider an edge $J_{0} \in Reg(\mu X_{\beta}(\lambda_{0}))$ of the graph
${\mathcal G}_{\mu X_{\beta}(\lambda_{0})}$ joining the vertexes
$\zeta(1)$ and $\zeta(2)$. We denote also by $J_{0}$
the component of $Reg(\mu X_{\beta}(\lambda_{0}),\rho)$
associated to $J_{0}$ by proposition \ref{pro:gofidi} where
$M_{\beta}=[|w| \leq \rho]$. Let $J_{1}$ be the element of
$Reg^{*}(I_{\beta}, \mu X, K_{X}^{\mu})$ such that
$J_{1}(0,\lambda_{0}) \subset J_{0}$.
By lemma \ref{lem:goins} applied to
$Re(\mu \lambda^{d_{\beta \zeta(1)}} X_{\beta \zeta(1),E})$ in
$T_{\beta \zeta(1)}$ we deduce that
$\alpha^{\mu X}(J_{1}) \subset Sing X$. By the open character of the singular
points in $(r,\lambda) \in [0,\delta_{0}) \times K_{X}^{\mu}$ we obtain that
$\alpha^{\mu X}(J_{1})$ is contained in an irreducible component
$\gamma_{1}$ of $Sing X$. Analogously $\omega^{\mu X}(J_{1})$ is contained
in an irreducible component $\gamma_{2}$ of $Sing X$.
Denote by $J_{2}$ the edge of ${\mathcal NG}_{\beta}(\mu X, K_{X}^{\mu})$
joining $\gamma_{1}$ and $\gamma_{2}$.

The set ${\mathbb C} \setminus J_{0}$ has two connected components
$H_{1} \ni \zeta(1)$ and $H_{2} \ni \zeta(2)$ (lemma \ref{lem:reg2}). Denote
$Sg_{j}=H_{j} \cap S_{\beta}$ for $j \in \{1,2\}$.
We obtain that there is no edge different than $J_{2}$ of
${\mathcal G}_{\beta}(\mu X, K_{X}^{\mu})$ joining
a vertex of ${\mathcal NG}_{\beta \upsilon}(\mu X, K_{X}^{\mu})$
and a vertex of ${\mathcal NG}_{\beta \kappa}(\mu X, K_{X}^{\mu})$
for $\upsilon \in Sg_{1}$ and $\kappa \in Sg_{2}$. Moreover the
restriction of ${\mathcal G}_{\beta}(\mu X, K_{X}^{\mu})$ to
$Sing_{V} X_{\beta \upsilon, E}$ is
${\mathcal G}_{\beta \upsilon}(\mu X, K_{X}^{\mu})$
for all $\upsilon \in S_{\beta}$. Then the acyclicity of every
${\mathcal NG}_{\beta \upsilon}(\mu X, K_{X}^{\mu})$
for all $\upsilon \in S_{\beta}$ imply
that ${\mathcal NG}_{\beta}(\mu X, K_{X}^{\mu})$ is acyclic.
Finally, since
${\mathcal NG}_{\mu X_{\beta}(\lambda_{0})}$ and
${\mathcal NG}_{\beta \upsilon}(\mu X, K_{X}^{\mu})$
are connected for all $\upsilon \in S_{\beta}$ then
${\mathcal NG}_{\beta}(\mu X, K_{X}^{\mu})$ is connected.
\end{proof}
The properties of ${\mathcal G}(\mu X, K_{X}^{\mu})$ are inherited
of the properties of the polynomial vector fields associated to $X$.
\begin{pro}
Let $X \in \Xt$. Fix $\mu \in {\mathbb S}^{1}$ and a compact connected
set $K_{X}^{\mu} \subset {\mathbb S}^{1} \setminus B_{X}^{\mu}$.
Then we have
\[ \sharp (Reg(\epsilon, \mu X, K_{X}^{\mu}) \cap
(\alpha^{\mu X}, \omega^{\mu X})^{-1}(\gamma, \gamma))
= 2 \nu_{X}(\gamma) \]
for all $\gamma \in Sing_{V} X$.
Moreover we have
$\sharp Reg(\epsilon, \mu X, K_{X}^{\mu}) = 2 \nu(X) - N(X) +1$.
\end{pro}
\subsection{Analyzing the regions}
Let $X \in \Xt$.  Fix $\mu \in e^{i(0,\pi)}$ and a compact connected
set $K_{X}^{\mu} \subset {\mathbb S}^{1} \setminus B_{X}^{\mu}$.
Consider a region $H \in Reg_{1}(\epsilon, \mu X,K_{X}^{\mu})$.
We denote by $T_{\mu X,H}^{\epsilon}(x)$ the unique tangent
point in $T_{\mu X}^{\epsilon}(x) \cap \overline{H(x)}$ for
all $x \in [0,\delta_{0})K_{X}^{\mu}$.
Let $\psi$ be an integral of the time form of $X$ defined in
a neighborhood of $T_{\mu X, H}^{\epsilon}(0)$. By analytic continuation
we obtain an integral of the time form
$\psi_{H,L}^{X}=\psi_{H,R}^{X}$ of $X$ in $H=H_{L}=H_{R}$ such that it is
holomorphic in $H \setminus [x=0]$ and continuous in $H$.
Moreover $(x,\psi_{H,L}^{X})=(x,\psi_{H,R}^{X})$ is injective in $H$
since $\psi_{H,L}^{X}(H(x))$ is simply connected for all
$x \in [0,\delta_{0}) K_{X}^{\mu}$.

Let $H \in Reg_{2}(\epsilon, \mu X,K_{X}^{\mu})$.
Let $L_{\mu X, H}^{\epsilon}(x)$ be the point in
$T_{\mu X}^{\epsilon}(x) \cap \overline{H(x)}$ such that
$Re(X)$ points towards $H$ for all $x \in [0,\delta_{0})K_{X}^{\mu}$.
We define $H_{L}(0)$ the connected component of $H(0)$ such that
$L_{\mu X,H}^{\epsilon}(0) \in \overline{H_{L}(0)}$.
We denote by $R_{\mu X,H}^{\epsilon}(x)$ the other point in
$T_{\mu X}^{\epsilon}(x) \cap \overline{H(x)}$ for
$x \in [0,\delta_{0})K_{X}^{\mu}$.
We define $H_{L}= H_{L}(0) \cup (H \setminus [x=0])$ and
$H_{R} = H \setminus H_{L}(0)$.
Let $\psi_{\kappa}$ be a holomorphic integral of the time form of
$X$ defined in a neighborhood of $\kappa_{\mu X,H}^{\epsilon}(0)$
for $\kappa \in \{L,R\}$. We obtain an integral $\psi_{H,\kappa}^{X}$
of the time form of $X$ in $H_{\kappa}$ obtained by analytic continuation
of $\psi_{\kappa}$ for $\kappa \in \{L,R\}$.
The function $\psi_{H,\kappa}^{X}$ is  holomorphic in
$H \setminus [x=0]$ and continuous in $H_{\kappa}$ for $\kappa \in \{L,R\}$.
Moreover $(x,\psi_{H,L}^{X})$ and $(x,\psi_{H,R}^{X})$ are injective in $H_{L}$
and $H_{R}$ respectively. The theorem of the residues implies that
\[ \psi_{H,L}^{X}(x,y) - \psi_{H,R}^{X}(x,y) - 2 \pi i
\sum_{P \in J(x)} Res(X,P) \]
is bounded in $H \setminus [x=0]$ where
$J(x)$ is the subset of $(Sing X)(x)$ of points contained in the
same connected component of $B(0,\epsilon) \setminus H(x)$ than
$\omega^{\mu X}(H(x))$. Since $H(0)$ is disconnected the function
$x \to \sum_{P \in J(x)} Res(X,P)$ is not bounded in
$x \in (0,\delta_{0})K_{X}^{\mu}$. Indeed
$x \to \sum_{P \in J(x)} Res(X,P)$ can be extended to a pure
meromorphic function defined in a neighborhood of $x=0$.

%
%

 We call subregion of a region
$H \in Reg(\epsilon, \mu X, K_{X}^{\mu})$ to every set of the form
$H \cap E_{\beta}$ or $H \cap I_{\beta}$ where $E_{\beta}$ is an
exterior set and $I_{\beta}$ is an intermediate set.
We say that all the subregions of
$H \in  Reg_{1}(\epsilon, \mu X, K_{X}^{\mu})$
are both $L$-subregions and $R$-subregions.
Consider $H \in Reg_{2}(\epsilon, \mu X, K_{X}^{\mu})$. There exists a
magnifying glass set $M_{\beta(0)}$ such that the curves
$\alpha^{\mu X}(H)$ and $\omega^{\mu X}(H)$ are contained in
$M_{\beta(0)}$
but they are in different connected components of
$M_{\beta(0)} \setminus I_{\beta(0)}$.
A subregion of $H$ contained in $M_{\beta(0)}$ is both a $L$-subregion
and an $R$-subregion. A subregion in the same connected component
of $\overline{H \setminus M_{\beta(0)}}$ than $L_{\mu X, H}^{\epsilon}$
is called a L-subregion. A subregion of $H$ in the same connected component
of $\overline{H \setminus M_{\beta(0)}}$ than $R_{\mu X, H}^{\epsilon}$
is called a R-subregion. We define
${H}^{L}$ the union of the
L-subregions of $H$ whereas ${H}^{R}$
is the union of the R-subregions of $H$. Clearly we have $H=H^{L} \cup H^{R}$
by lemma \ref{lem:goins}.
\section{Extension of the Fatou coordinates}
\label{sec:Fatou}
A diffeomorphism $\varphi \in \diff{tp1}{2}$ is a small deformation
of its convergent normal form ${\rm exp}(X)$ in suitable domains.
The dynamical splitting associated to $X$ provides
information about the dynamics of $\varphi$. That is going to lead us to
define the analogue of the Ecalle-Voronin invariants. For such a purpose we
need to mesure the ``distance'' from ${\rm exp}(X)$ to $\varphi$.
\subsection{Comparing $\varphi \in \diff{tp1}{2}$ and a convergent normal
form}
Let ${\rm exp}(X)$ be a convergent normal form of $\varphi$. We consider
$\sigma_{z}(x,y) = (x, y + z (y \circ \varphi -y ))$
for $z \in B(0,2)$. Let $\psi$ be an integral of the time form
of $X$, i.e. $X(\psi)=1$. We define
$\Delta_{\varphi}= \psi \circ \sigma_{1}(P)-(\psi(P)+1)$ for
$P \not \in Fix \varphi$ in a neighborhood of
$(0,0)$ as follows: we choose a determination
$\psi_{|x=x(P)}$ in the neighborhood of $P$, we define
$\psi \circ \sigma_{1}(P)$ as the evaluation at $\sigma_{1}(P)$
of the analytic continuation of $\psi_{|x=x(P)}$ along the path
$\gamma:[0,1] \to [x=x(P)]$ given by $\gamma(z)=\sigma_{z}(P)$.
The value of $\Delta_{\varphi}$ does not depend on the determination
of $\psi$ that we chose. Clearly $\Delta_{\varphi}$ is holomorphic in
$U \setminus Fix \varphi$ for some neighborhood $U$ of $(0,0)$.
Indeed we have:
\begin{lem}
Let $\varphi \in \diff{tp1}{2}$ (with fixed convergent normal form).
Then the function $\Delta_{\varphi}$ belongs to the ideal
$(y \circ \varphi -y)$ of the ring ${\mathbb C}\{x,y\}$.
\end{lem}
The result is a consequence of Taylor's formula applied to
\[ \Delta_{\varphi} = \psi \circ \varphi - \psi \circ {\rm exp}(X)
\sim (\partial \psi / \partial y) \circ {\rm exp}(X)
(y \circ \varphi - y \circ {\rm exp}(X)) = O (y \circ \varphi - y) . \]
\begin{pro}
\label{pro:bddconf}
Let $\varphi \in \diff{tp1}{2}$ with fixed convergent normal form
${\rm exp}(X)$. Fix $\mu \in {\mathbb S}^{1}$ and a compact connected
set $K_{X}^{\mu} \subset {\mathbb S}^{1}$. Consider
$H \in Reg(\epsilon, \mu X, K_{X}^{\mu})$. Then we have
\[ \Delta_{\varphi} = O(X(y)) = O \left({
\frac{1}{(1+|\psi_{H,\kappa}^{X}|)^{1+1/\nu(\varphi)}}
}\right)   \]
in $H^{\kappa}$ for all $\kappa \in \{L,R\}$.
\end{pro}
\begin{proof}
Denote $f=X(y)$.
Let us prove the result for a $L$-subregion $J$ without lack of generality.
Lemma \ref{lem:goins} implies that
there exists a sequence $B(0)$, $\hdots$, $B(k)=J$ of $L$-subregions
of $H$ such that
\begin{itemize}
\item $B(2j) \subset E_{\beta(2j)}$ for all $0 \leq 2j \leq k$.
\item $B(2j+1) \subset I_{\beta(2j)}$ for all $0 \leq 2j+1 \leq k$.
\item $\beta(0)=0$ and $\beta(2j+2) = \beta(2j) \upsilon(j)$ for
some $\upsilon(j) \in {\mathbb C}$ and all $0 \leq 2j+2 \leq k$.
\end{itemize}
Denote $K(2j)=E_{\beta(2j)}$, $h(2j)=d_{\beta(2j)}$,
$K(2j+1)=I_{\beta(2j)}$ and $h(2j+1)=m_{\beta(2j)}$.
Denote $\partial_{e} B(0)= [|y| \leq \epsilon] \cap B(0)$ and
$\partial_{e} B(j) = \overline{B(j)} \cap \partial{K(j-1)}$ for $j \geq 1$.
We define the property $Pr(j)$ as
\[
Pr(j): \left\{{
\begin{array}{l}
\sup |(\psi_{H,L}^{X})_{|\partial_{e} B(j)}| \leq
M_{j}/|x|^{h(j)} \ {\rm for \ some} \ M_{j} \in {\mathbb R}^{+}
\ {\rm if} \ j \leq k \\
f = O (1/(1+|\psi_{H,L}^{X}|)^{1+1/\nu(\varphi)})
\ {\rm in} \ B(0) \cup \hdots \cup B(j-1).
\end{array} }\right.
\]
We have that $Pr(k+1)$ implies the result in the proposition for $J$.
The result is true for $j=0$. It is enough to prove that
$Pr(j) \implies Pr(j+1)$ for all $0 \leq j \leq k$.

  From the construction of the splitting we obtain that
$f \in ({x}^{h(j)+ [(j+1)/2]})$
in $K(j)$ for all $0 \leq j \leq k$ (let us remark that $[(j+1)/2]$
is the integer part of $(j+1)/2$). Denote $Y=(X/x^{h(j)})_{|K(j)}$.
There exists a holomorphic integral
$\psi_{j}$ of the time form of $Y$ in a neighborhood of the
simply connected set $\overline{B(j)}$ such that $|\psi_{j}| \leq M_{j}'$ in
$\partial_{e} B(j)$ for some $M_{j}'>0$. Suppose
that $K(j)=[\eta \geq |t| \geq \rho |x|]$ is a parabolic exterior set,
since $\nu(Y) \leq \nu(X)$ we obtain
\begin{equation}
\label{equ:equsou}
 f = O \left({
\frac{x^{h(j) +[(j+1)/2]}}{(1+|\psi_{j}|)^{1 + 1/\nu(Y)}}
}\right) = O \left({
\frac{x^{h(j) +[(j+1)/2]}}{(1+|\psi_{j}|)^{1 + 1/\nu(X)}} }\right)
\end{equation}
by remark \ref{rem:bditfpet}.
The inequality  $|x^{h(j)} \psi_{H,L}^{X} - \psi_{j}| \leq M_{j} + M_{j}'$
implies
\[ f = O \left({
\frac{x^{h(j) +[(j+1)/2] -h(j) (1+1/\nu(X))}}
{(1+|\psi_{H,L}^{X}|)^{1 + 1/\nu(X)}}
}\right) = O \left({
\frac{1}{(1+|\psi_{H,L}^{X}|)^{1 + 1/\nu(X)}} }\right)  \]
in $B(j)$ since $h(j) \leq [(j+1)/2] \nu(X)$ by construction. Moreover
if $j < k$ then we have $|\psi_{j}|=O(1/|x|^{\nu(Y)})$ in
$\partial_{e} B(j+1)$. We deduce that there exists $M_{j+1}>0$
such that
$|\psi_{H,L}^{X}| \leq M_{j+1}/|x|^{h(j+1)}$ in
$\partial_{e} B(j+1)$ since $h(j+1)=h(j)+\nu(Y)$,

Suppose that $K(j)=[\eta \geq |t|]$ is a non-parabolic exterior set,
this implies $j=k$. We have that
$\psi_{k}(r,\lambda,t) \lambda^{- d_{\beta(k)}} \mu^{-1} -
C(r,\lambda) \ln (t-\gamma(x))$
is bounded in $B(k)$ where  $t=\gamma(x)$ is the only irreducible
component of $Sing X_{\beta(j),E}$
by the discussion in subsection \ref{subsub:npes}.
There exists $\upsilon>0$ such that
$\arg(C(r,\lambda))$ in $(-\pi/2+\upsilon,\pi/2-\upsilon)$
for all $(r,\lambda) \in [0,\delta_{0}) \times K_{X}^{\mu}$ if
$B(k)$ is a basin of repulsion, otherwise we have that
$\arg(C(r,\lambda)) \in (\pi/2+\upsilon,3 \pi/2 - \upsilon)$
for all $(r,\lambda) \in [0,\delta_{0}) \times K_{X}^{\mu}$. We deduce that
\[ f=O(x^{h(k) + [(k+1)/2]} (t-\gamma(x))) =
O(x^{h(k) + [(k+1)/2]} e^{-K |\psi_{k}|}) \]
in $B(j)$ for some $K>0$. This implies equation
\ref{equ:equsou} and then $Pr(k+1)$.

Finally suppose that $K(j)$ is an intermediate set. We have that
$\psi_{j}$ is bounded in $B(j)$. Thus there exists $M_{j+1}>0$
such that $|\psi_{H,L}^{X}| \leq M_{j+1}/|x|^{h(j)}=
 M_{j+1}/|x|^{h(j+1)}$ in $\overline{B(j)}$ and then in
$\partial_{e}(B(j+1))$. Moreover $f= O(x^{h(j) +[(j+1)/2]})$
implies equation \ref{equ:equsou} and then $Pr(j+1)$.
\end{proof}
\subsection{Constructing Fatou coordinates}
Let $\varphi \in \diff{tp1}{2}$ with convergent normal form
$\alpha={\rm exp}(X)$. Fix $\mu =i {e}^{i \theta}$
with $\theta \in (-\pi/2,\pi/2)$ and a
compact connected set
$K_{X}^{\mu} \subset {\mathbb S}^{1} \setminus B_{X}^{\mu}$. Consider
$H \in Reg(\epsilon, \mu X, K_{X}^{\mu})$.
Let $P \in H$, suppose $P \in H^{L}$ without lack of generality.
The trajectory $\Gamma=\Gamma(\mu X,P,|y| \leq \epsilon)$
is contained in $H^{L}$.
Let $B(P)$ be the strip ${\rm exp}([0,1]X)(\Gamma)$ and
$\Gamma'=\alpha(\Gamma)$. The distance between the lines
$\psi_{H,L}^{X}(\Gamma)$ and $\psi_{H,L}^{X}(\Gamma')$ is $\cos \theta$.
Since $\psi_{H,L}^{X} \circ \varphi = \psi_{H,L}^{X} \circ \alpha +
\Delta_{\varphi}$ then $\Gamma$ and
$\varphi(\Gamma)$ enclose a strip $B_{1}(P)$ whenever
$sup_{B(0,\delta_{0}) \times B(0,\epsilon)} |\Delta_{\varphi}| <
(\cos \theta) / 3$.
Since $\Delta_{\varphi}(0,0)=0$ this condition is fulfilled  by taking
$\mu$ away from $-1$ and $1$ and a small neighborhood
$B(0,\delta_{0}) \times B(0,\epsilon)$ of $(0,0)$.

Let $\tilde{B}(P)$ be the complex space obtained from $B(P)$ by identifying
$\Gamma$ and $\Gamma'$. Let $\tilde{B}_{1}(P)$ be the complex space obtained
from $B_{1}(P)$ by identifying $\Gamma$ and $\varphi(\Gamma)$.
The space $\tilde{B}(P)$ is biholomorphic to ${\mathbb C}^{*}$
by $e^{2 \pi i z} \circ \psi_{H,L}^{X}$. A natural compactification
$\overline{B}(P)$ is obtained by adding $0 \sim \omega^{\mu X}(P)$ and
$\infty \sim \alpha^{\mu X}(P)$. Analogously we will obtain a biholomorphism
from $\overline{B}_{1}(P)$ to
$\pn{1}$. The space of orbits of $\varphi_{|H_{L}(x(P))}$ is then
rigid, that will allow us to define analytic invariants of $\varphi$.
Let us remark that $\tilde{B}_{1}(P)$ is the restriction
of the space of orbits of $\varphi$ to $H_{L}(x(P))$ for all choices of
$H$ and $P$ if and only if
$\nu_{X}(\gamma) \geq 1$ for all $\gamma \in Sing_{V} X$ \cite{topology}.
In general the complete space of orbits is messier, we obtain
further identifications via return maps.

We consider the coordinates given by $\psi_{H,L}^{X}$. We define
\[ \sigma_{0}(z) =
z + \eta(\cos \theta Re((z-\psi_{H,L}^{X}(P)) e^{-i \theta}))
\Delta_{\varphi} \circ {\alpha}^{\circ (-1)}
\circ (x, \psi_{H,L}^{X})^{\circ (-1)}(x(P), z) \]
where $\eta: {\mathbb R} \to [0,1]$ is a $C^{\infty}$ function
such that $\eta(b)=0$ for all $b \leq 1/3$ and $\eta(b)=1$ for
all $b \geq 2/3$. This definition implies that
$\sigma = (\psi_{H,L}^{X})^{\circ (-1)} \circ \sigma_{0} \circ
\psi_{H,L}^{X}$ satisfies
\[ \sigma_{{\rm exp}([0,1/3]X)(\Gamma)} \equiv Id \ \ {\rm and} \ \
\sigma_{{\rm exp}([-1/3,0]X)(\Gamma')} \equiv
\varphi \circ {\alpha}^{\circ (-1)}. \]
The mappings $\sigma_{0}$ and $\sigma$ depend on the choice of
the base point $P$.
The function $\Delta_{\varphi} \circ \alpha^{\circ (-1)} \circ
(\psi_{H,L}^{X})^{\circ (-1)}$ is holomorphic. By Cauchy's integral
formula we obtain
\[ \frac{\partial( \Delta_{\varphi} \circ \alpha^{\circ (-1)} \circ
(\psi_{H,L}^{X})^{\circ (-1)})}{\partial z}(z_{0}) =
\frac{1}{2 \pi i} \int_{|z-z_{0}|=1}
\frac{\Delta_{\varphi} \circ \alpha^{\circ (-1)} \circ
(\psi_{H,L}^{X})^{\circ (-1)}}{(z-z_{0})^{2}} dz . \]
Denote $B^{2}(P) = {\rm exp}(\overline{B}(0,2)X)(B(P))$.
By prop. \ref{pro:bddconf} there exists $C>1$ such that
\begin{equation}
\label{equ:boujac}
\left|{
 \frac{\partial( \Delta_{\varphi} \circ \alpha^{\circ (-1)} \circ
(\psi_{H,L}^{X})^{\circ (-1)})}{\partial z}
}\right| (z)  \leq C \min \left({
\frac{1}{(1+|z|)^{1+1/\nu(\varphi)}} ,
\sup_{B^{2}(P)} |\Delta_{\varphi}| }\right)
\end{equation}
for all $z \in \psi_{H,L}^{X}(B(P))$. The jacobian matrix
${\mathcal J} \sigma_{0}$ of $\sigma_{0}$ is a $2 \times 2$ real matrix.
The coefficients of ${\mathcal J} \sigma_{0} - Id$ are bounded
by an expression like the one in the right hand side of equation
\ref{equ:boujac}, maybe for a bigger $C>1$.
We obtain that $\sup_{B(0,\delta_{0}) \times B(0,\epsilon)}|\Delta_{\varphi}|$
small implies  ${\mathcal J} \sigma_{0} \sim Id$ and then $\sigma$ is
a $C^{\infty}$ diffeomorphism from $\tilde{B}(P)$ onto $\tilde{B}_{1}(P)$.

The mapping $\xi=e^{2 \pi i z} \circ \psi_{H,L}^{X} \circ \sigma^{\circ (-1)}$
is a $C^{\infty}$ diffeomorphism from $\tilde{B}_{1}(P)$ onto
${\mathbb C}^{*}$.
The function $\psi_{H,L}^{X} \circ \sigma^{\circ (-1)}$ is a Fatou coordinate,
even if not holomorphic in general, of $\varphi$ in $B_{1}(P)$.
The complex dilatation $\chi_{\sigma_{0}}$ of $\sigma_{0}$ satisfies
\[ |\chi_{\sigma_{0}}|(z) =
\frac{
\left|{
\frac{\partial \sigma_{0}}{\partial \overline{z}}
}\right|
}{\left|{ \frac{\partial \sigma_{0}}{\partial z}
}\right|} (z)
\leq  K(H) \min \left({
\frac{1}{(1+|z|)^{1+1/ \nu(\varphi)}} ,
\sup_{{\rm exp}(\overline{B}(0,2)X)(H(x(P)))} |\Delta_{\varphi}| }\right) \]
for all $z \in \psi_{H,L}^{X}(B(P))$ and some $K(H)>1$ independent of
$P \in H$.
Since $\xi^{\circ (-1)}$ is equal to
$(\psi_{H,L}^{X})^{\circ (-1)} \circ \sigma_{0} \circ ((1/2\pi i) \ln z)$ then
\begin{lem}
\label{lem:norinf}
We have
\[ |\chi_{\xi^{\circ (-1)}}|(z) \leq K(H) \min \left({
\frac{1}{(1+2^{-1} \pi^{-1}|\ln z|)^{1+1/ \nu(\varphi)}} ,
\sup_{{\rm exp}(\overline{B}(0,2)X)(H(x(P)))} |\Delta_{\varphi}| }\right) \]
for all $z \in {e}^{2 \pi i w} \circ \psi_{H,L}^{X}(B(P))$.
\end{lem}
The mapping $\xi$ and then $\chi_{\xi^{\circ (-1)}}$ depend on the
base point $P$.
We look for a quasi-conformal mapping $\tilde{\rho}: \pn{1} \to \pn{1}$
such that $\chi_{\tilde{\rho}}=\chi_{\xi^{\circ (-1)}}$. Since we can
suppose $||\chi_{\xi^{\circ (-1)}}||_{\infty}=
\sup_{{\mathbb C}^{*}} |\chi_{\xi^{\circ (-1)}}| <1/2<1$
then such a mapping
exists by the Ahlfors-Bers theorem. The choice of
$\tilde{\rho}$ is unique if $\tilde{\rho}$ fulfills $\tilde{\rho}(0)=0$,
$\tilde{\rho}(1)=1$ and $\tilde{\rho}(\infty)=\infty$.
By construction $\tilde{\rho} \circ \xi$ is a biholomorphism
from $\tilde{B}_{1}(P)$ to ${\mathbb C}^{*}$.

We define
\[ J(r) = \frac{2}{\pi} \int_{|z| < r}
\frac{K(H)}{(1+ 2^{-1} \pi^{-1} |\ln |z||)^{1+1/\nu(\varphi)}}
\frac{1}{|z|^{2}} d \sigma  \]
for $r \in {\mathbb R}^{+}$. We have that $J(r) < \infty$ for
all $r \in {\mathbb R}^{+}$.
\begin{lem}
\label{lem:Le-Vi}
The mapping $\tilde{\rho}$ is conformal at $0$ and at $\infty$.
Moreover we have
\[ \left|{ \frac{\tilde{\rho}(z)}{z} -
\frac{\partial \tilde{\rho}}{\partial z}(0) }\right| \leq \left|{
\frac{\partial \tilde{\rho}}{\partial z}(0)
}\right| \iota(|z|) \ {\rm and} \
 \left|{ \frac{z}{\tilde{\rho}(z)} -
{\frac{\partial \tilde{\rho}}{\partial z}(\infty) }^{-1}
}\right| \leq {\left|{
\frac{\partial \tilde{\rho}}{\partial z}(\infty)
}\right|}^{-1} \iota(1/|z|) \]
where $\iota$ depends on $K(H)$, it satisfies
$\lim_{|z| \to 0} \iota(|z|)=0$. We have
\[ min_{|z|=1} |\tilde{\rho}(z)| e^{-J(1)} \leq
|\partial \tilde{\rho}/\partial z|(0), |\partial \tilde{\rho}/
\partial z|(\infty) \leq  max_{|z|=1} |\tilde{\rho}(z)| e^{J(1)} . \]
\end{lem}
\begin{proof}
We define
\[ I(r) = \frac{1}{\pi}
\int_{|z|<r} \frac{1}{1-|\chi_{\tilde{\rho}}|}
\frac{|\chi_{\tilde{\rho}}(z)|}{|z|^{2}} d \sigma \]
for all $r \in {\mathbb R}^{+}$. We have $I(r) \leq J(r)$ for
all $r \in {\mathbb R}^{+}$.
To get the conformality of $\tilde{\rho}$ at $z=0$ it is enough to
prove that $I(r) < \infty$ for all $r \in {\mathbb R}^{+}$
(theorem 6.1 in page 232 of \cite{Le-Vi}). This is clear since
$J(r) < \infty$ for all $r \in {\mathbb R}^{+}$.
The inequality is obtained for
a function $\iota$ such that $\lim_{|z| \to 0} \iota(|z|)=0$, it depends on
$J$ and then on $K(H)$. The proof for $z=\infty$
is obtained by applying the result in \cite{Le-Vi}
to $1/\tilde{\rho}(1/z)$.
\end{proof}
We denote by $[z_{0},z_{1}]$ the spherical distance for
$z_{0},z_{1} \in \pn{1}$.
\begin{lem} (\cite{Ah-Be}, lemma 17, page 398).
\label{lem:ah-be}
Let $\chi$ be a measurable complex-valued function in $\pn{1}$.
Suppose $||\chi||_{\infty} <1$. Then
there exists a unique quasi-conformal mapping
$\upsilon: \pn{1} \to \pn{1}$ such that $\chi_{\upsilon}=\chi$,
$\upsilon(0)=0$, $\upsilon(1)=1$, $\upsilon(\infty)=\infty$ and
$[\upsilon(z),z] \leq C_{0} ||\chi||_{\infty}$
for all $z \in \pn{1}$ and some $C_{0}>0$ not depending on $\chi$.
\end{lem}
\begin{cor}
\label{cor:sphere}
$[\tilde{\rho}(z),z] \leq C_{0} ||\chi_{\xi^{\circ (-1)}}||_{\infty}$
for all $z \in \pn{1}$
\end{cor}
We define $\rho=\tilde{\rho}/(\partial \tilde{\rho}/\partial z)(0)$.
The quasi-conformal mapping $\rho$ is the only solution of
$\chi_{\rho} = \chi_{\xi^{\circ (-1)}}$ such that
$\rho(0)=0$, $\rho(\infty)=\infty$ and $(\partial \rho/\partial z)(0)=1$.
\begin{lem}
\label{lem:codeinf}
$\lim_{||\chi_{\xi^{\circ (-1)}}||_{\infty} \to 0}
(\partial \tilde{\rho} / \partial z)(z_{0}) = 1$ for $z_{0} \in \{0,\infty\}$.
In particular we have $\lim_{||\chi_{\xi^{\circ (-1)}}||_{\infty} \to 0}
(\partial \rho / \partial z)(\infty) = 1$.
\end{lem}
\begin{proof}
Denote $\chi = \chi_{\xi^{\circ (-1)}}$.
For $||\chi||_{\infty}$ small enough there exists $C_{1}>0$ such that
$|\tilde{\rho}(z)-z| \leq C_{1} ||\chi||_{\infty}$
for all $z \in \overline{B}(0,1)$ by corollary \ref{cor:sphere}.
This leads us to
\[ \left|{ \frac{\partial{\tilde{\rho}}}{\partial z}(0)-1 }\right|
\leq (1+C_{1}) e^{J(1)} \iota(|z|) + \frac{C_{1}}{|z|} ||\chi||_{\infty} \]
for all $z \in \overline{B}(0,1) \setminus \{0\}$
(lemma \ref{lem:Le-Vi}).
By evaluating at $z=\sqrt{||\chi}||_{\infty}$ we obtain that
$\lim_{||\chi||_{\infty} \to 0}
(\partial \tilde{\rho} / \partial z)(0) = 1$.
Analogously we get
$\lim_{||\chi||_{\infty} \to 0}
(\partial \tilde{\rho} / \partial z)(\infty) = 1$.
Since we have $\rho=\tilde{\rho}/(\partial \tilde{\rho}/\partial z)(0)$
then $\lim_{||\chi||_{\infty} \to 0}
(\partial \rho / \partial z)(\infty) = 1$.
\end{proof}
\begin{lem}
\label{lem:rhsi1}
$\lim_{||\chi_{\xi^{\circ (-1)}}||_{\infty} \to 0}
\sup_{z \in \pn{1}} |\rho(z)/z-1| = 0$.
\end{lem}
\begin{proof}
Denote $\chi=\chi_{\xi^{\circ (-1)}}$.
Let $b>0$. By lemma \ref{lem:Le-Vi} there exists
$r_{0} \in {\mathbb R}^{+}$ such that
$|\rho(z)/z-1| < b$ for all $z \in B(0,r_{0})$. We also obtain
\[ \left|{ \frac{z}{\rho(z)} - {\frac{\partial \rho}{\partial z}(\infty)}^{-1}
}\right| \leq \left|{ {\frac{\partial \rho}{\partial z}(\infty)}^{-1} }\right|
\iota (1/|z|) . \]
Since $\lim_{||\chi||_{\infty} \to 0}
(\partial \rho/\partial z)(\infty)=1$
then there exist $a_{0}>0$ and $r_{1}>0$ such that
$|\rho(z)/z-1| < b$ for all $z \in {\mathbb C} \setminus B(0,r_{1})$
if $||\chi||_{\infty} < a_{0}$. There exists $a_{1}>0$ and $C_{1}>0$
such that $|\tilde{\rho}(z)-z| < C_{1} ||\chi||_{\infty}$ for all
$z \in \overline{B}(0,r_{1}) \setminus B(0,r_{0})$ if
$||\chi||_{\infty} < a_{1}$. We deduce that
\[ \left|{
\frac{\rho(z)}{z} -1
}\right| \leq \left|{
1 - 1/(\partial \tilde{\rho}/\partial z)(0)
}\right| + \frac{C_{1} ||\chi||_{\infty}}
{|\partial \tilde{\rho}/\partial z|(0) |z|} \]
for all $z \in \overline{B}(0,r_{1}) \setminus B(0,r_{0})$ and
$||\chi||_{\infty} < a_{1}$. By lemma \ref{lem:codeinf} there exists
$a \in {\mathbb R}^{+}$ such that $|\rho(z)/z-1| < b$ for all
$z \in \pn{1}$ if $||\chi||_{\infty} < a$.
\end{proof}
Now we can define the function
\[ \psi_{H,L,P}^{\varphi} = \frac{1}{2 \pi i} \ln z \circ
\rho \circ {e}^{2 \pi i z} \circ \psi_{H,L}^{X} \circ \sigma^{\circ (-1)} . \]
It is an injective Fatou coordinate of $\varphi$ in the
neighborhood of $B_{1}(P)$. By using
$\psi_{H,L,P}^{\varphi} \circ \varphi = \psi_{H,L,P}^{\varphi}+1$
we can extend $\psi_{H,L,P}^{\varphi}$ to $H_{L}(x(P))$.

It looks like $\psi_{H,L,P}^{\varphi}$ depends on the choice of the
base point $P \in H^{L}$. Nevertheless the functions $\psi_{H,L,P}^{\varphi}$
paste together to provide a Fatou coordinate $\psi_{H,L}^{\varphi}$, it is
continuous in $H_{L}$ and holomorphic in $\dot{H}$.
\begin{lem}
\label{lem:bdddif}
Denote $\xi_{0}= e^{2 \pi i z}  \circ \psi_{H,L}^{X}$.
There exists $C>0$ independent of $P \in H^{L}$ such that
\[ \left|{ \psi_{H,L,P}^{\varphi} - \psi_{H,L}^{X}}\right| \leq
\frac{1}{\pi}  {\left|{ \left|{ \frac{\rho}{z} -1
}\right| }\right|}_{\infty} +
\frac{C}{{(1+|\psi_{H,L}^{X}|)}^{1+1/\nu(\varphi)}} \]
in $B_{1}(P)$. Moreover we have
\[ \lim_{Z \in B_{1}(Q), \xi_{0}(Z) \to z_{0}}
\psi_{H,L,P}^{\varphi}(Z) - \psi_{H,L}^{X}(Z) =
\frac{1}{2 \pi i} \ln \frac{\partial \rho}{\partial z}(z_{0}) \]
for all $z_{0} \in \{0, \infty \}$ and all $Q \in H_{L}(x(P))$.
\end{lem}
\begin{proof}
Denote $\chi=\chi_{\xi^{\circ (-1)}}$ and $\kappa=\rho/z-1$.
We have $\lim_{||\chi||_{\infty} \to 0} ||\kappa||_{\infty}=0$
(lemma \ref{lem:rhsi1}). Thus we obtain
\[ \left|{
\psi_{H,L,P}^{\varphi} - \psi_{H,L}^{X} \circ \sigma^{\circ (-1)}
}\right| =\frac{1}{2 \pi} \left|{ \ln (1+\kappa(z)) \circ
{e}^{2 \pi i z} \circ \psi_{H,L}^{X} \circ \sigma^{\circ (-1)}
}\right| \leq
\frac{||\kappa||_{\infty}}{\pi}   \]
for $||\kappa||_{\infty}$ small enough.
On the other hand we get
\[ \left|{
\psi_{H,L}^{X} \circ \sigma^{\circ (-1)} - \psi_{H,L}^{X} }\right|
=  \left|{
\sigma_{0}^{\circ (-1)} \circ  \psi_{H,L}^{X} - \psi_{H,L}^{X} }\right|
\leq \frac{C}{{(1+|\psi_{H,L}^{X}|)}^{1+1/\nu(\varphi)}} \]
for some $C>0$ and all $z \in B_{1}(P)$ (prop. \ref{pro:bddconf}).
Analogously we obtain
\[ \lim_{Z \in B_{1}(P), \xi_{0}(Z) \to z_{0}}
\psi_{H,L,P}^{\varphi}(Z) - \psi_{H,L}^{X}(Z) =
\frac{1}{2 \pi i} \ln \frac{\partial \rho}{\partial z}(z_{0}) \]
for $z_{0} \in \{0,\infty\}$. We can suppose
$\sup_{B(0,\delta_{0}) \times B(0,\epsilon)} |\Delta_{\varphi}| < 1/2$.
As a consequence given $Q \in H_{L}(x(P))$ there exists
$k(Q) \in {\mathbb N}$ such that every $Z \in B_{1}(Q)$
is of the form $\varphi^{\circ (j(Z))}(P')$ for some
$P' \in B_{1}(P)$ and $j(Z) \in [-k(Q),k(Q)]$.
Moreover if $j(Z) \geq 0$ then $\varphi^{\circ (l)}(P') \in H_{L}(x(P))$
for $0 \leq l < j(Z)$ whereas for $j(Z) <0$ we have that
$\varphi^{\circ (-l)}(P') \in H_{L}(x(P))$
for $0 \leq l < -j(Z)$.

Fix $Q \in H_{L}(x(P))$. Consider $Z \in B_{1}(Q)$, we can suppose
$j(Z)>0$ without lack of generality. This leads us to
\[ \psi_{H,L,P}^{\varphi}(Z) - \psi_{H,L}^{X}(Z)
= (\psi_{H,L,P}^{\varphi}(P') - \psi_{H,L}^{X}(P')) -
\sum_{l=0}^{j(Z)-1} \Delta_{\varphi} \circ \varphi^{\circ (l)}(P') . \]
Since
$|\psi_{H,L}^{X}(\varphi^{\circ (l)}(P')) -
\psi_{H,L}^{X}(Z) + j(Z) -l| < k(Q)/2$ for all $0 \leq l \leq j(Z)$ then
\[ \left|{
\sum_{l=0}^{j(Z)-1} \Delta_{\varphi} \circ \varphi^{\circ (l)}(P')
}\right| \leq
\frac{k(Q) C}{{(1- k(Q)/2 + |Img(\psi_{H,L}^{X}(Z))|)}^{1+1/\nu(\varphi)}} \]
Now $\xi_{0}(Z) \to 0,\infty$ implies
$|Img (\psi_{H,L}^{X}(Z))| \to \infty$. Thus we obtain
\[ \lim_{Z \in B_{1}(P), \xi_{0}(Z) \to z_{0}}
\psi_{H,L,P}^{\varphi}(Z) - \psi_{H,L}^{X}(Z) =
\lim_{Z \in B_{1}(Q), \xi_{0}(Z) \to z_{0}}
\psi_{H,L,P}^{\varphi}(Z) - \psi_{H,L}^{X}(Z) \]
for $z_{0} \in \{0,\infty\}$.
\end{proof}
We prove next that $\psi_{H,L,P}$ depends only on $x(P)$.
\begin{lem}
Let $x_{0} \in [0,\delta_{0})K_{X}^{\mu}$. We have
$\psi_{H,L,P}^{\varphi} \equiv \psi_{H,L,Q}^{\varphi}$ in
$H_{L}(x_{0})$ for all $P,Q \in H^{L}(x_{0})$.
We also have
$\psi_{H,L,P}^{\varphi} - \psi_{H,L}^{X} \equiv
\psi_{H,R,Q}^{\varphi}  - \psi_{H,R}^{X}$ if $x_{0} \neq 0$
and $(P,Q) \in H^{L}(x_{0}) \times H^{R}(x_{0})$. Then
$(\partial \rho/\partial z)(\infty)$ depends only on $H$ and $x(P)$.
\end{lem}
\begin{proof}
Let $P,Q \in H^{L}(x_{0})$. We have
$\psi_{H,L,P}^{\varphi} - \psi_{H,L,Q}^{\varphi} \in
\vartheta(\tilde{B}_{1}(P))$ since
\[ (\psi_{H,L,P}^{\varphi} - \psi_{H,L,Q}^{\varphi}) \circ \varphi
\equiv \psi_{H,L,P}^{\varphi} - \psi_{H,L,Q}^{\varphi} .\]
We define
$h =(\psi_{H,L,P}^{\varphi} - \psi_{H,L,Q}^{\varphi})
\circ {(\psi_{H,L,P}^{\varphi})}^{\circ (-1)}
\circ 1/(2 \pi i) \ln z$ in ${\mathbb C}^{*}$. The function
extends to a holomorphic function in $\pn{1}$ such that
$h(0)=0$ by lemma \ref{lem:bdddif}. Therefore
we obtain $h \equiv 0$ and then
$\psi_{H,L,P}^{\varphi} \equiv \psi_{H,L,Q}^{\varphi}$.

We have $(\psi_{H,L}^{X} - \psi_{H,R}^{X})(x_{0},y) \equiv b(x_{0})$
in $H(x_{0})$ for some $b(x_{0}) \in {\mathbb C}$. We define
$g =(\psi_{H,L,P}^{\varphi} - \psi_{H,R,Q}^{\varphi})
\circ {(\psi_{H,L,P}^{\varphi})}^{\circ (-1)}
\circ 1/(2 \pi i) \ln z$ in
${\mathbb C}^{*}=  (e^{2 \pi i z} \circ \psi_{H,L,P})(H(x_{0}))$.
By lemma \ref{lem:bdddif}
the complex function $g$ admits a continuous extension to $\pn{1}$
such that $g(0)=b(x_{0})$.
We are done since then $g \equiv b(x_{0})$.
\end{proof}
Here it is important the choice $\rho(0)=0$, $\rho(\infty)=\infty$,
$\rho'(0)=1$. By replacing $\rho$ by the canonical choice
$\tilde{\rho}(0)=0$, $\tilde{\rho}(1)=1$, $\tilde{\rho}(\infty)=\infty$
in the definition of $\psi_{H,L,P}^{\varphi}$ we would have
$\psi_{H,L,P}^{\varphi} \not \equiv \psi_{H,L,Q}^{\varphi}$ in general.

Denote by $\psi_{H,L}^{\varphi}$ any of the functions
$\psi_{H,L,P}^{\varphi}$ defined in $H_{L}$. The definition of
$\psi_{H,R}^{\varphi}$ is analogous.
We denote by $\psi_{H}^{\varphi} - \psi_{H}^{X}$ the function
defined in $H$ which is given by the expression
$\psi_{H,l}^{\varphi}-\psi_{H,l}^{X}$ in $H_{l}$ for $l \in \{L,R\}$.
The definitions of $\psi_{H,L}^{\varphi}$, $\psi_{H,R}^{\varphi}$
and $\psi_{H}^{\varphi}-\psi_{H}^{X}$ allow to deduce asymptotic
properties of those functions when approaching the fixed points without
checking out that they are stable by iteration.
\begin{pro}
Let $\varphi \in \diff{tp1}{2}$ with fixed convergent normal form
${\rm exp}(X)$. Fix $\mu \in e^{i(0,\pi)}$ and a compact connected
set $K_{X}^{\mu} \subset {\mathbb S}^{1} \setminus B_{X}^{\mu}$. Let
$H \in Reg(\epsilon, \mu X, K_{X}^{\mu})$; the mappings
$(x,\psi_{H,L}^{\varphi})$ and
$(x,\psi_{H,R}^{\varphi})$ are holomorphic in $\dot{H}$ and continuous
and injective in $H_{L}$ and $H_{R}$ respectively.
\end{pro}
\begin{proof}
Consider $P=(x_{0},y_{0}) \in H^{L}$.
The mapping $\sigma_{0}(x,z)$ depends holomorphically on $x$.
There exists a continuous section $P(x_{1}) \in [x=x_{1}]$ for
$x_{1}$ in a neighborhood $V$ of $x_{0}$ in  $[0,\delta_{0}) K_{X}^{\mu}$
such that $\psi_{H,L}^{X}(P(x_{1}))=\psi_{H,L}^{X}(P)$ and $P(x_{0})=P$.
The mapping $\sigma = \psi_{H,L}^{X} \circ \sigma_{0} \circ
(\psi_{H,L}^{X})^{\circ (-1)}$ maps $B(P(x))$ onto
$B_{1}(P(x))$ and establishes a $C^{\infty}$ diffeomorphism
from $\tilde{B}(P(x))$ onto $\tilde{B}_{1}(P(x))$ for all
$x \in V$. The complex dilation $\chi_{\xi^{\circ (-1)}}$
depends holomorphically on $x \in \dot{V}$
and continuously on $x \in V$.
Hence the dependance of the canonical solution
$\tilde{\rho}:\pn{1} \to \pn{1}$ of
$\chi_{\tilde{\rho}} = \chi_{\xi^{\circ (-1)}}$ with respect to
$x$ is continuous in $V$ and holomorphic in $\dot{V}$.
In particular the function
$x \to (\partial \tilde{\rho}/\partial z)(x,0)$ is holomorphic in
$\dot{V}$ and continuous in $V$.
We deduce that
$\rho(x,z) = \tilde{\rho}(x,z)/(\partial \tilde{\rho} / \partial z)(x,0)$
depends continuously on  $x \in V$ and holomorphically on $x \in \dot{V}$.
Then $\psi_{H,L}^{\varphi}$ is continuous
in $\cup_{x \in V} B_{1}(P(x))$ and holomorphic in
a neighborhood of  $\cup_{x \in \dot{V}} B_{1}(P(x))$. Since
$P$ can be any point of $H^{L}$ then $\psi_{H,L}^{\varphi}$ is holomorphic
in $\dot{H}$ and continuous in $H_{L}$. Moreover
$(x,\psi_{H,\upsilon}^{\varphi})$ is injective in
$H_{\upsilon}$ for $\upsilon \in \{L,R\}$ since $\psi_{H,\upsilon}^{\varphi}$
is injective in the fundamental domains of type $B_{1}(P)$.
\end{proof}
\begin{cor}
\label{cor:lifipo}
Let $\varphi \in \diff{tp1}{2}$ with fixed convergent normal form
${\rm exp}(X)$. Fix $\mu \in e^{i(0,\pi)}$ and a compact connected
set $K_{X}^{\mu} \subset {\mathbb S}^{1} \setminus B_{X}^{\mu}$. Let
$H \in Reg(\epsilon, \mu X, K_{X}^{\mu})$.
The function  $x \to (\partial \rho / \partial z)(H,x,\infty)$
is well-defined and continuous in $[0,\delta_{0}) K_{X}^{\mu}$.
It is holomorphic in $(0,\delta_{0})\dot{K}_{X}^{\mu}$ and
$(\partial \rho / \partial z)(H,0,\infty)=1$. Moreover
we have $(\partial \rho / \partial z)(H,x,\infty) \equiv 1$ if
$H \in Reg_{1}(\epsilon, \mu X,K_{X}^{\mu})$.
\end{cor}
\begin{proof}
By the proof of the previous proposition we have that
$x \to (\partial \tilde{\rho} / \partial z)(x,0)$
and $x \to (\partial \tilde{\rho} / \partial z)(x,\infty)$ are
continuous in $[0,\delta_{0})K_{X}^{\mu}$ and holomorphic
in $(0,\delta_{0}) \dot{K}_{X}^{\mu}$. The same property is clearly
fulfilled by $x \to (\partial \rho / \partial z)(x,\infty)$.

Consider $P = {\rm exp}(s X)(L_{\mu X}^{H}(0))$
if $H \in Reg_{2}(\epsilon, \mu X,K_{X}^{\mu})$ for all
$s \in {\mathbb R}^{+}$.
For $H \in Reg_{1}(\epsilon, \mu X,K_{X}^{\mu})$ consider
$P  = {\rm exp}(s X)(T_{\mu X}^{H}(0))$ for $s \in {\mathbb R}^{+}$
if $Re(- i \mu X)$ points towards $H$ at $T_{\mu X}^{H}(0)$, otherwise
we denote $P  = {\rm exp}(-s X)(T_{\mu X}^{H}(0))$ for
$s \in {\mathbb R}^{+}$. Then $P$ is well-defined and belongs to
$H^{L}(0)=H_{L}(0)$ for all $s \in {\mathbb R}^{+}$. Moreover
$\inf_{Q \in B(P)} |\psi_{H,L}^{X}(Q)|$ tends to $\infty$ when
$s \to \infty$. We obtain  $||\chi_{\xi^{\circ (-1)}}||_{\infty} \to 0$
when $s \to 0$ by lemma \ref{lem:norinf}. This
implies $(\partial \rho / \partial z)(0,\infty)=1$ by lemma \ref{lem:codeinf}.
The prove of $(\partial \rho / \partial z)(x,\infty) \equiv 1$
in the case $H \in Reg_{1}(\epsilon, \mu X,K_{X}^{\mu})$ is analogous.
\end{proof}
\begin{pro}
\label{pro:bddcon}
Let $\varphi \in \diff{tp1}{2}$ with fixed convergent normal form
${\rm exp}(X)$. Fix $\mu \in e^{i(0,\pi)}$ and a compact connected set
$K_{X}^{\mu} \subset {\mathbb S}^{1} \setminus B_{X}^{\mu}$. Let
$H \in Reg(\epsilon, \mu X, K_{X}^{\mu})$;
the function $\psi_{H}^{\varphi} - \psi_{H}^{X}$ is continuous in
$H \cup [Fix \varphi \cap \partial H]$.
\end{pro}
\begin{proof}
The function $\psi_{H}^{\varphi} - \psi_{H}^{X}$ is
clearly continuous in $H$. We define
\[ (\psi_{H}^{\varphi} - \psi_{H}^{X})(\alpha^{\mu X}(H(x)))
= \frac{1}{2 \pi i} \ln \frac{\partial \rho}{\partial z}(H, x,\infty)
\ \ {\rm and} \ \
(\psi_{H}^{\varphi} - \psi_{H}^{X})(\omega^{\mu X}(H(x)))=0 \]
for all $x \in [0,\delta_{0})K_{X}^{\mu}$ where $\ln 1=0$.
The function
$(\psi_{H}^{\varphi} - \psi_{H}^{X})_{|Fix \varphi \cap \partial H}$
is continuous by corollary \ref{cor:lifipo}. Let $P \in H^{L}$. From
\[ \psi_{H,L,P}^{\varphi}-\psi_{H,L}^{X} =
(\psi_{H,L,P}^{\varphi} - \psi_{H,L}^{X} \circ \sigma^{\circ (-1)})
+ (\psi_{H,L}^{X} \circ \sigma^{\circ (-1)} - \psi_{H,L}^{X}) \]
we deduce that
\[ \left|{ \psi_{H}^{\varphi}-\psi_{H}^{X} -
\frac{1}{2 \pi i}  \left({ \ln  \left({ \frac{\rho}{z} }\right) \circ
{e}^{2 \pi i z} \circ \psi_{H,L}^{X} \circ \sigma^{\circ (-1)}
}\right) }\right| \leq
\frac{C}{{(1+|\psi_{H,L}^{X}|)}^{1+1/\nu(\varphi)}}   \]
in $B_{1}(P)$ for some $C>0$ independent of $P \in H^{L}$. We can suppose that
the function $\iota$ provided by lemma \ref{lem:Le-Vi} is increasing.
By varying $P$ we obtain
\[ |\psi_{H}^{\varphi} - \psi_{H}^{X}| \leq
\frac{1}{\pi} \iota \circ e^{z}(-\pi Img (\psi_{H,L}^{X})) +
\frac{C}{{(1+|\psi_{H,L}^{X}|)}^{1+1/\nu(\varphi)}} \]
in $H^{L} \cap [Img \psi_{H,L}^{X} > J_{0}]$ for some $J_{0}>0$.
An analogous expression can be obtained in $H^{R}$
by replacing $\psi_{H,L}^{X}$ with $\psi_{H,R}^{X}$.
Lemma \ref{lem:Le-Vi} implies the existence of an increasing
$\iota'$ independent of $P \in H^{L}$ such that
$\lim_{|z| \to 0} \iota'(|z|)=0$ and
$|\rho/z (\partial \rho/\partial z)(\infty)^{-1} - 1| \leq \iota'(1/|z|)$.
We deduce that
\[ \left|{
\psi_{H}^{\varphi} - \psi_{H}^{X} - \frac{1}{2 \pi i}
\ln \frac{\partial \rho}{\partial z} (x,\infty) }\right|
\leq \frac{1}{\pi} \iota' \circ e^{z}(\pi Img (\psi_{H,L}^{X})) +
\frac{C}{{(1+|\psi_{H,L}^{X}|)}^{1+1/\nu(\varphi)}} \]
in $H^{L} \cap [Img \psi_{H,L}^{X} < -J_{1}]$ for some $J_{1}>0$.
Again an analogous expression can be obtained for $H^{R}$. We deduce that
\[ \lim_{(x,y) \in H^{\kappa}, \ |Img(\psi_{H,\kappa}^{X}(x,y))| \to \infty, \
(x,y) \to (x_{0},y_{0})} (\psi_{H}^{\varphi}-\psi_{H}^{X})(x,y) =
(\psi_{H}^{\varphi}-\psi_{H}^{X})(x_{0},y_{0})  \]
for $\kappa \in \{L,R\}$ and all
$(x_{0},y_{0}) \in Fix \varphi \cap \partial H$. This implies the
continuity at $(x_{0},y_{0})$
except when $H \in Reg_{1}(\epsilon,\mu X, K_{X}^{\mu})$ or $x_{0} = 0$.
In particular we can suppose
$(\psi_{H}^{\varphi}-\psi_{H}^{X})(x_{0},y_{0})=0$.
We have
\[ \lim_{(x,y) \in H_{\kappa},\ (x,y) \to (x_{0},y_{0})}
|\psi_{H,\kappa}^{X}(x,y)| = \infty \ \ \forall \kappa \in \{L,R\} . \]
It is enough to prove that
$(\psi_{H}^{\varphi}-\psi_{H}^{X}
)_{(H^{\kappa} \cap [|Img \psi_{H,\kappa}^{X}| < D]) \cup \{ (x_{0},y_{0} )\}}$
is continuous at $(x_{0},y_{0})$ for all $D \in {\mathbb R}^{+}$.
Suppose $\kappa=L$ without lack of generality. There exists a function
$\upsilon_{D}:{\mathbb R}^{+} \to {\mathbb R}_{\geq 0}$
such that $\lim_{b \to \infty} \upsilon_{D}(b)=\infty$ and holding that
given $P$ in $H^{L} \cap [|Img \psi_{H,L}^{X}| < D]$ then $B(P)$ is
contained in $[|\psi_{H,L}^{X}| > \upsilon_{D}(|Re (\psi_{H,L}^{X}(P))|)]$.
The value $||\rho/z-1||_{\infty}$ tends to $0$ when
$||\chi_{\xi^{\circ (-1)}}||_{|\infty}$ by lemma \ref{lem:rhsi1}.
Moreover we have $||\chi_{\xi^{\circ (-1)}}||_{|\infty} \leq
K(H)/(1+ \upsilon_{D} (|Re (\psi_{H,L}^{X}(P))|))^{1+1/\nu(\varphi)}$
by lemma \ref{lem:norinf}. The lemma \ref{lem:bdddif} implies that
\[ \lim_{(x,y) \in H^{L} \cap  [|Img \psi_{H,L}^{X}| < D]
, (x,y) \to (x_{0},y_{0})} (\psi_{H}^{\varphi}-\psi_{H}^{X})(x,y)=0=
(\psi_{H}^{\varphi}-\psi_{H}^{X})(x_{0},y_{0}) \]
and then the result is proved.
\end{proof}
The previous proposition implies that by considering a smaller domain
of definition $|y| \leq \epsilon$ we can suppose that
$sup_{Q \in H} |\psi_{H}^{\varphi}-\psi_{H}^{X}|(Q)$ is as small as
desired for all $H \in Reg(\epsilon, \mu X, K_{X}^{\mu})$ since
$(\psi_{H}^{\varphi}-\psi_{H}^{X})(0,0)=0$.
\begin{cor}
\label{cor:Lavvf}
Let $\varphi \in \diff{tp1}{2}$ with fixed convergent normal form
${\rm exp}(X)$. Fix $\mu \in e^{i(0,\pi)}$ and a compact connected set
$K_{X}^{\mu} \subset {\mathbb S}^{1} \setminus B_{X}^{\mu}$. Let
$H \in Reg(\epsilon, \mu X, K_{X}^{\mu})$. There exists a unique
vector field $X_{H}^{\varphi}=X_{H}^{\varphi}(y) \partial/\partial y$
(continuous in $H$ and holomorphic in $\dot{H})$ such that
$X_{H}^{\varphi}(\psi_{H}^{\varphi}) \equiv 1$.
Moreover $X_{H}^{\varphi}(y)/X(y)-1$ is a continuous function in
$H \cup [\overline{H} \cap Fix \varphi]$
vanishing at $\overline{H} \cap Fix \varphi$.
\end{cor}
\begin{rem}
The Lavaurs vector field $X_{H}^{\varphi}$
has asymptotic development $\log \varphi$
until the first non-zero term in the neighborhood of the fixed points.
That is a consequence of
$(\log \varphi)(y) - X(y) \in (y \circ \varphi -y)^{2}$
and the previous corollary.
\end{rem}
\begin{rem}
\label{rem:smalld}
The construction of $\psi_{H,L,P}^{\varphi}$ or $\psi_{H,R,P}^{\varphi}$
in $B_{1}(P)$ depends only on getting small values of
$||\chi_{\xi^{\circ (-1)}}||_{\infty}$. This condition
is automatically fulfilled for
$\sup_{B(0,\delta) \times B(0,\epsilon)} |\Delta_{\varphi}|$ small enough
(lemma \ref{lem:norinf}).
\end{rem}
\section{Defining the analytic invariants}
\label{sec:defanai}
Now we define an extension of the Ecalle-Voronin
invariants for $\varphi \in \diff{tp1}{2}$. It is the key to prove
the main theorems in this paper.
\subsection{Normalizing the Fatou coordinates}
Let $\varphi \in \diff{tp1}{2}$ with fixed convergent normal form
${\rm exp}(X)$. Fix $\mu \in e^{i(0,\pi)}$ and a compact connected
set $K_{X}^{\mu}$ contained in ${\mathbb S}^{1} \setminus B_{X}^{\mu}$.
There are $2 \nu(\varphi)$ continuous sections $T_{X}^{\epsilon,1}$,
$\hdots$, $T_{X}^{\epsilon, 2 \nu(\varphi)}$ of the set
$T_{X}^{\epsilon}$.
We will always suppose that $T_{X}^{\epsilon,1}$, $\hdots$,
$T_{X}^{\epsilon, 2 \nu(\varphi)}$,
$T_{X}^{\epsilon, 2 \nu(\varphi)+1}=T_{X}^{\epsilon,1}$
are ordered in counter clock-wise sense. For all
$j \in {\mathbb Z}/(2 \nu(\varphi) {\mathbb Z})$ there exists a function
$\theta_{j}: B(0,\delta) \to {\mathbb R}^{+}$ such that
\[ T_{X}^{\epsilon, j+1}(x) =
T_{X}^{\epsilon,j}(x) e^{i \theta_{j}(x)} \ \ {\rm and} \ \
T_{X}^{\epsilon, j}(x) e^{i (0, \theta_{j}(x))} \cap
T_{X}^{\epsilon}(x)=\emptyset  \ \ \forall x \in B(0,\delta). \]
There exists a unique $T_{\mu X}^{\epsilon,j}(x)$ in
$T_{X}^{\epsilon,j}(x) e^{i (0,\theta_{j}(x))}$.
Denote by $\upsilon_{j}(x)$ the only value in $(0,2 \pi)$ such that
$T_{\mu X}^{\epsilon,j+1}(x) = T_{\mu X}^{\epsilon,j}(x)
e^{i \upsilon_{j}(x)}$.
We define $H(j)$ as the element of $Reg(\epsilon, \mu X, K_{X}^{\mu})$
such that $T_{\mu X}^{\epsilon,j}(x) \in \partial H(j)(x)$ for all
$x \in [0,\delta_{0})K_{X}^{\mu}$.
We define $H(j)_{S}=H(j)_{L}$ if $Re(X)$ points towards
$H(j)$ at $T_{\mu X}^{\epsilon, j}(0)$, otherwise we define
$H(j)_{S}=H(j)_{R}$. The region
$H \in Reg_{k}(\epsilon, \mu X, K_{X}^{\mu})$ appears $k$ times in the
sequence $H(1)$, $\hdots$, $H(2 \nu(\varphi))$.
We denote by $H_{\infty}(j)$ the element
of $Reg_{\infty}(\epsilon, \mu X,K_{X}^{\mu})$ such that
$T_{X}^{\epsilon, j+1}(x)$ belongs to $\partial (H_{\infty}(j)(x))$
for all $x \in [0,\delta_{0}) K_{X}^{\mu}$.
%
%

We define the function
$\zeta_{\varphi}(x) = - \pi i \nu(\varphi)^{-1}
\sum_{P \in (Fix \varphi)(x)} Res(\varphi,P)$. It is holomorphic in
a neighborhood of $0$.
Fix $j_{0} \in \{1, \hdots, 2\nu(\varphi) \}$. Consider an integral
$\psi_{j_{0}}^{X}$ of the time form of $X$ defined in the neighborhood
of $T_{\mu X}^{\epsilon,j_{0}}(0)$. We can extend it to
$H(j_{0})_{S}$ by
analytic continuation. In an analogous way we can define
$\psi_{j_{0}+k}^{X}$ in $H(j_{0} + k)_{S}$ for all $k \in {\mathbb Z}$;
we choose $\psi_{j_{0}+k}^{X}(T_{\mu X}^{\epsilon,j_{0}+k}(0))$ to be
the result of evaluating the analytic extension of
$\psi_{j_{0}}^{X} + k \zeta_{\varphi}$ along the curve
$t \to T_{\mu X}^{\epsilon,j_{0}}(0) e^{i t \kappa}$ for $t \in [0,1]$ where
$\kappa = \sum_{l=0}^{k-1} \upsilon_{j_{0}+l}(0)$ if $k>0$ and
$\kappa = - \sum_{l=1}^{-k} \upsilon_{j_{0}-l}(0)$ for $k<0$.
If $Re(X)$ points towards $H(j)$ at
$T_{\mu X}^{\epsilon, j}(0)$
then we define $\psi_{H(j),L}^{X}=\psi_{j}^{X}$, otherwise
we define $\psi_{H(j),R}^{X}=\psi_{j}^{X}$. By construction we have
$\psi_{j + 2\nu(\varphi)}^{X} \equiv \psi_{j}^{X}$ for $j \in {\mathbb Z}$.

We choose an element $\gamma_{1} \equiv (y=\alpha_{1}(x))$ of $Sing_{V} X$,
we call $\gamma_{1}$ the {\it privileged curve} associated to $X$
(or $\varphi$). We have
$X = u(x,y) \prod_{j=1}^{N} (y-\alpha_{j}(x))^{n_{j}} \partial / \partial y$
for some unit $u \in {\mathbb C}\{x,y\}$. Denote by $\gamma_{j}$
the curve $y=\alpha_{j}(x)$ for $2 \leq j \leq N(\varphi)$.
We look for functions $c_{1}$, $\hdots$, $c_{N}$ contained in
$C^{0}([0,\delta_{0})K_{X}^{\mu}) \cap \vartheta ((0,\delta_{0})
\dot{K}_{X}^{\mu})$ such that
\begin{itemize}
\item $c_{1} \equiv 0$
\item Given $\gamma_{j} \stackrel{H}{\to} \gamma_{k}$
of ${\mathcal G}(\mu X, K_{X}^{\mu})$ then
$(c_{j}-c_{k})(x) \equiv 1/(2 \pi i)
\ln (\partial \rho/\partial z)(H,x,\infty)$.
\end{itemize}
By cor. \ref{cor:lifipo}
the reflexive edges of ${\mathcal G}(\mu X, K_{X}^{\mu})$ do not
impose any restriction. There is a unique solution
$c_{1}$, $\hdots$, $c_{N}$ since ${\mathcal NG}(\mu X, K_{X}^{\mu})$
is connected and acyclic. We say that $c_{1}, \hdots, c_{N}$ is a sequence of
{\it privileged functions} associated to $(X,\varphi,K_{X}^{\mu}, \gamma_{1})$.

Denote $\gamma_{k(j)}=\omega^{\mu X}(H(j))$.
We define a Fatou coordinate $\psi_{j}^{\varphi}$ of
$\varphi$ in the set $H(j)_{S}$ given by
$\psi_{j}^{\varphi}(x,y) = \psi_{j}^{X}(x,y) +c_{k(j)}(x)$
We obtain that
$(\psi_{j}^{\varphi}-\psi_{j}^{X})_{|\gamma_{k}} \equiv c_{k}$
for $\gamma_{k} \in \{ \alpha^{\mu X}(H(j)), \omega^{\mu X}(H(j)) \}$.
Since given $\psi_{j}^{\varphi}$ the
function $\psi_{j}^{\varphi} + c(x)$ is also a Fatou coordinate we
normalize by fixing a privileged curve and the
sequence of privileged functions attached to such a choice.
\subsection{Defining the changes of charts}
Our aim is to define
\[ \xi_{\varphi,K_{X}^{\mu}}^{j}(x,z)=\psi_{j+1}^{\varphi} \circ
(x,\psi_{j}^{\varphi})^{\circ (-1)} \]
for $j \in {\mathbb Z}/(2 \nu(\varphi) {\mathbb Z})$. A priori it seems
that this does not make any sense since the domains of definition of
$\psi_{j}^{\varphi}$ and $\psi_{j+1}^{\varphi}$ are disjoint. Nevertheless
we can extend those domains, the function $\xi_{\varphi,K_{X}^{\mu}}^{j}$
will be defined in a strip.

We denote $D(\varphi)= {\mathbb Z}/(2 \nu(\varphi) {\mathbb Z})$. We define
\[ D_{1}(\varphi)= \{ j \in {\mathbb Z}/(2 \nu(\varphi) {\mathbb Z}) :
Re(X) \ {\rm points \ at} \
T_{\mu X}^{\epsilon, j}(0) \ {\rm towards} \ H(j) \} . \]
The condition $j \in D_{1}(\varphi)$ is equivalent to
$Re(- \mu X)$ pointing towards $|y|<\epsilon$ at
$(\partial H_{\infty}(j) \cap [|y|=\epsilon]) \setminus
T_{\mu X}^{\epsilon}$. We denote
$D_{-1}(\varphi)=D(\varphi) \setminus D_{1}(\varphi)$.

Suppose without lack of generality that $j \in D_{-1}(\varphi)$.
There exists a constant $W \in {\mathbb R}^{+}$ such that
$|Re(\psi_{j}^{X}(B) - \psi_{j}^{X}(A))|< W$ for all $A,B \in H_{\infty}(j)(x)$
and all $x \in [0,\delta_{0}) K_{X}^{\mu}$.
Denote $Im(x)=Img(\psi_{j}^{X}(T_{X}^{\epsilon, j+1}(x)))$.
We obtain that every
$Q \in H_{\infty}(j) \cap [Img \psi_{j}^{X} > Im]$
fulfills $[-W,W] \subset It(X,Q,|y|<\epsilon)$; we obtain
\[ {\rm exp}((0,W)X)(Q) \cap H(j+1) \neq \emptyset  \ {\rm and} \
{\rm exp}((-W,0)X)(Q) \cap H(j) \neq \emptyset . \]
Denote
$\Gamma_{l}(x)=\Gamma(\mu X, T_{\mu X}^{\epsilon, l}(x), |y| \leq \epsilon)$.
We define the strip $St_{j}(x)$ enclosed by
$\Gamma_{j}$ and $\varphi^{\circ (-1)}(St_{j}(x))$ whereas
$St_{j+1}(x)$ is the strip enclosed by $\Gamma_{j+1}(x)$ and
$\varphi(\Gamma_{j+1}(x))$ for all $x \in [0,\delta_{0}) K_{X}^{\mu}$.

The functions $\psi_{l}^{\varphi}-\psi_{l}^{X}$ are bounded in
$H(l)_{S}$ and continuous at the curve  $\omega^{\mu X}(H(j))$
for $l \in \{j,j+1\}$ (prop. \ref{pro:bddcon}).
Suppose that $\sup_{B(0,\delta) \times B(0,\epsilon)} |\Delta| < 1/2$.
It is easy to see that $\psi_{j}^{\varphi}$ can be defined by iteration
in the set $E_{j}$ given by
\[ E_{j}(x) =
([St_{j+1}(x) \cup \overline{H_{\infty}(j)(x))}] \setminus Fix \varphi)
\cap [Img(\psi_{j}^{X}) > Im(x) + 1 + W] \]
for $x \in [0,\delta_{0}) K_{X}^{\mu}$. The function $\psi_{j}^{\varphi}(x,.)$
is injective in the simply connected set
$H(j)_{S}(x) \cup E_{j}(x)$. Moreover
since we only need a finite number of iterations the function
$\psi_{j}^{\varphi}-\psi_{j}^{X}$ is still bounded in $E_{j}$ and
continuous at the curve $\omega^{\mu X}(H(j))$.
There exists $I \in {\mathbb R}^{+}$ such that
$\xi_{\varphi,K_{X}^{\mu}}^{j}$ is defined in
\[ [\cup_{x \in [0,\delta_{0})K_{X}^{\mu}}
\{ x \} \times \psi_{j}^{\varphi}(St_{j+1}(x))] \cap [Img(z)>I] . \]
Since we have
$\xi_{\varphi,K_{X}^{\mu}}^{j}(x,z+1) =\xi_{\varphi,K_{X}^{\mu}}^{j}(x,z) + 1$
then $\xi_{\varphi,K_{X}^{\mu}}^{j}$ is defined in $Img z>I$.
The value of
\[ \psi_{j+1}^{\varphi}-\psi_{j}^{\varphi} =
(\psi_{j+1}^{\varphi}-\psi_{j+1}^{X}) - (\psi_{j}^{\varphi}-\psi_{j}^{X})
+ (\psi_{j+1}^{X} - \psi_{j}^{X}) \]
at the curve $\gamma_{k(j)}=\omega^{\mu X}(H(j))$ is
$c_{k(j)} - c_{k(j)} + \zeta_{\varphi} \equiv \zeta_{\varphi}$, thus
$\xi_{\varphi,K_{X}^{\mu}}^{j}$ admits a expression of the type
$\xi_{\varphi, K_{X}^{\mu}}^{j}(x,z) =
z + \zeta_{\varphi}(x) +
\sum_{l=1}^{\infty} a_{j, l, K_{X}^{\mu}}^{\varphi}(x) {e}^{2 \pi i l z}$.
In particular the function $a_{j,l,K_{X}^{\mu}}^{\varphi}$ is continuous in
$[0,\delta_{0})K_{X}^{\mu}$ and holomorphic in
$(0,\delta_{0}) \dot{K}_{X}^{\mu}$ for all $l \in {\mathbb N}$.
The case $j \in D_{1}(\varphi)$ is analogous. The previous discussion implies:
\begin{pro}
\label{pro:ecvoinv}
Let $\varphi \in \diff{tp1}{2}$ with fixed convergent normal form
${\rm exp}(X)$. Fix $\mu \in e^{i(0,\pi)}$ and a compact connected set
$K_{X}^{\mu} \subset {\mathbb S}^{1} \setminus B_{X}^{\mu}$.
Then there exists $I \in {\mathbb R}^{+}$ such that
for all $s \in \{-1,1\}$ and $j \in D_{s}(\varphi)$ we have
\begin{itemize}
\item $\xi_{\varphi,K_{X}^{\mu}}^{j} \circ (x,z+1) \equiv (z+1) \circ
\xi_{\varphi,K_{X}^{\mu}}^{j}$.
\item $\xi_{\varphi,K_{X}^{\mu}}^{j} \in
C^{0}([0,\delta_{0}) K_{X}^{\mu} \times [s  Imgz < -I]) \cap
\vartheta((0,\delta_{0}) \dot{K}_{X}^{\mu} \times [s  Imgz < -I])$.
\item $\lim_{|Img(z)| \to \infty} \xi_{\varphi,K_{X}^{\mu}}^{j}(x,z) -
(z + \zeta_{\varphi}(x)) = 0$.
\item $\xi_{\varphi,K_{X}^{\mu}}^{j}$ is of the form
$z + \zeta_{\varphi}(x) + \sum_{l=1}^{\infty} a_{j,l,K_{X}^{\mu}}^{\varphi}(x)
{e}^{-2 \pi i s l z}$.
\end{itemize}
\end{pro}
Let $orb_{H,j}(\varphi)$ be the space of orbits of
$\varphi_{|H(j)_{S}}$ for $H(j) \in Reg(\epsilon, \mu X, K_{X}^{\mu})$.
The mapping $\Theta_{j}:orb_{H,j}(\varphi) \to [0,\delta_{0}) \times \pn{1}$
given by $\Theta_{j} \equiv (x,e^{2 \pi i z} \circ \psi_{j}^{\varphi})$
is continuous everywhere and holomorphic outside of $x=0$.
We define the $\mu$-{\it space of orbits} of $\varphi$ at $K_{X}^{\mu}$
as the variety obtained by taking an atlas composed of
charts $W_{j} \sim [0,\delta_{0}) \times \pn{1}$
for $j \in {\mathbb Z}/(2 \nu(\varphi) {\mathbb Z})$
and the changes of charts
$\Theta_{j+1} \circ \Theta_{j}^{\circ (-1)}$ identifying
subsets of $orb_{H,j}(\varphi)$ and $orb_{H,j+1}(\varphi)$
for all $j \in {\mathbb Z}/(2 \nu(\varphi) {\mathbb Z})$.

Let $j \in D_{s}(\varphi)$. The trajectory $t \to
{\rm exp}(s t X)(T_{\mu X}^{\epsilon, j}(0))$ (for $t \in {\mathbb
R}^{+}$) adheres to a direction $\Lambda(\varphi, j)
\in D_{s}(\varphi_{|x=0})$ when $t \to \infty$. The mapping
$\Lambda(\varphi)$
is a bijection from ${\mathbb Z}/ (2 \nu (\varphi) {\mathbb Z})$ to
$D(\varphi_{|x=0})$. The restriction
of the changes of charts to $x=0$ provide the Ecalle-Voronin
invariants of $\varphi_{|x=0}$.
\begin{cor}
\label{cor:EV0}
Let $\varphi \in \diff{tp1}{2}$ with fixed convergent normal form
${\rm exp}(X)$. Fix $\mu \in e^{i(0,\pi)}$ and a compact connected set
$K_{X}^{\mu} \subset {\mathbb S}^{1} \setminus B_{X}^{\mu}$. Then
the functions
$\xi_{\varphi, K_{X}^{\mu}}^{j}(0,z)$
($j \in {\mathbb Z}/ (2 \nu (\varphi) {\mathbb Z})$)
are the changes of charts of $\varphi_{|x=0}$. Indeed we have
$\xi_{\varphi, K_{X}^{\mu}}^{j}(0,z) \equiv
\xi_{\varphi_{|x=0}}^{\Lambda(\varphi, j)}(z)$
for all $j \in {\mathbb Z}/ (2 \nu (\varphi) {\mathbb Z})$.
\end{cor}
We have extended the Ecalle-Voronin invariants to all the lines $x=cte$
in a neighborhood of $x=0$
even if in general they do not support elements of $\diff{1}{}$.
\subsection{Nature of the invariants}
Let $X \in \Xt$.
In our sectors $[0,\delta_{0}) K_{X}^{\mu}$ the direction
$\mu \in {\mathbb S}^{1}$
providing the real flow $Re(\mu X)$ is fixed. The analogue
in \cite{MRR} is allowed to vary continuously. Such
a thing is also possible with our approach.

More precisely we want to find connected sets $E \subset {\mathbb S}^{1}$
and a continuous function $\mu: E \to e^{i(0,\pi)}$ such that
$\mu(\lambda) \not \in B_{X, \lambda}$ (see subsection \ref{subsubsec:diruns})
for all $\lambda \in E$. A maximal set with respect to the previous property
will be called a {\it maximal sector}. The idea is that for every compact
connected set $K$ contained in a maximal sector there exists
$\delta_{0}(K)>0$ such that
$Re(\mu(\lambda)X)_{
(r,\lambda,y) \in [0,\delta_{0}(K)) \times K \times B(0,\epsilon)}$
has a simple stable behavior. Thus the maximal sectors provide
sectorial domains of stability.

Let $\varphi \in \diff{tp1}{2}$ with convergent normal form ${\rm exp}(X)$.
Consider $x \in \lambda {\mathbb R}^{+}$ and $\mu_{0},\mu_{1}$ in the
same connected component $J$ of $e^{i(0,\pi)} \setminus B_{X, \lambda}$.
We claim that there exists a compact connected neighborhood
$K=K_{X}^{\mu_{0}} = K_{X}^{\mu_{1}}$
of $\lambda$ in ${\mathbb S}^{1}$ such that
$\xi_{\varphi,K_{X}^{\mu_{0}}}^{j} \equiv \xi_{\varphi,K_{X}^{\mu_{1}}}^{j}$
for all $j \in {\mathbb Z}$.
Then we can define changes of charts $\xi_{\varphi, K}^{j}(x,z)$
which are continuous in
$[x \in [0,\delta_{0}(K)) K] \cap [s Img z < - I]$ and
holomorphic in its interior for $j \in D_{s}(\varphi)$ .
Roughly speaking
this is a consequence of the continuous dependance of the regions
on $\mu \in J$.
The choices of $\mu$-spaces of orbits of $\varphi$ at
$x \in \lambda {\mathbb R}^{+}$
are at most the number of connected components of
$e^{i(0,\pi)} \setminus B_{X, \lambda}$.

It is remarkable that the dependance of $B_{X, \lambda}$ with
respect to $\lambda$ is not product-like. For instance
$B_{\beta, \lambda {e}^{i \theta}}(X)=e^{-i m_{\beta} \theta}
B_{\beta, \lambda}(x)$ for a magnifying glass $M_{\beta}$ associated
to $X$. Hence the points of $B_{X, \lambda}$ turn at different speeds.

There are much simpler cases. The diffeomorphisms considered in
\cite{MRR} are of the form
$\varphi(x,y) = (x,y - x + c_{1}(x) y^{2} + O(y^{3}))$
where $c_{1}(0) \neq 0$. They consider also its ramified
version $\tilde{\varphi}=(w^{1/2},y) \circ \varphi \circ (w^{2},y)$.
Given a convergent normal form ${\rm exp}(X)$ of $\tilde{\varphi}$ it is
easy to see that
$\sharp \tilde{B}_{X, \lambda} =1$ for all $\lambda \in {\mathbb S}^{1}$.
Then there are two choices of $\mu$-space of orbits in general.
There are two maximal sectors, they  are of the form
$\lambda_{0} e^{i(0,2 \pi)}$ for
$\lambda_{0} \in B_{X}^{1}$.  In the $x$ coordinate
only one sector is required to cover
${\mathbb S}^{1}$, it describes an angle as close to $4 \pi$ as desired.
We obtain the same division in the parameter space than in \cite{MRR};
nevertheless our techniques can be applied to every unfolding of tangent
to the identity germs and not only to the generic ones.
\subsection{Embedding in a flow}
Let $\varphi \in \diff{tp1}{2}$ with fixed convergent normal form
${\rm exp}(X)$. We say that a sequence
$K_{X}^{\mu_{1}}$, $\hdots$, $K_{X}^{\mu_{l}}$ of compact connected subsets of
${\mathbb S}^{1}$ is a EV-covering if
\begin{itemize}
\item $\mu_{j} \in e^{i(0,\pi)}$ and
$K_{X}^{\mu_{j}} \subset {\mathbb S}^{1} \setminus B_{X}^{\mu_{j}}$ for all
$j \in \{1, \hdots, l\}$.
\item $\cup_{j=1}^{l} \dot{K}_{X}^{\mu_{j}} = {\mathbb S}^{1}$.
\end{itemize}
Such a covering exists. We have $B_{X}^{i} \cap B_{X}^{\kappa}= \emptyset$
for $\kappa \in {\mathbb S}^{1}$ in the neighborhood of $i$. Fix such
$\kappa$, then we can choose a  EV-covering such that
$\{ \mu_{1}, \hdots, \mu_{l} \} \subset \{i, \kappa \}$.
In the trivial type case we can choose $K_{X}^{i} = {\mathbb S}^{1}$
as the only element of the EV-covering.
\begin{rem}
\label{rem:evndepnf}
The definition of EV-covering does not depend on the choice of the convergent
normal form but on
$Fix \varphi$ and $Res(\varphi)$ (remark \ref{rem:spndepnf}).
\end{rem}
\begin{pro}
\label{pro:nembflo}
Let $\varphi \in \diff{tp1}{2}$ with fixed convergent normal form
${\rm exp}(X)$. Fix $\mu \in e^{i(0,\pi)}$ and a compact connected set
$K_{X}^{\mu} \subset {\mathbb S}^{1} \setminus B_{X}^{\mu}$. Then
$\varphi$ is analytically trivial if and only if
$\xi_{\varphi, K_{X}^{\mu}}^{j} \equiv z + \zeta_{\varphi}$
for all $j \in {\mathbb Z}/(2 \nu(\varphi) {\mathbb Z})$.
\end{pro}
\begin{proof}
The sufficient condition is obvious.
The functions $\psi_{H}^{\varphi} - \psi_{H}^{X}$ paste together
for $H \in Reg(\epsilon, \mu X, K_{X}^{\mu})$ in a function
$J$ defined in
$[0,\delta_{0})K_{X}^{\mu} \times B(0,\epsilon')$
for some $0 < \epsilon' < \epsilon$. Moreover $J$ is
continuous in $[0,\delta_{0})K_{X}^{\mu} \times B(0,\epsilon')$ and analytic
in its interior (prop. \ref{pro:bddcon}) and
satisfies $J - J \circ \varphi = \Delta_{\varphi}$.
By Cauchy's integral formula we obtain
$|\partial J /\partial y| \leq M$ in $|y| < \epsilon'/2$ for some $M>0$.
We define the vector field
\[ X(K_{X}^{\mu}) = \frac{X(y)}{1+ X(y) \partial J/\partial y}
\frac{\partial}{\partial y} =
\frac{X(y)}{1+ X(J)} \frac{\partial}{\partial y},  \]
it coincides with $X_{H}^{\varphi}$ (see cor. \ref{cor:Lavvf}) for
$H \in Reg(\epsilon, \mu X, K_{X}^{\mu})$.
Since $X(K_{X}^{\mu})(\psi^{\varphi})=1$ then
$\varphi={\rm exp}(X(K_{X}^{\mu}))$
in $[0,\delta_{0})K_{X}^{\mu} \times B(0,\epsilon_{1})$ for some
$0 < \epsilon_{1} < \epsilon'/2$.

Consider a minimal EV-covering
$K_{1}=K_{X}^{\mu}$, $K_{2}=K_{X}^{\mu_{2}}$, $\hdots$,
$K_{l}=K_{X}^{\mu_{l}}$.
Consider $K_{b}$ such that $\dot{K}_{1} \cap \dot{K}_{b} \neq \emptyset$.
We define $\psi_{H,L}^{\varphi}=\psi_{H,L}^{X} + J$ and
$\psi_{H,R}^{\varphi}=\psi_{H,R}^{X} + J$
in $[0,\delta_{0})(K_{1} \cap K_{b}) \times B(0,\epsilon_{1})$ for all
$H \in Reg(\epsilon, \mu_{b} X, K_{b})$. Since
$J - J \circ \varphi = \Delta_{\varphi}$ then
$\psi_{H,L}^{\varphi}$ and $\psi_{H,R}^{\varphi}$ are
Fatou coordinates of $\varphi$ for
$H \in Reg(\epsilon, \mu_{b} X, K_{b})$. We obtain
$\xi_{\varphi,K_{b}}^{j} \equiv z + \zeta_{\varphi}$
for $j \in {\mathbb Z}/(2 \nu(\varphi) {\mathbb Z})$ in
$x \in \dot{K}_{1} \cap \dot{K}_{b}$ and then in $x \in K_{b}$ by
analytic continuation. Analogously to $X(K_{1})$ we can construct a vector
field $X(K_{b})$
such that $\varphi={\rm exp}(X(K_{b}))$ in
$[0,\delta_{0}) K_{b} \times B(0,\epsilon_{b})$
for some $\epsilon_{b}>0$. Moreover the construction implies that
$X(K_{1}) \equiv X(K_{b})$ in
$[0,\delta_{0})(\dot{K}_{1} \cap \dot{K}_{b}) \times
B(0,\min(\epsilon_{1},\epsilon_{b}))$. Finally we obtain
$Y \in {\mathcal X} \cn{2}$ of the form
$Y = X(y) (1 + X(y) A) \partial / \partial y$ for some
$A \in {\mathbb C} \{x,y\}$
such that $\varphi={\rm exp}(Y)$.
Since $Y$ is nilpotent then $\log \varphi=Y$.
\end{proof}
\section{Complete system of analytic invariants. Trivial type case}
\label{sec:trityp}
We suppose  throughout this section
that $I(Fix \varphi)=(y)$ for all $\varphi \in \diff{p1}{2}$
with $Fix \varphi$ of trivial type. This is possible up to change
of coordinates $(x,y+h(x))$.
\begin{lem}
\label{lem:tflog}
Let $\varphi \in \diff{p1}{2}$ such that $Fix \varphi$ is of trivial type.
Then $(\log \varphi)(y)$ belongs to $\vartheta(B(0,\delta))[[y]]$ for some
$\delta \in {\mathbb R}^{+}$.
\end{lem}
\begin{proof}
Suppose that $\varphi$ is defined in $B(0,\delta) \times B(0,\epsilon)$.
Denote by $\Theta$ the operator $\varphi - Id$. We have
$(\log \varphi)(y) - \sum_{j=1}^{l} (-1)^{j+1}
\Theta^{\circ (j)}(y)/j \in ({y}^{\nu(\varphi)+l+1})$
by the proof of proposition \ref {pro:excofor}.
We are done since $\Theta^{\circ (j)}(y)$ is holomorphic in the neighborhood
of $B(0,\delta) \times \{0\}$ for all $j \in {\mathbb N}$.
\end{proof}
Let $\varphi_{1},\varphi_{2} \in \diff{p1}{2}$ with common convergent
normal form such that
$Fix \varphi_{1}$ is of trivial type. We define
$\hat{\sigma}(\varphi_{1},\varphi_{2})=
\hat{\sigma}(\varphi_{1},\varphi_{2},Fix \varphi_{1})$.
We say that $\hat{\sigma}(\varphi_{1},\varphi_{2})$ is the
{\it privileged formal conjugation} between $\varphi_{1}$ and
$\varphi_{2}$. By construction we obtain that
$y \circ \hat{\sigma}(\varphi_{1},\varphi_{2}) - y \in
({y}^{\nu(\varphi_{1})+2})$.
\begin{lem}
\label{lem:tfpri}
Let $\varphi_{1},\varphi_{2} \in \diff{p1}{2}$ with common convergent
normal form and $Fix \varphi_{1}$ of trivial type. Then
$y \circ \hat{\sigma}(\varphi_{1},\varphi_{2}) \in \vartheta(B(0,\delta))[[y]]$
for some $\delta \in {\mathbb R}^{+}$.
\end{lem}
\begin{proof}
We have $(\log \varphi_{1})(y), (\log \varphi_{2})(y) \in
\vartheta(B(0,\delta))[[y]]$ for some
$\delta \in {\mathbb R}^{+}$ by lemma \ref{lem:tflog}. Consider
$\hat{\beta} \in {\mathbb C}[[x,y]]$ such that
$\partial \hat{\beta}/\partial y = 1/(\log \varphi_{1})(y) -
1/(\log \varphi_{2})(y)$ and $\hat{\beta}(x,0) \equiv 0$.
We deduce that $\hat{\beta}$ and then
$y \circ \hat{\sigma}(\varphi_{1},\varphi_{2})$
belong to $\vartheta(B(0,\delta))[[y]]$ by proposition \ref{pro:despcon}.
\end{proof}
\begin{lem}
\label{lem:tfsim}
Let $\varphi \in \diff{p1}{2}$ such that $Fix \varphi$ is of trivial type.
Then $y \circ \hat{\tau}_{0}(\varphi)$ belongs to
$\vartheta(B(0,\delta))[[y]]$ for some $\delta \in {\mathbb R}^{+}$.
\end{lem}
\begin{proof}
Denote $\nu=\nu(\varphi)$ and
$X=y^{\nu+1}/(1+Res(\varphi,(x,0)) y^{\nu}) \partial / \partial y$.
Consider a convergent normal form ${\rm exp}(X_{1})$ of $\varphi$.
We obtain that ${\rm exp}(X)$ and ${\rm exp}(X_{1})$ are conjugated
by some
$\sigma_{0} \in \diff{p}{2}$ (prop. \ref{pro:anconfl}). We have that
\[ \hat{\tau}_{0}(\varphi) = \hat{\sigma}({\rm exp}(X_{1}),\varphi)
\circ \sigma_{0} \circ
(x,e^{2 \pi i / \nu(\varphi)} y) \circ \sigma_{0}^{\circ (-1)} \circ
\hat{\sigma}({\rm exp}(X_{1}),\varphi)^{\circ (-1)} .\]
Now $y \circ  \hat{\sigma}({\rm exp}(X_{1}),\varphi)$ belongs to
$\vartheta(B(0,\delta))[[y]]$ for some $\delta \in {\mathbb R}^{+}$
by the previous lemma. Therefore $y \circ \hat{\tau}_{0}(\varphi)$
belongs to $\vartheta(B(0,\delta))[[y]]$.
\end{proof}

Let $\varphi \in  \diff{p1}{2}$. Suppose that $Fix \varphi$ is of trivial
type. We define $\varphi_{w}$ as the germ of $\varphi_{|x=w}$ in the
neighborhood of $y=0$.

Fix a convergent normal form ${\rm exp}(X)$ of $\varphi$.
We choose an EV-covering with a unique element $K_{X}^{i}={\mathbb S}^{i}$.
We denote $\xi_{\varphi, K_{X}^{i}}^{j}$ by either $\xi_{\varphi}^{j}$ or
$\xi_{\varphi}^{\Lambda(\varphi,j)}$.
We have that $\xi_{\varphi}^{\lambda}$ is holomorphic in
$B(0,\delta_{0}) \times [sImgz < -I]$ for some $I \in {\mathbb R}^{+}$
and all $\lambda \in D_{s}(\varphi_{0})$ by proposition
\ref{pro:ecvoinv}.
We obtain that $\xi_{\varphi}^{j}$ is of the form
\[ \xi_{\varphi}^{j}(x,z) =
z - \pi i Res(\varphi,(x,0))/\nu(\varphi) +
\sum_{k=1}^{\infty} a_{j, k}^{\varphi}(x) {e}^{-2 \pi i s k z} \]
where $\sum_{k=1}^{\infty} a_{j, k}^{\varphi}(x) w^{k}$ is an
analytic function in a neighborhood
of $(x,w)=(0,0)$. A different choice of convergent normal form or
homogeneous coordinates provides new Fatou coordinates
$\psi_{\varphi}^{j}(x,z)+t(x)$ for some $t \in {\mathbb C}\{x\}$
independent of $j \in D(\varphi)$. The changes of charts
are unique up to conjugation with $z+t(x)$ for some $t \in {\mathbb C}\{x\}$.

Let $\varphi_{1}, \varphi_{2} \in \diff{p1}{2}$ with common
convergent normal form ${\rm exp}(X)$. We always suppose that their Fatou
coordinates are calculated with respect to a common system of
homogeneous coordinates. Since $\hat{\tau}_{0}(\varphi_{2})$ and
$\hat{\sigma}(\varphi_{1},\varphi_{2})$ depend analytically on $x$
by lemmas \ref{lem:tfsim} and \ref{lem:tfpri} then there are
parameterized versions of the results in
subsections \ref{subsec:topbeh}, \ref{subsec:anapro} and
\ref{subsec:anacla}. We obtain:
\begin{pro}
\label{pro:ancottc}
Let $\varphi_{1},\varphi_{2} \in  \diff{p1}{2}$ with common convergent
normal form ${\rm exp}(X)$.
Suppose that $Fix \varphi_{1}$ is of trivial type.
Then $\varphi_{1} \sim \varphi_{2}$ if and only if there exists
$(k,t) \in {\mathbb Z}/(\nu(X) {\mathbb Z}) \times {\mathbb C}\{x\}$
such that
\begin{equation}
\label{equ:eqca2}
 \xi_{\varphi_{2}}^{j + 2k} (x,z+t(x)) = (z+t(x)) \circ
\xi_{\varphi_{1}}^{j}(x,z) \ \ \forall j \in D(\varphi_{1}) .
\end{equation}
The equation \ref{equ:eqca2} is equivalent to
$Z_{\varphi_{2}}^{\kappa,t} \circ
\hat{\sigma}(\varphi_{1},\varphi_{2}) \in \diff{}{2}$
where $\kappa = e^{2 \pi i k/ \nu(X)}$.
\end{pro}
Let $\varphi_{1}, \varphi_{2} \in \diff{p1}{2}$ with
$Fix \varphi_{1}=Fix \varphi_{2}$ of trivial type.
We say that $m_{\varphi_{1}}(w) = m_{\varphi_{2}}(w)$ if
$(\varphi_{1})_{w} \sim (\varphi_{2})_{w}$.
We denote $Inv(\varphi_{1}) \sim Inv(\varphi_{2})$ if there exists
$(k(x),d(x)) \in {\mathbb Z}/(\nu(X) {\mathbb Z}) \times [|Img(z)| < I]$
such that
\[ \xi_{\varphi_{2}}^{j + 2k(x)} (x,z+d(x)) = (z+d(x)) \circ
\xi_{\varphi_{1}}^{j}(x,z) \ \ \forall j \in D(\varphi_{1}) . \]
for all $x \neq 0$ in a neighborhood of $0$ and some $I \in {\mathbb R}^{+}$.
Consider the set
\[ E_{s}(\varphi)= \{ (j,k) \in D_{s}(\varphi) \times {\mathbb N}
\ s. t. \ a_{j,k}^{\varphi} \not \equiv  0\} . \]
We define $E(\varphi_{1})=E_{-1}(\varphi) \cup E_{1}(\varphi_{1})$.
The definitions  of $E_{-1}(\varphi_{1})$, $E_{1}(\varphi_{1})$ and
$E(\varphi)$ do not depend on the choice of homogeneous coordinates.
\begin{pro}
\label{pro:anainv}
Let $\varphi_{1}, \varphi_{2} \in \diff{p1}{2}$ with
$Fix \varphi_{1} = Fix \varphi_{2}$ of trivial type. Then we have
$\varphi_{1} \sim \varphi_{2}$ if and only if
$Inv(\varphi_{1}) \sim Inv(\varphi_{2})$.
\end{pro}
This result provides a complete system of analytic invariants
in the trivial type case; it is composed by the changes of charts modulo
uniform changes of coordinates.
\begin{proof}
The condition $Inv(\varphi_{1}) \sim Inv(\varphi_{2})$ implies in
particular $Res(\varphi_{1}) \equiv Res(\varphi_{2})$. Let $\alpha_{j}$ be
a convergent normal form of $\varphi_{j}$ for $j \in \{1,2\}$.
Thus $\alpha_{1}$ and $\alpha_{2}$ are conjugated by
$\sigma \in \diff{p}{2}$ (prop. \ref{pro:anconfl}).
By replacing $\varphi_{2}$ with
$\sigma^{\circ (-1)} \circ \varphi_{2} \circ \sigma$ and
$\xi_{\varphi_{2}}^{j}(x,z)$ with
$(z+t(x)) \circ \xi_{\varphi_{2}}^{j + 2k_{0}} \circ (x,z-t(x))$
for all $j \in D(\varphi_{2})$ and some
$(k_{0},t) \in {\mathbb Z}/(\nu(X) {\mathbb Z}) \times {\mathbb C}\{x\}$
we can suppose that $\varphi_{1}$ and $\varphi_{2}$ have common
convergent normal form ${\rm exp}(X)$. We can also suppose that
$\log \varphi_{1} \not \in {\mathcal X} \cn{2}$, otherwise we get
$\varphi_{1} \sim \varphi_{2}$ (prop. \ref{pro:anconfl}).
Suppose there exists
$I \in {\mathbb R}^{+}$ such that
\[ \xi_{\varphi_{2}}^{j + 2 k(x)} (x,z+d(x))=(z+d(x)) \circ
\xi_{\varphi_{1}}^{j}(x,z) \ \ \forall j \in D(\varphi_{1}) \]
for some $(k(x), d(x)) \in {\mathbb Z}/(\nu(X) {\mathbb Z})
\times [|Img z| < I]$ and all $x \neq 0$.
We choose $k \in {\mathbb Z}/(\nu(X) {\mathbb Z})$ such that
$[k(x)=k]$ is uncountable in every neighborhood of $0$ and
$x_{0} \in [k(x)=k] \setminus \{0\}$ such that
$a_{j,l}^{\varphi_{1}}(x_{0}) \neq 0$ for all
$(j,l) \in E(\varphi_{1})$.
Fix $(j_{0},l_{0}) \in E_{s}(\varphi_{1})$; since
$e^{-2 \pi l_{0} I}  \leq
|a_{j_{0} + 2 k,l_{0}}^{\varphi_{2}}/
a_{j_{0},l_{0}}^{\varphi_{1}}|(0) \leq
e^{2 \pi l_{0} I}$ then
there exists a holomorphic function $m$ defined in an open set containing
$0$ and $x_{0}$ such that $m(x_{0})=d(x_{0})$ and
$a_{j_{0} + 2 k,l_{0}}^{\varphi_{2}}/a_{j_{0},l_{0}}^{\varphi_{1}}
=e^{2 \pi i s l_{0} m}$. Since for all
$(j,l) \in E_{s'}(\varphi_{1})$ and $s' \in \{-1,1\}$ we have
$(a_{j_{0} + 2 k,l_{0}}^{\varphi_{2}}/
a_{j_{0},l_{0}}^{\varphi_{1}})^{s l} =
(a_{j + 2 k,l}^{\varphi_{2}}/
a_{j,l}^{\varphi_{1}})^{s' l_{0}}$ then
$m$ does not depend on the choice of $(j_{0},l_{0})$.
We obtain that
$Z_{\varphi_{2}}^{\kappa,m} \circ \hat{\sigma}(\varphi_{1},\varphi_{2})$
is an analytic mapping conjugating
$\varphi_{1}$ and $\varphi_{2}$ by proposition \ref{pro:ancottc}
where $\kappa = e^{2 \pi i k/ \nu(X)}$.
\end{proof}
\begin{cor}
\label{cor:rig}
Let $\varphi_{1},\varphi_{2} \in  \diff{p1}{2}$ such that
$Fix \varphi_{1} = Fix \varphi_{2}$.
Suppose that $Fix \varphi_{1}$ is of trivial type and that
$\log (\varphi_{1})_{0} \not \in {\mathcal X} \cn{}$.
Then $\varphi_{1} \sim \varphi_{2}$
$\Leftrightarrow$
$m_{\varphi_{1}} \equiv m_{\varphi_{2}}$.
\end{cor}
\begin{proof}
There exists $(j_{0},l_{0}) \in E_{s}(\varphi_{1})$ such that
$a_{j_{0},l_{0}}^{\varphi_{1}}(0) \neq 0$. There also exists
$(j_{1},l_{0}) \in E_{s}(\varphi_{2})$ such that
$a_{j_{1},l_{0}}^{\varphi_{2}}(0) \neq 0$
since $m_{\varphi_{1}}(0) = m_{\varphi_{2}}(0)$. We have
\[ \left\{
{\begin{array}{c}
a_{j_{0} + 2 k(x) ,l_{0}}^{\varphi_{2}}(x) =
a_{j_{0},l_{0}}^{\varphi_{1}}(x) e^{2 \pi i s l_{0} d(x)} \\
a_{j_{1}, l_{0}}^{\varphi_{2}}(x) =
a_{j_{1} - 2 k(x), l_{0}}^{\varphi_{1}}(x)
e^{2 \pi i s l_{0} d(x)} .
\end{array} }\right. \]
for some $(k(x),d(x)) \in {\mathbb Z}/(\nu(X) {\mathbb Z}) \times
{\mathbb C}$ and all $x$ in a neighborhood of $0$ by hypothesis.
We deduce that $Img d$ is bounded. Thus we obtain
$\varphi_{1} \sim \varphi_{2}$ (prop. \ref{pro:anainv}).
\end{proof}
We say that $\eta$ is a {\it r-mapping} if
$\eta$ is a biholomorphism from $B(0,r)$ onto
$\eta(B(0,r))$. If $\eta(B(0,r))$ is contained in $B(0,R)$ then
we say that $\eta$ is a rR-mapping.

Next we provide a geometrical interpretation of the system of invariants.
\begin{pro}
\label{pro:mod}
Let $\varphi_{1},\varphi_{2} \in  \diff{p1}{2}$ such that
$Fix \varphi_{1} = Fix \varphi_{2}$.
Suppose that $Fix \varphi_{1}$ is of trivial type.
Then $\varphi_{1} \sim \varphi_{2}$
if there exist $r \in {\mathbb R}^{+}$ and a r-mapping
$\eta_{x}$ conjugating  $(\varphi_{1})_{x}$ and $(\varphi_{2})_{x}$
for all $x$ in a pointed neighborhood of $0$.
\end{pro}
We do not ask $\eta_{x}$ to have any kind of good dependance
with respect to $x$.
\begin{proof}
By proposition \ref{pro:cftg} we have that
$\nu(\varphi_{1})=\nu(\varphi_{2})$ and
$Res(\varphi_{1}) \equiv Res(\varphi_{2})$.
Let $\alpha_{j}$ be a convergent normal form of $\varphi_{j}$
for $j \in \{1,2\}$. Let $\zeta \in \diff{p}{2}$ be the mapping
conjugating $\alpha_{1}$ and $\alpha_{2}$ provided by
proposition \ref{pro:anconfl}.
The mapping $\zeta_{x}^{\circ (-1)} \circ \eta_{x}$ conjugates
diffeomorphisms with common convergent normal form $(\alpha_{1})_{x}$.
We obtain
$|(\partial (\zeta_{x}^{\circ (-1)} \circ \eta_{x})/\partial y)(0)|=1$
(prop. \ref{pro:comonevar}).
Denote $b(x)=(\partial \eta_{x}/\partial y)(0)$. We have that
$\eta_{x}(ry)/(r b(x))$ is a Schlicht  function
for all $x$ in a pointed neighborhood of $0$. By the Koebe's distortion
theorem (see \cite{Conway}, page 65) we get
\[ \sup_{y \in B(0,r_{1})} |\eta_{x}(y)| \leq
r |b(x)|
\sup_{y \in B(0,r_{1}/r)} \left|{ \frac{\eta_{x}(r y)}{r b(x)} }\right| \leq
r \left|{ \frac{\partial (y \circ \zeta)}{\partial y}(x,0) }\right|
\frac{r_{1}/r}{(1-r_{1}/r)^{2}} \]
for all $r_{1} < r$ and all $x$ in a pointed neighborhood of $0$.
We deduce that $\zeta_{x}^{\circ (-1)} \circ \eta_{x}$ is a rR-mapping
for some $R \in {\mathbb R}^{+}$ by considering a smaller $r>0$ if
necessary. By replacing $\varphi_{2}$ with
$\zeta^{\circ (-1)} \circ {\varphi}_{2} \circ \zeta$
and $\eta_{x}$ with  $\zeta_{x}^{\circ (-1)} \circ \eta_{x}$ we can
suppose that $\varphi_{1}$ and $\varphi_{2}$ have common normal form.

The mapping $\eta_{w}$ is of the form
$Z_{(\varphi_{2})_{w}}^{\kappa(w),d(w)} \circ
\hat{\sigma}((\varphi_{1})_{w}, (\varphi_{2})_{w})$ since
it conjugates $(\varphi_{1})_{w}$ and $(\varphi_{2})_{w}$
where $(\kappa(w),d(w))$ belongs to
$<e^{2 \pi i/\nu(\varphi_{1})}> \times {\mathbb C}$
for all $w$ in a pointed neighborhood of $0$.
We want to estimate $d(x)$. We have
\[ y \circ \eta_{w} - y \circ (Z_{\varphi_{2}}^{\kappa(w),0})_{|x=w}
- \kappa(w) d(w) \frac{(\log \varphi_{2})(y)}{y^{\nu(\varphi_{1})+1}}(w,0)
y^{\nu(\varphi_{1})+1} \in (y^{\nu(\varphi_{1})+2})  \]
for all $w \neq 0$. The series
$[(\log \varphi_{2})(y)/y^{\nu(\varphi_{1})+1}](x,0)$ is a unit
of ${\mathbb C}\{x\}$ by lemma \ref{lem:tflog}.
Moreover $y \circ Z_{\varphi_{2}}^{\kappa,0} \in
\vartheta(B(0,\delta))[[y]]$ for all
$\kappa \in <e^{2 \pi i/\nu(\varphi_{1})}>$ and some
$\delta \in {\mathbb R}^{+}$ by lemma \ref{lem:tfsim}. Since
\[ \left|{ \frac{1}{(\nu(\varphi_{1})+1)!}
\frac{\partial^{\nu(\varphi_{1})+1} \eta_{x}}
{\partial y^{\nu(\varphi_{1})+1}} }\right| (0) =
\left|{ \frac{1}{2 \pi i} \int_{|y|=r/2}
\frac{y \circ \eta_{x}(y)}{y^{\nu(\varphi_{1})+2}} }\right| \leq
\frac{2^{\nu(\varphi_{1})+1} R}{r^{\nu(\varphi_{1})+1}}  \]
then $d$ is bounded. We deduce
$\varphi_{1} \sim \varphi_{2}$ by proposition \ref{pro:anainv}.
\end{proof}
\section{Applications}
\label{sec:app}
In this section we complete the task of classifying analytically the
elements of $\diff{p1}{2}$. Moreover given
$\varphi_{1} \sim \varphi_{2}$ we provide the formal
power series developments of the conjugating diffeomorphisms.
We also relate the analytic class of $\varphi \in \diff{p1}{2}$ and
the analytic classes of the elements of
${\{ \varphi_{|x=x_{0}} \}}_{x_{0} \in B(0,\delta_{0})}$.
\subsection{Uniform conjugations}
We denote by $\diff{xp1}{2}$ the subset obtaining by removing from
$\diff{tp1}{2}$ the elements with fixed points set of trivial type.
We want to identify how an analytic conjugation between elements of
$\diff{xp1}{2}$ acts on
the changes of charts. We remind the reader that $N(X)$ is the number
of points in $(Sing X)(x_{0})$ for $x_{0}$ generic in a neighborhood
of $0$.
\begin{lem}
\label{lem:quitR}
Let $X \in {\mathcal X} \cn{2}$ with $N(X) \geq 2$.
Fix $r \geq 0$.
There exists a function
$R:(0,r) \to {\mathbb R}^{+}$ with $\lim_{b \to 0} R(b)=0$
such that all r-mapping
$\kappa$ holding $\kappa_{|(Sing X)(x_{0})} \equiv Id$
is a $r_{1}R(r_{1})$-mapping for $x_{0} \neq 0$ in a
neighborhood $V(r_{1})$ of $0$.
\end{lem}
\begin{proof}
Let $\gamma_{1}(x_{0})$ and $\gamma_{2}(x_{0})$ be two different points of
$(Sing X)(x_{0})$. We define
\[ \kappa_{1}(y) =
\frac{\kappa((r-|\gamma_{1}(x_{0})|)y + \gamma_{1}(x_{0})) -\gamma_{1}(x_{0})}
{(r-|\gamma_{1}(x_{0})|) (\partial \kappa/\partial y)(\gamma_{1}(x_{0}))} .\]
Then $\kappa_{1}$ is a Schlicht function.
Denote $\upsilon(x_{0}) =
(\gamma_{2}(x_{0}) - \gamma_{1}(x_{0}))/(r-|\gamma_{1}(x_{0})|)$. We have
$\kappa_{1} (\upsilon(x_{0})) = \upsilon(x_{0})/
(\partial \kappa/\partial y)(\gamma_{1}(x_{0}))$.
Koebe's distortion theorem (see \cite{Conway}, page 65) implies
$| (\partial \kappa / \partial y)(\gamma_{1}(x_{0})) | \leq
(1 + |\upsilon(x_{0})|)^{2}$. We have
\[ \sup_{y \in B(0,r_{1})} |\kappa(y)| \leq
(r - |\gamma_{1}(x_{0})|) (\partial \kappa/\partial y)(\gamma_{1}(x_{0}))
\sup_{y \in B(0,A(r_{1}))}
|\kappa_{1}(y)| + |\gamma_{1}(x_{0})| \]
where $A(r_{1})=(r_{1}+|\gamma_{1}(x_{0})|)/(r-|\gamma_{1}(x_{0})|)$.
Since again by Koebe's distortion theorem we have
$\sup_{y \in B(0,A(r_{1}))} |\kappa_{1}(y)| \leq
A(r_{1})/(1-A(r_{1}))^{2}$  then the
value $R(r_{1})$ can be chosen as close to $r_{1}/(1-r_{1}/r)^{2}$ as desired.
\end{proof}
The last lemma implies that in our context
the existence of r and rR conjugating mappings
are equivalent concepts.
\begin{lem}
\label{lem:modbdd}
Let $\varphi_{1}, \varphi_{2} \in \diff{xp1}{2}$ with common convergent
normal form ${\rm exp}(X)$.
There exist an open set $0 \in V \subset {\mathbb C}$ and
$D(r,R) \in {\mathbb R}^{+}$
such that a rR-mapping $\kappa$ conjugating
$(\varphi_{1})_{|x=x_{0}}$ and
$(\varphi_{2})_{|x=x_{0}}$ can be expressed in the form
$y + X(y)(x_{0},y) J_{\kappa}(y)$ where $\sup_{B(0,r)} |J_{\kappa}| < D(r,R)$
for all $x_{0} \in V \setminus \{ 0 \}$.
\end{lem}
\begin{proof}
Denote
$X(y)=u(x,y) (y-\gamma_{1}(x))^{n_{1}} \hdots (y - \gamma_{N}(x))^{n_{N}}$
where
$u \in {\mathbb C} \{x,y\}$ is a unit.
By hypothesis we have
$\kappa = y + (y-\gamma_{1}(x_{0})) \hdots (y - \gamma_{N}(x_{0})) A(y)$
for some $A \in \vartheta(B(0,r))$. By the modulus maximum principle we obtain
\[ \sup_{B(0,r)} |A| = \lim_{s \to r} \sup_{y \in B(0,s)}
\frac{|\kappa(y) -y|}{|(y-\gamma_{1}(x_{0})) \hdots (y - \gamma_{N}(x_{0}))|}
\leq \frac{r+R}{(r/2)^{N}}\]
for all $x_{0}$ in a pointed neighborhood of $0$. We have that
\[ \left|{ \frac{\partial \kappa}{\partial y }(\gamma_{j}(x_{0})) -1 }\right|
\leq \frac{2^{N} (r+R)}{r^{N}}
\prod_{k \in \{1,\hdots, N\} \setminus \{ j \}}
|\gamma_{j}(x_{0}) - \gamma_{k}(x_{0})| . \]
Fix $j \in \{1,\hdots,N\}$. We claim that
$(y-\gamma_{j}(x_{0}))^{n_{j}}$ divides $\kappa$.
We can suppose $n_{j}>1$. Denote by $\zeta_{1}$,
$\zeta_{2}$ and $\upsilon$ the germs of
diffeomorphism induced by $(\varphi_{1})_{|x=x_{0}}$,
$(\varphi_{2})_{|x=x_{0}}$  and $\kappa$ respectively
in the neighborhood of $x_{0}$. We have
$\upsilon = Z_{\zeta_{2}}^{\lambda,t} \circ \hat{\sigma}(\zeta_{1},\zeta_{2})$
for some $t \in {\mathbb C}$ and
$\lambda =(\partial \kappa/\partial y)(\gamma_{j}(x_{0}))
\in <e^{2 \pi i/(n_{j}-1)}>$
(prop. \ref{pro:comonevar}).
This implies $\lambda = 1$ for
$x_{0}$ in a neighborhood of $0$ since $N \geq 2$. Thus
$y \circ \kappa -y \in (y-\gamma_{j}(x_{0}))^{n_{j}}$. Denote
$J_{\kappa}=(\kappa-y)/X(y)_{|x=x_{0}}$, it belongs to $\vartheta(B(0,r))$.
Analogously than
for $A$ we obtain
$\sup_{B(0,r)} |J_{\kappa}| \leq D(r,R)$ for some
$D(r,R) \in {\mathbb R}^{+}$ and all $x_{0} \neq 0$.
\end{proof}
\begin{lem}
\label{lem:closide}
Let $\varphi_{1}, \varphi_{2} \in \diff{xp1}{2}$ with common convergent
normal form ${\rm exp}(X)$.
Fix $r,R$ in ${\mathbb R}^{+}$ and $0 < r_{1} <r$.
There exist $M(r,R,r_{1}) \in {\mathbb R}^{+}$ and a neighborhood
$V \subset {\mathbb C}$ of $0$
such that a rR-mapping $\kappa$ conjugating
$\varphi_{1}(x_{0},y)$ and
$\varphi_{2}(x_{0},y)$ satisfies
$\sup_{B(0,r_{1})} | \partial \kappa / \partial y - 1 | \leq M(r,R,r_{1})$
for all $x_{0} \in V \setminus \{ 0 \}$. Moreover we have $\lim_{r_{1} \to 0}
M(r,R,r_{1}) =0$.
\end{lem}
\begin{proof}
Denote $X(y)=u(x,y)
\prod_{j=1}^{N} (y-\gamma_{j}(x))^{n_{j}}$ where
$u \in {\mathbb C} \{x,y\}$ is a unit. By lemma \ref{lem:modbdd} we have that
$\kappa$ is of the form
$y + A(y) \prod_{j=1}^{N} (y-\gamma_{j}(x))^{n_{j}}$
for some $A \in \vartheta(B(0,r))$. We have
$\sup_{B(0,r)} |A| \leq H(r,R)$ for some $H(r,R) \in {\mathbb R}^{+}$ and all
$x_{0}$ in a pointed neighborhood of $0$. Fix $0 < r_{1} < r$.
Cauchy's integral formula implies
$\sup_{y \in B(0,r_{1})} |\partial A/ \partial y| \leq H(r,R)/(r-r_{1})$.
Thus we get
\[ \left|{ \frac{\partial \kappa}{\partial y }(y) -1 }\right| \leq
H(r,R)  (\nu(X)+1) {\left({ 2r_{1} }\right)}^{\nu(X)}  +
\frac{H(r,R)}{r-r_{1}} {\left({ 2r_{1} }\right)}^{\nu(X)+1}  \]
for $y \in B(0,r_{1})$.
We define $M(r,R,r_{1})$ as the right hand side of the previous formula.
Clearly we have $\lim_{r_{1} \to 0} M(r,R,r_{1})=0$.
\end{proof}
Given a rR-conjugation $\kappa$ we can suppose that
$\sup_{B(0,r)} |\partial \kappa / \partial y -1|$
is as small as desired just by considering a smaller $r>0$
independent on $x_{0}$.
%
We will make this kind of assumption
without stressing it every time. We define
$\kappa_{t}(y) = y + t (\kappa(y)-y)$
for $y \in B(0,r)$ and $t \in {\mathbb C}$.
\begin{lem}
Let $\varphi_{1}, \varphi_{2} \in \diff{xp1}{2}$ with common convergent
normal form ${\rm exp}(X)$.
Fix $r,R \in {\mathbb R}^{+}$. There exist $0<r_{1}<r$ and an open set
$0 \in V \subset {\mathbb C}$ such that for all
rR-mapping $\kappa$ conjugating $\varphi_{1}(x_{0},y)$ and
$\varphi_{2}(x_{0},y)$ and all $x_{0} \in V \setminus \{0\}$ we have that
$\kappa_{t}$ is a $r_{1}R$-mapping for all $t \in B(0,2)$.
\end{lem}
\begin{proof}
We can choose $0 < r_{1} < \min(r,R/7)$ such that
$\sup_{B(0,r_{1})} |\partial \kappa / \partial y -1| \leq 1/4$
by lemma \ref{lem:closide}. Therefore we obtain
$\sup_{B(0,r_{1})} |\kappa| \leq 2 r_{1}$ for all $x_{0}$ in a pointed
neighborhood $V(r_{1})$ of $0$. This implies
$\sup_{B(0,r_{1})} |\kappa_{t}| \leq 7 r_{1} < R$ for all $t \in B(0,2)$.
Moreover since
$\sup_{B(0,r_{1})} |\partial \kappa_{t} / \partial y -1| \leq 1/2$
then $\kappa_{t}$ is injective and hence a $r_{1}R$-mapping for all
$t \in B(0,2)$.
\end{proof}
Let $\psi^{X}$ be a holomorphic integral of the time form of $X$.
We can define the function $\psi^{X} \circ \kappa(x,y)  - \psi^{X}(x,y)$ in
an analogous way than $\Delta_{\varphi}$.
The continuous path that we use to extend $\psi^{X}$ is parameterized by
$t \to \kappa_{t}(x,y)$ for $t \in [0,1]$. The function
$\psi^{X} \circ \kappa - \psi^{X}$ is well-defined
and holomorphic in $B(0,r) \setminus Sing X$.
\begin{lem}
\label{lem:bddfldi}
Let $\varphi_{1}, \varphi_{2} \in \diff{xp1}{2}$ with common convergent
normal form
${\rm exp}(X)$. Fix $r,R$ in ${\mathbb R}^{+}$. Then there exist $0<r_{1}<r$ and
$C(r,R)>0$ such that we have
$\sup_{B(0,r_{1})} |\psi^{X} \circ \kappa - \psi^{X}| \leq C(r,R)$
for all rR-mapping $\kappa$ conjugating $\varphi_{1}(x_{0},y)$ and
$\varphi_{2}(x_{0},y)$ and all
$x_{0}$ in a pointed neighborhood of $0$. In particular we obtain that
$\psi^{X} \circ \kappa - \psi^{X}$ belongs to $\vartheta(B(0,r_{1}))$.
\end{lem}
\begin{proof}
Denote
$X(y)=u(x,y) (y-\gamma_{1}(x))^{n_{1}} \hdots (y - \gamma_{N}(x))^{n_{N}}$
where $u \in {\mathbb C} \{x,y\}$ is a unit. There exits a positive real
number $H(r,R)$ such that
\[ \left|{ \frac{\partial \kappa_{t}}{\partial t}(y) }\right|
\leq \frac{H(r,R)}{|u \circ \kappa_{t}(y)|} |X(y) \circ \kappa_{t}(y)|
\left|{ \frac{\prod_{j=1}^{N} (y-\gamma_{j}(x_{0}))^{n_{j}}}
{\prod_{j=1}^{N} (y-\gamma_{j}(x_{0}))^{n_{j}} \circ \kappa_{t}(y)} }\right| \]
for all $y \in B(0,r) \setminus (Sing X)(x_{0})$.
Denote $C(r,R)=2^{\nu(X)+1} H(r,R)/\inf_{B(0,R)} |u|$. Since $\nu(X) \geq 1$
there exists $0<r_{1}<r_{2}<r$ and a neighborhood $V$ of $0$ such that
${\rm exp}(B(0,C(r,R))X)(V \times B(0,r_{1})) \subset V \times B(0,r_{2})$ and
\[ |y - \gamma_{j}(x_{0})| \leq 2
|(y-\gamma_{j}(x_{0})) \circ \kappa_{t}|
\ \ \forall (y,t,j) \in B(0,r_{2}) \times [0,1] \times \{1,\hdots,N\} . \]
We obtain
\[ \left|{ \partial \kappa_{t} / \partial t }\right|(y)
\leq C(r,R) |X(y) \circ \kappa_{t}(y)|
\ \ \forall (y,t) \in (B(0,r_{2}) \setminus (Sing X)(x_{0})) \times [0,1]. \]
We deduce that  $|\psi^{X} \circ \kappa - \psi^{X}|(y) \leq C(r,R)$
for all $y \in B(0,r_{1}) \setminus (Sing X)(x_{0})$. By Riemann's theorem
$\psi^{X} \circ \kappa -\psi^{X}$ belongs to $\vartheta(B(0,r_{1}))$.
\end{proof}
The next results are important. Later on they will allow us to establish
the connection between the formal and analytic conjugations.
\begin{lem}
\label{lem:firsterms}
Let $Y \in {\mathcal X}\cn{}$. Consider an integral of the time form
$\psi$ of $Y$. Suppose that $\kappa \in \diff{}{}$ satisfies
that $\psi \circ \kappa - \psi$ belongs to ${\mathbb C}\{y\}$.
Then we have
$(\partial \kappa / \partial y)(0) =
e^{(\psi \circ \kappa - \psi)(0) (\partial Y(y) / \partial y)(0)}$.
Supposed $(\partial Y(y)/\partial y)(0)=0$ we also obtain
\[ \frac{\partial^{\nu(Y)+1} \kappa}
{\partial y^{\nu(Y)+1}}(0) = (\psi \circ \kappa - \psi)(0)
\frac{\partial^{\nu(Y)+1} Y(y)} {\partial y^{\nu(Y)+1}}(0) \]
and $(\partial^{j} \kappa/\partial y^{j})(0)=0$ for all $2 \leq j \leq \nu(Y)$.
\end{lem}
\begin{proof}
Denote $\lambda = (\partial Y(y)/\partial y)(0)$.
We have that $\psi \circ \kappa - \psi$ is of the form
$d+ L(y)$ for some $d \in {\mathbb C}$ and $L \in (y)$.
Suppose $\lambda \neq 0$.
Then $\psi$ is of the form $(\ln y) / \lambda + B(y)$
in the neighborhood of $0$ where  $B \in {\mathbb C}\{y\}$. Therefore we obtain
$d = (\ln (\partial \kappa / \partial y) (0))/\lambda$.
Suppose $\lambda = 0$. We obtain
$\kappa(y) = {\rm exp} ((d+t) Y(y)\partial / \partial y) (y,L(y))$.
This implies $\kappa(y) = y + d  Y(y) + O(y^{\nu(Y)+2})$.
The result is a consequence of last formula.
\end{proof}
Every $\phi \in \diff{}{}$ such that
$(\partial \phi/\partial y)(0)$ is not in
${e}^{2 \pi i{\mathbb Q}} \setminus \{1\}$ has a convergent
normal form. If the linear part is the identity is a consequence of
proposition \ref{pro:excofor}. Otherwise it is clear since
$\phi$ is formally linearizable.
\begin{cor}
\label{cor:conmpsi}
Let $\phi \in \diff{}{} \setminus \{ Id \}$ such that
$(\partial \phi/\partial y)(0) \not \in
{e}^{2 \pi i{\mathbb Q}} \setminus \{1\}$.
Consider a convergent normal form ${\rm exp}(Y)$ of $\phi$.
Let $\psi$ a holomorphic integral of the time form of $Y$.
Suppose that $\psi \circ \upsilon - \psi$ belongs to
${\mathbb C}\{y\} \cap (y)$ for some $\upsilon \in Z(\phi)$.
Then we have $\upsilon \equiv Id$.
\end{cor}
\begin{proof}
By lemma \ref{lem:firsterms} we have $j^{1} \upsilon \equiv Id$.
Moreover, if $(\partial \phi/\partial y)(0) \neq 1$ then
$\upsilon \equiv Id$ (prop. \ref{pro:comonevar}).
Suppose $(\partial \phi/\partial y)(0) = 1$, then we have
$y \circ \upsilon - y \in (y^{\nu(Y)+2})$ (lemma \ref{lem:firsterms}).
We obtain $\upsilon = \hat{\sigma}(\phi,\phi)=Id$.
\end{proof}
\begin{lem}
\label{lem:coeffc}
Let $\varphi_{1}, \varphi_{2} \in \diff{p1}{2}$ with common normal form
${\rm exp}(X)$. Fix $\gamma \equiv (y=\gamma_{1}(x)) \in Sing_{V} X$ and
$\hat{c} \in {\mathbb C}[[x]]$. Then we have
\[ \frac{\partial ({\rm exp}(\hat{c} \log \varphi_{2}) \circ
\hat{\sigma}(\varphi_{1},\varphi_{2},\gamma))}{\partial y}(x,\gamma_{1}(x))
\equiv
e^{\hat{c}(x) \frac{\partial X(y)}{\partial y}(x,\gamma_{1}(x))} . \]
Supposed $(\partial X(y)/\partial y)(x, \gamma_{1}(x)) \equiv 0$ we
also obtain
\[ \frac{\partial^{\nu_{X}(\gamma)+1}
({\rm exp}(\hat{c} \log \varphi_{2}) \circ
\hat{\sigma}(\varphi_{1},\varphi_{2},\gamma))}
{\partial y^{\nu_{X}(\gamma)+1}}(x,\gamma_{1}(x)) \equiv
\hat{c}(x) \frac{\partial^{\nu_{X}(\gamma)+1} X(y)}
{\partial y^{\nu_{X}(\gamma)+1}} (x,\gamma_{1}(x)) \]
and $(\partial^{j} ({\rm exp}(\hat{c} \log \varphi_{2}) \circ
\hat{\sigma}(\varphi_{1},\varphi_{2},\gamma))/\partial y^{j})
(x,\gamma_{1}(x)) \equiv 0$ for all  $2 \leq j \leq \nu_{X}(\gamma)$.
\end{lem}
\begin{proof}
Since $y \circ \hat{\sigma}(\varphi_{1},\varphi_{2},\gamma) -y
\in I(\gamma)^{\nu_{X}(\gamma)+2}$
then it is enough to prove the result for
${\rm exp}(\hat{c}(x) \log \varphi_{2})$.
We denote $\hat{X}=\log \varphi_{2}$, the equation
\[ \sum_{j=0}^{\infty} \frac{\hat{c}(x)^{j}}{j!}
\frac{\partial \hat{X}^{\circ (j)}(y)}{\partial y}(x,\gamma_{1}(x))
\equiv  \sum_{j=0}^{\infty} \frac{\hat{c}(x)^{j}}{j!}
\frac{\partial X(y)}{\partial y}(x,\gamma_{1}(x))^{j} \equiv
e^{\hat{c}(x) \frac{\partial X(y)}{\partial y}(x,\gamma_{1}(x))}  \]
implies the first part of the lemma. Suppose
$(\partial X(y)/\partial y)(x, \gamma_{1}(x)) \equiv 0$. Since
$\hat{X}(y)-X(y) \in (y \circ \varphi_{2}-y)^{2} \subset
(y-\gamma_{1}(x))^{2 \nu_{X}(\gamma)+2}$ then we obtain
\[ y \circ {\rm exp}(\hat{c}(x) \log \varphi_{2}) - y =
\hat{c}(x) X(y) + O((y-\gamma_{1}(x))^{\nu_{X}(\gamma)+2}) . \]
The rest of the proof is trivial.
\end{proof}
\subsection{Analytic classification and centralizer}
Let $\varphi_{1}, \varphi_{2} \in \diff{xp1}{2}$ with common convergent
normal form. Given a formal conjugation $\hat{\eta} \in \diffh{p}{2}$
we express the condition $\hat{\eta} \in \diff{}{2}$ in terms of the
changes of charts.
\begin{pro}
\label{pro:conimeq}
  Let $\varphi_{1}, \varphi_{2} \in \diff{xp1}{2}$ with common convergent
normal form
${\rm exp}(X)$. Fix $\mu \in e^{i(0,\pi)}$ and a compact connected set
$K_{X}^{\mu} \subset {\mathbb S}^{1} \setminus B_{X}^{\mu}$. Fix a privileged
curve $y =\gamma_{1}(x)$ associated to $X$. Consider a r-mapping $\kappa$
conjugating $(\varphi_{1})_{|x=x_{0}}$ and
$(\varphi_{2})_{|x=x_{0}}$. Then we have
\[ \xi_{\varphi_{2}, K_{X}^{\mu}}^{j} (x_{0},z) =
(z + c(x_{0})) \circ \xi_{\varphi_{1}, K_{X}^{\mu}}^{j}(x_{0},z)
\circ (z-c(x_{0}))
\ \ \forall j \in {\mathbb Z}/(2 \nu(X) {\mathbb Z}) \]
for all $x_{0} \in (0,\delta_{0}) K_{X}^{\mu}$
where $c(x_{0})=(\psi^{X} \circ \kappa - \psi^{X})(x_{0},\gamma_{1}(x_{0}))$.
\end{pro}
\begin{proof}
Suppose that $\kappa$ is a rR-mapping by taking
a smaller $0<r<\epsilon$ if necessary (lemma \ref{lem:quitR}).
Denote $X= u(x,y) \prod_{k=1}^{N} (y-\gamma_{k}(x))^{n_{k}} \partial /\partial y$
where $u \in {\mathbb C}\{x,y\}$ is a unit. Let $c_{1}^{l}, \hdots, c_{N}^{l}$
be the privileged functions associated to 
$(X,\varphi_{l},K_{X}^{\mu},\gamma_{1})$ for $l \in \{1,2\}$.
Consider the sections $T_{\mu X}^{\epsilon,1}$, $\hdots$,
$T_{\mu X}^{\epsilon, 2 \nu(X)}$.
Denote by $H(j)$ the unique element of
$Reg(\epsilon, \mu X, K_{X}^{\mu})$ such that
$T_{\mu X}^{\epsilon,j}(x) \in \partial H(j)(x)$ for all
$x \in [0,\delta_{0})K_{X}^{\mu}$. Let $0 < r_{1} <r$ and
$C(r,R) \in {\mathbb R}^{+}$ be the constants provided by lemma
\ref{lem:bddfldi}. We choose $r_{1}$ such that
${\rm exp}(B(0,C(r,R)) X)(|y|<r_{1}) \subset (|y|<\epsilon)$, we obtain
$\kappa(H(j)') \subset H(j)$ for all
$j \in {\mathbb Z}$ where $H(j)'$ is the element of
$Reg(r_{1}, \mu X, K_{X}^{\mu})$
contained in $H(j)$.

We define ${\phi}_{j}^{\varphi_{1}} =
\psi_{j}^{\varphi_{2}} \circ \kappa$ for $j \in {\mathbb Z}$.
Since
\[ {\phi}_{j}^{\varphi_{1}} - \psi_{j}^{X}=
(\psi_{j}^{\varphi_{2}}  - \psi_{j}^{X}) \circ \kappa +
(\psi_{j}^{X} \circ \kappa - \psi_{j}^{X}) \]
then ${\phi}_{j}^{\varphi_{1}} - \psi_{j}^{X}$ is continuous
in $H(j)'(x_{0}) \cup (\partial H(j)'(x_{0}) \cap Sing X)$ by
proposition \ref{pro:bddcon} and lemma \ref{lem:bddfldi}.
Therefore $({\phi}_{j}^{\varphi_{1}}-{\psi}_{j}^{\varphi_{1}})(x_{0},y)$
is continuous in $\partial H(j)'(x_{0}) \cap Sing X$
and then constant. Clearly $\phi_{j}^{\varphi_{1}}$ can be extended by
iteration to
a Fatou coordinate of $\varphi_{1}$  in $H(j)(x_{0})$.
We have that $\alpha^{\mu X}(H(j))$ and $\omega^{\mu X}(H(j))$ are
equal to  curves $y=\gamma_{k(j,\alpha)}(x)$ and
$y=\gamma_{k(j, \omega)}(x)$ respectively. We obtain
\[ \lim_{y \to \gamma_{k}(x_{0})} ({\phi}_{j}^{\varphi_{1}} -
\psi_{j}^{X})(x_{0},y)
= c_{k}^{2}(x_{0}) + (\psi^{X} \circ \kappa -\psi^{X})
(x_{0},\gamma_{k}(x_{0}))  \]
where $k \in \{k(j,\alpha), k(j,\omega) \}$. We deduce that
\[ \lim_{y \to \gamma_{k}(x_{0})} ({\phi}_{j}^{\varphi_{1}} -
\psi_{j}^{\varphi_{1}})(x_{0},y)
= c_{k}^{2}(x_{0}) - c_{k}^{1}(x_{0}) + (\psi^{X} \circ \kappa -\psi^{X})
(x_{0},\gamma_{k}(x_{0}))  \]
for $k \in \{k(j,\alpha), k(j,\omega) \}$. Since
$({\phi}_{j}^{\varphi_{1}}-{\psi}_{j}^{\varphi_{1}})(x_{0},y)$ is constant then
\[ c_{k(j,\upsilon)}^{2}(x_{0}) - c_{k(j,\upsilon)}^{1}(x_{0}) +
(\psi^{X} \circ \kappa -\psi^{X}) (x_{0},\gamma_{k(j,\upsilon)}(x_{0})) \]
does not depend on $\upsilon \in \{ \alpha, \omega \}$.
The graph ${\mathcal G}(\mu X, K_{X}^{\mu})$ is connected
(prop. \ref{pro:acyconn}), hence $c_{k}^{2}(x_{0}) - c_{k}^{1}(x_{0}) +
(\psi^{X} \circ \kappa -\psi^{X})(x_{0},\gamma_{k}(x_{0}))$
does not depend on $k \in \{1,\hdots,N\}$. In particular we obtain that
$({\phi}_{j}^{\varphi_{1}}-{\psi}_{j}^{\varphi_{1}})(x_{0},y)$ is equal to the
constant function $c(x_{0})$ for all
$j \in  {\mathbb Z}/(2 \nu(X) {\mathbb Z})$. By construction we get
\[ \xi_{\varphi_{2}, K_{X}^{\mu}}^{j}(x_{0},z) = \phi_{j+1}^{\varphi_{1}} \circ
{\left({ \phi_{j}^{\varphi_{1}} }\right)}^{\circ (-1)}(x_{0},z)
=(z + c(x_{0})) \circ  \xi_{\varphi_{1}, K_{X}^{\mu}}^{j}(x_{0},z) \circ
(z - c(x_{0})) \]
for all $j \in  {\mathbb Z}/(2 \nu(X) {\mathbb Z})$ as we wanted to prove.
\end{proof}
\begin{pro}
\label{pro:iinvic}
  Let $\varphi_{1}, \varphi_{2} \in \diff{xp1}{2}$ with common convergent
normal form
${\rm exp}(X)$. Fix $\mu \in e^{i(0,\pi)}$ and a compact connected set
$K_{X}^{\mu} \subset {\mathbb S}^{1} \setminus B_{X}^{\mu}$. Fix a privileged
curve $y=\gamma_{1}(x)$ associated to $X$ and a constant $M>0$. Suppose that
\[ \xi_{\varphi_{2}, K_{X}^{\mu}}^{j} (x_{0},z) =
(z + c(x_{0})) \circ \xi_{\varphi_{1}, K_{X}^{\mu}}^{j}(x_{0},z) \circ
(z-c(x_{0})) \ \ \forall j \in  {\mathbb Z}/(2 \nu(X) {\mathbb Z}) \]
for some  $x_{0} \in [0,\delta_{0}) K_{X}^{\mu}$ and $c(x_{0}) \in B(0,M)$.
Then there exists a r-mapping $\kappa$ such that
$\kappa \circ (\varphi_{1})_{|x=x_{0}} =(\varphi_{2})_{|x=x_{0}} \circ \kappa$.
The constant $r \in {\mathbb R}^{+}$ does not depend on $x_{0}$. Moreover we
get $(\psi^{X} \circ \kappa - \psi^{X})(x_{0},\gamma_{1}(x_{0}))=c(x_{0})$.
\end{pro}
\begin{proof}
Consider the notations in proposition \ref{pro:conimeq}.
We want to define
\[ \kappa(y) = { \left({ \psi_{j}^{\varphi_{2}} }\right) }^{\circ (-1)} \circ
(x_{0}, z + c(x_{0}))
\circ \psi_{j}^{\varphi_{1}}(x_{0},y) \]
for $j \in {\mathbb Z}$. There exists $A \in {\mathbb R}^{+}$ such that
$\sup_{H(j)} |\psi_{j}^{\varphi_{l}} - \psi_{j}^{X}| \leq A$ for
$l \in \{1,2\}$ (prop. \ref{pro:bddcon}). We have
${\rm exp}(B(2A+M)X)(|y|<R) \subset (|y|<\epsilon)$
for some $R \in {\mathbb R}^{+}$
Let $E$ be the union of the elements of $Reg(R,\mu X, K_{X}^{\mu})$.
We deduce that $\kappa$ is well-defined in $E(x_{0})$
and satisfies $\sup_{E(x_{0})} |\psi^{X} \circ \kappa - \psi^{X}| < 2A+M$,
in particular we have $\kappa(E(x_{0})) \subset B(0,\epsilon)$.
Denote $D= \max_{l \in \{1,2\}, s \in \{-1,1\}}
\sup_{B(0,R)} |\Delta_{\varphi_{l}^{\circ (s)}}|$.
There exist $0<r<R$ and $B \in {\mathbb N}$ such that
for all $J \in Reg_{\infty}(r, \mu X, K_{X}^{\mu})$ we have
\begin{itemize}
\item $\cup_{k \in \{-B,\hdots,B\}} \{ \varphi_{1}^{\circ (k)}(P) \}
\subset (y|<R)$
for all $P \in \overline{J} \setminus Sing X$.
\item $\exists 0 \leq k_{0},k_{1} \leq B$ such that
$\{ \varphi_{1}^{\circ (-k_{0})}(P), \varphi_{1}^{\circ (k_{1})}(P) \}
\subset E$ \ $\forall P \in \overline{J}  \setminus Sing X$.
\item ${\rm exp}((2A+M+2BD)X)(|y|<r) \subset (|y|<R)$.
\end{itemize}
We can define $\kappa$ in $\overline{J}(x_{0}) \setminus Sing X$ as
either
$\varphi_{2}^{\circ (k_{0})} \circ \kappa \circ \varphi_{1}^{\circ (-k_{0})}$
or
$\varphi_{2}^{\circ (-k_{1})} \circ \kappa \circ \varphi_{1}^{\circ (k_{1})}$.
By the construction and the hypothesis $\kappa$ is a well-defined holomorphic
mapping in $B(0,r) \setminus (Sing X)(x_{0})$ conjugating
$(\varphi_{1})_{|x=x_{0}}$ and $(\varphi_{2})_{|x=x_{0}}$. Moreover, we have
$sup_{B(0,r)} |\psi^{X} \circ \kappa - \psi^{X}|<2A+M+2BD$. As a consequence we
can extend $\kappa$ to $B(0,r)$ in a continuous (and then holomorphic) way by
defining $\kappa_{|(Sing X)(x_{0})} \equiv Id$.
The mapping $\kappa$ satisfies $\kappa(B(0,r)) \subset B(0,R)$.
Analogously by defining
\[ \kappa^{\circ (-1)}(y) =
{ \left({ \psi_{j}^{\varphi_{1}} }\right) }^{\circ (-1)} \circ
(x_{0}, z - c(x_{0}))
\circ \psi_{j}^{\varphi_{2}}(x_{0},y) \]
for $j \in {\mathbb Z}$ we obtain a mapping
$\kappa^{\circ (-1)}:B(0,r') \to B(0,R')$ conjugating
$(\varphi_{2})_{|x=x_{0}}$ and $(\varphi_{1})_{|x=x_{0}}$. By taking
$R \leq r'$ in the construction
of $\kappa$ we obtain that $\kappa$ is a rR-mapping.
\end{proof}
The next theorem is the analogue of proposition \ref{pro:ancottc}
in the non-trivial type case.
\begin{teo}
\label{teo:forcen}
  Let $\varphi_{1}, \varphi_{2} \in \diff{xp1}{2}$ with common convergent
normal form
${\rm exp}(X)$. Fix $\mu \in e^{i(0,\pi)}$ and a compact connected set
$K_{X}^{\mu} \subset {\mathbb S}^{1} \setminus B_{X}^{\mu}$.
Consider a privileged curve $\gamma \equiv (y=\gamma_{1}(x))$
in $Sing_{V}X$. Then $\varphi_{1} \sim \varphi_{2}$
if and only if there exists $d \in {\mathbb C}\{x\}$ such that
\[ \xi_{\varphi_{2},K_{X}^{\mu}}^{j}(x,z) \equiv
(z + d(x)) \circ  \xi_{\varphi_{1},K_{X}^{\mu}}^{j} \circ
(x,z-d(x)) \ \  \forall j \in  {\mathbb Z}/(2 \nu(X) {\mathbb Z}). \]
The previous equation is equivalent to
${\rm exp}(d(x) \log \varphi_{2}) \circ
\hat{\sigma}(\varphi_{1},\varphi_{2},\gamma) \in \diff{}{2}$.
\end{teo}
\begin{proof}
Implication $\Rightarrow$. Let $\sigma \in \diff{p}{2}$ conjugating
$\varphi_{1}$ and $\varphi_{2}$. We denote
$c(x) \equiv (\psi^{X} \circ \sigma - \psi^{X})(x,\gamma_{1}(x))$,
we have $c \in {\mathbb C}\{x\}$ (lemma \ref{lem:bddfldi}). We deduce that
\[ \xi_{\varphi_{2}, K_{X}^{\mu}}^{j} (x,z) \equiv
(z + c(x)) \circ \xi_{\varphi_{1}, K_{X}^{\mu}}^{j}(x,z)
\circ (x,z-c(x)) \ \ \forall j \in  {\mathbb Z}/(2 \nu(X) {\mathbb Z}) \]
by proposition \ref{pro:conimeq}.
The mapping $\sigma$ is of the form
${\rm exp}(\hat{c}(x) \log \varphi_{2}) \circ
\hat{\sigma}(\varphi_{1},\varphi_{2},\gamma)$ (lemma \ref{lem:comtwodif}).
Lemmas \ref{lem:firsterms} and \ref{lem:coeffc} imply
$\hat{c} \equiv c$.

Implication $\Leftarrow$.
Fix an EV-covering $K_{X}^{\mu_{1}}=K_{X}^{\mu}$, $K_{X}^{\mu_{2}}$,
$\hdots$, $K_{X}^{\mu_{l}}$. Supposed
\begin{equation}
\label{equ:forcen}
\xi_{\varphi_{2}, K_{X}^{\mu_{p}}}^{j} (x,z) \equiv
(z + d(x)) \circ \xi_{\varphi_{1}, K_{X}^{\mu_{p}}}^{j}(x,z)
\circ (x,z-d(x)) \ \ \forall j \in  {\mathbb Z}/(2 \nu(X) {\mathbb Z})
\end{equation}
for some $p \in \{1, \hdots , l\}$ the proof of prop. \ref{pro:iinvic}
provides a continuous mapping $\sigma_{p}(x,y)$
in the set $[0,\delta_{0})K_{X}^{\mu} \times B(0,r)$ such that it
is holomorphic in the interior and conjugates $\varphi_{1}$ and $\varphi_{2}$.
Moreover $\sigma_{p}(x_{0},y)$ is a rR-mapping for all
$x_{0} \in [0,\delta_{0})K_{X}^{\mu}$ and some $r,R \in {\mathbb R}^{+}$.
We obtain
$(\psi^{X} \circ \sigma_{p} - \psi^{X})(x,\gamma_{1}(x)) \equiv d(x)$.

The existence of $\sigma_{1}$ and proposition \ref{pro:conimeq} imply that
\[ \xi_{\varphi_{2},K_{X}^{\mu_{q}}}^{j}(x_{0},z) \equiv
(z + d(x_{0})) \circ  \xi_{\varphi_{1},K_{X}^{\mu_{q}}}^{j}(x_{0},z)
\circ  (z-d(x_{0})) \]
for all $j \in {\mathbb Z}$ and for all
$x_{0} \in (0,\delta_{0}) (\dot{K}_{X}^{\mu_{1}} \cap \dot{K}_{X}^{\mu_{q}})$.
By analytic continuation we obtain the same result for
$x_{0} \in [0,\delta_{0}) \dot{K}_{X}^{\mu_{q}}$ if
$\dot{K}_{X}^{\mu_{1}} \cap \dot{K}_{X}^{\mu_{q}} \neq \emptyset$.
The iteration of this process shows that the equation \ref{equ:forcen}
is fulfilled for all $q \in \{1,\hdots,l\}$ and
$x_{0} \in [0,\delta_{0}) K_{X}^{\mu_{q}}$.

Suppose $\dot{K}_{X}^{\mu_{p}} \cap \dot{K}_{X}^{\mu_{q}} \neq \emptyset$
for $p,q \in \{1,\hdots,l\}$.
Denote $h = (\sigma_{q})^{\circ (-1)} \circ \sigma_{p}$. We obtain
$h \circ \varphi_{1} = \varphi_{1} \circ h$ in
$x \in [0,\delta_{0}) (\dot{K}_{X}^{\mu_{p}} \cap \dot{K}_{X}^{\mu_{q}})$ and
$(\psi^{X} \circ h - \psi^{X})(x,\gamma_{1}(x)) \equiv 0$.
Corollary \ref{cor:conmpsi} implies
$h(x,y) \equiv Id$ and then
$\sigma_{p} \equiv \sigma_{q}$ in
$[0,\delta_{0})  (\dot{K}_{X}^{\mu_{p}} \cap \dot{K}_{X}^{\mu_{q}})
\times B(0,r)$. Thus all the $\sigma_{b}$ ($b \in \{1,\hdots,l\}$)
paste together in a mapping $\sigma$ such that it is
continuous in $B(0,\delta_{0}) \times B(0,r)$
and holomorphic in $(B(0,\delta_{0}) \setminus \{0\}) \times B(0,r)$.
By Riemann's theorem $\sigma$ is an element of $\diff{p}{2}$
conjugating $\varphi_{1}$ and $\varphi_{2}$.
Moreover we have
$\sigma = {\rm exp}(d(x) \log \varphi_{2}) \circ
\hat{\sigma}(\varphi_{1},\varphi_{2},\gamma)$ by the first part of the
proof.
%
%
\end{proof}
%
\begin{pro}
\label{pro:ratcen}
Let $\varphi \in \diff{p1}{2}$ such that
$\log \varphi \not \in \Xnt$ and $Fix \varphi$ is not of trivial type.
Then there exists $q \in {\mathbb N}$
such that $Z(\varphi) = <{\rm exp}(q^{-1} \log \varphi)>$.
\end{pro}
\begin{proof}
We can suppose $\varphi \in \diff{xp1}{2}$ up to a ramification
$(x^{k},y)$. Let ${\rm exp}(X)$ be a convergent normal form of $\varphi$.
A diffeomorphism $\eta \in Z(\varphi)$ is of the form
${\rm exp}(c(x) \log \varphi)$ by lemma \ref{lem:comtwodif}.
Consider a privileged $y=\gamma_{1}(x)$ in $Sing_{V} X$.
We have $(\psi^{X} \circ \eta - \psi^{X})(x,\gamma_{1}(x)) \equiv c(x)$
by lemmas \ref{lem:firsterms} and \ref{lem:coeffc}.
Fix $\mu \in e^{i(0,\pi)}$ and a compact connected set
$K_{X}^{\mu} \subset {\mathbb S}^{1} \setminus B_{X}^{\mu}$.
Denote
$E=\{ l \in {\mathbb N} : \exists j \in {\mathbb Z} \ s.t.\
a_{j,l,K_{X}^{\mu}}^{\varphi} \not \equiv 0\}$.
The set $E$ is not empty (prop. \ref{pro:nembflo}). Denote $q = \gcd E$.
The continuous functions $c(x)$  satisfying the equation
\[ \xi_{\varphi, K_{X}^{\mu}}^{j} (x,z) =
(z + c(x)) \circ \xi_{\varphi, K_{X}^{\mu}}^{j}(x,z) \circ (x,z-c(x))
\ \ \forall j \in  {\mathbb Z}/(2 \nu(\varphi) {\mathbb Z}) . \]
are the constant functions of the form $p/q$ for some
$p \in {\mathbb Z}$. Thus the result is a consequence of theorem
\ref{teo:forcen}.
\end{proof}
\subsection{Complete system of analytic invariants}
We can introduce a complete system of analytic invariants for
elements $\varphi \in \diff{p1}{2}$. The presentation is slightly
simpler if $\varphi_{|x=0}$ is not analytically trivial. In such a
case we obtain the generalization of Mardesic-Roussarie-Rousseau's
invariants.

Let $\varphi_{1}, \varphi_{2} \in \diff{p1}{2}$ with
$Fix \varphi_{1}=Fix \varphi_{2}$ and
$Res(\varphi_{1}) \equiv Res(\varphi_{2})$.
Suppose that $Fix \varphi_{1}$ is not of trivial type.
Let ${\rm exp}(X)$ be a convergent normal form of $\varphi_{1}$.
There exists $k \in {\mathbb N}$ such that $Y=(x^{k},y)^{*} X$ belongs to
$\Xt$. Fix a privileged curve $\gamma \in Sing_{V} Y$. Consider an
EV-covering $K_{1}=K_{Y}^{\mu_{1}}$, $\hdots$, $K_{l}=K_{Y}^{\mu_{l}}$.
We say that $m_{\varphi_{1}}(x_{0}) = m_{\varphi_{2}}(x_{0})$ for
$x_{0}$ in $B(0,\delta_{0}) \setminus \{0\}$ if there exist
$c(x_{0}) \in {\mathbb C}$ and $b(x_{0}) \in \{1,\hdots,l\}$ such that
$x_{0} \in {\mathbb R}^{+} \dot{K}_{b(x_{0})}$ and
\begin{equation}
\label{equ:rigid}
\xi_{\varphi_{2},K_{b(x_{0})}}^{j}(x_{0},z) \equiv
(z + c(x_{0})) \circ  \xi_{\varphi_{1},K_{b(x_{0})}}^{j}(x_{0},z) \circ
(x_{0},z-c(x_{0})) \ \ \forall j \in {\mathbb Z}.
\end{equation}
The definition makes sense since an EV-covering depends only on
$Fix \varphi$ and $Res(\varphi)$ for $\varphi \in \diff{tp1}{2}$
by remark \ref{rem:evndepnf}.
We denote $m_{\varphi_{1}}(0) = m_{\varphi_{2}}(0)$ if we have
$(\varphi_{1})_{|x=0} \sim (\varphi_{2})_{|x=0}$. We say that
$Inv(\varphi_{1}) \sim Inv(\varphi_{2})$ if
$m_{\varphi_{1}}(x_{0}) = m_{\varphi_{2}}(x_{0})$ for all
$x_{0}$ in a pointed neighborhood of $0$ and we can choose
$c:B(0,\delta_{0}) \setminus \{0\} \to {\mathbb C}$ such that
$Img(c)$ is bounded. Both invariants $m_{\varphi}$ and $Inv(\varphi)$
can be expressed in terms of $\mu$-spaces of orbits.
In this section we prove that
$\varphi_{1} \sim \varphi_{2}$ is equivalent to
$Inv(\varphi_{1}) \sim Inv(\varphi_{2})$.
%
%
\begin{lem}
\label{lem:nomon}
Let $f(x)$ be a multi-valuated holomorphic
function of $B(0,\delta) \setminus \{0\}$ such that
$f(e^{2 \pi i}x) -f(x) \equiv C$ for some $C \in {\mathbb R}$.
Suppose that $|Img f(x)|$ is bounded in a neighborhood of $0$.
Then $f$ belongs to $\vartheta(B(0,\delta))$.
\end{lem}
\begin{proof}
We define $F=f(x)-(C/2 \pi i) \ln x$, we obtain
$F \in \vartheta(B(0,\delta) \setminus \{0\})$.
Moreover we have $Img F = Img f +(C/2 \pi) \ln |x|$.
Suppose $C=0$, then $f$ has a removable singularity at $x=0$
since $Img f$ is bounded.

Let us prove that $C \neq 0$ is not possible.
Since $\lim_{x \to 0} Img F \in \{ - \infty, +\infty \}$
then $F$ does not have an essential singularity.
Supposed $F$ has a pole at $x=0$ then it is
of the form $A e^{i \theta}/x^{l}+O(1/x^{l-1})$ for some
$(l,A,\theta) \in {\mathbb N} \times {\mathbb R}^{+} \times {\mathbb R}$.
This is not possible since then
$\lim_{r \to 0} Img f(r e^{\frac{i(\theta- \pi/2)}{l}}) = \infty$.
Finally we obtain $C=0$ since we have
$\lim_{x \to 0} Img f(x)= Img F(0) -(C/2 \pi) \lim_{x \to 0} \ln |x|$.
\end{proof}
All the elements of $\diff{p1}{2}$ can be interpreted as elements of
$\diff{tp1}{2}$ up to a ramification $(x^{m},y)$. The ramification
preserves the analytic classes of elements of $\diff{p1}{2}$.
\begin{lem}
\label{lem:prepram}
Let $\varphi_{1},\varphi_{2} \in \diff{p1}{2}$ with
$Fix \varphi_{1} = Fix \varphi_{2}$. Consider $m \in {\mathbb N}$.
Then $\varphi_{1} \sim \varphi_{2}$ if and only if
$(x^{1/m},y) \circ \varphi_{1} \circ (x^{m},y) \sim
(x^{1/m},y) \circ \varphi_{2} \circ (x^{m},y)$.
\end{lem}
\begin{proof}
The sufficient condition is obvious.

  Denote $\tilde{\varphi}_{j} =
(x^{1/m},y) \circ \varphi_{j} \circ (x^{m},y)$ for $j \in \{1,2\}$.
We have $Fix \varphi_{1}=Fix \varphi_{2}$ by hypothesis and
$Res(\varphi_{1}) \equiv Res(\varphi_{2})$
since the residues are analytic invariants. We can suppose that
$\varphi_{1}$ and $\varphi_{2}$ are not analytically trivial.
Otherwise both $\log \varphi_{1}$ or $\log \varphi_{2}$ belong to
${\mathcal X} \cn{2}$, we obtain
$\varphi_{1} \sim \varphi_{2}$ by proposition \ref{pro:anconfl}.

Denote $h=(e^{2 \pi i/m} x,y)$.
Let $\sigma_{0} \in \diff{p}{2}$ conjugating
$\tilde{\varphi}_{1}$, $\tilde{\varphi}_{2}$. Since we have
$h^{\circ (-1)} \circ \tilde{\varphi}_{j} \circ h
= \tilde{\varphi}_{j}$ for $j \in \{1,2\}$ then
$\sigma_{k} = h^{\circ (-k)} \circ {\sigma}_{0} \circ h^{\circ (k)}$
conjugates $\tilde{\varphi}_{1}$ and $\tilde{\varphi}_{2}$
for $k \in \{0,\hdots,m\}$. The diffeomorphism
$\sigma_{0}^{\circ (-1)} \circ \sigma_{1}$ belongs to
$Z_{up}(\tilde{\varphi}_{1})$, hence it is of the form
${\rm exp}(C \log \tilde{\varphi}_{1})$ for some
$C \in {\mathbb Q}$ by propositions \ref{pro:ratcen} and \ref{pro:ancottc}.
Since $\sigma_{k}^{\circ (-1)} \circ \sigma_{k+1}$
is equal to
$h^{\circ (-k)} \circ {\rm exp}(C \log \tilde{\varphi}_{1})
\circ h^{\circ (k)} =  {\rm exp}(C \log \tilde{\varphi}_{1})$ then
\[ Id = (\sigma_{0}^{\circ (-1)} \circ \sigma_{1}) \circ
(\sigma_{1}^{\circ (-1)} \circ \sigma_{2}) \circ \hdots \circ
(\sigma_{m-1}^{\circ (-1)} \circ \sigma_{m}) =
{\rm exp}(C m \log \tilde{\varphi}_{1}) . \]
We obtain $C=0$ by uniqueness of the infinitesimal generator.
Since $\sigma_{0}$ and  $(e^{2 \pi i/m}x,y)$ commute
we deduce that $\sigma = (x^{m},y) \circ \sigma_{0} \circ (x^{1/m},y)$
is an element of $\diff{p}{2}$ conjugating $\varphi_{1}$ and
$\varphi_{2}$.
\end{proof}
We can prove now that $Inv$ provides a complete system of analytic invariants.
\begin{teo}
\label{teo:imgbdd}
 Let $\varphi_{1}, \varphi_{2} \in \diff{p1}{2}$. Suppose that
$Fix \varphi_{1}= Fix \varphi_{2}$ and
$Res \varphi_{1} \equiv Res \varphi_{2}$. Then
$\varphi_{1} \sim \varphi_{2}$ is equivalent to
$Inv(\varphi_{1}) \sim Inv(\varphi_{2})$.
\end{teo}
\begin{proof}
We can suppose that $Fix \varphi_{1}$ is not of trivial type by
proposition \ref{pro:anainv}.
We consider the notations at the beginning of this section.
We can suppose that $\log \varphi_{1}$ and $\log \varphi_{2}$ are divergent,
otherwise we have that $\varphi_{1} \sim \varphi_{2}$
(prop. \ref{pro:anconfl}) and we can choose $c \equiv 0$.
Let $\alpha_{j}$ be a convergent normal form of $\varphi_{j}$ for
$j \in \{1,2\}$. There exists a mapping $\sigma_{0}$ conjugating
$\alpha_{1}$ and $\alpha_{2}$ (prop. \ref{pro:anconfl}). Up to
replace $\varphi_{2}$ with
$\sigma_{0}^{\circ (-1)} \circ \varphi_{2} \circ \sigma_{0}$ and
$\xi_{\varphi_{2},K_{b}}^{j}$ with
$(z-d(x)) \circ \xi_{\varphi_{2},K_{b}}^{j} \circ (x,z+d(x))$ for all
$(b,j) \in \{1,\hdots,l\} \times {\mathbb Z}$ and some
$d \in {\mathbb C}\{x\}$ we can suppose that $\varphi_{1}$ and $\varphi_{2}$
have common convergent normal form. Finally we can suppose that
$\varphi_{1}, \varphi_{2} \in \diff{tp1}{2}$ by lemma \ref{lem:prepram}.

The sufficient condition is a consequence of theorem \ref{teo:forcen}.
Since change of charts commute with $z \to z+1$ we can suppose that
$c$ is bounded by replacing $c(x)$ with $c(x) - [Re(c(x))]$ where
$[]$ is the integer part. There exists a $r$-mapping
conjugating $\varphi_{1}(x_{0},y)$ and $\varphi_{2}(x_{0},y)$
for all $x_{0}$ in a pointed neighborhood of $0$ and some
$r \in {\mathbb R}^{+}$ by proposition \ref{pro:iinvic}. We obtain
\[ \xi_{\varphi_{2},K_{b}}^{j}(x_{0},z) \equiv
(z + c(x_{0})) \circ  \xi_{\varphi_{1},K_{b}}^{j}(x_{0},z)
\circ  (z-c(x_{0})) \ \ \forall j \in {\mathbb Z} \ \
\forall b \in \{1,\hdots,l\} \]
for all $x_{0} \in (0,\delta_{0}) \dot{K}_{b}$ by proposition
\ref{pro:conimeq}.

Suppose $\sup_{B(0,\delta_{0}) \setminus \{0\}} |Img \ c| < M$.
Fix $p \in \{1,\hdots,l\}$. Consider the set
\[ E_{s}^{p}(\varphi_{1})=\{(j,m) \in D_{s}(\varphi_{1}) \times {\mathbb N} :
a_{j,m,K_{p}}^{\varphi_{1}} \not \equiv 0 \} . \]
We define
$E^{p}(\varphi_{1})=E_{-1}^{p}(\varphi_{1}) \cup E_{1}^{p}(\varphi_{1})$.
We have $E^{p}(\varphi_{1}) \neq \emptyset$ by proposition \ref{pro:nembflo}.
Let $x_{1} \in (0,\delta_{0}) \dot{K}_{X}^{\mu_{p}}$
such that $(j,m) \in E^{p}(\varphi_{1})$ implies
$a_{j,m,K_{p}}^{\varphi_{1}}(x_{1}) \neq 0$. We define
\[ d_{j,m} = \frac{1}{2 \pi i m s} \ln
\frac{a_{j,m,K_{p}}^{\varphi_{2}}}
{a_{j,m,K_{p}}^{\varphi_{1}}} \ s.t. \ d_{j,m}(x_{1})=c(x_{1}) \]
for all $(j,m) \in E_{s}^{p}(\varphi_{1})$. Since
$e^{-2 \pi m M} \leq |a_{j,m,K_{p}}^{\varphi_{2}} /
a_{j,m,K_{p}}^{\varphi_{1}} | \leq e^{2 \pi m M}$
in $(0,\delta_{0})K_{p}$
then $d_{j,m} \in \vartheta((0,\delta_{0})\dot{K}_{p})$
for all $(j,m) \in E^{p}(\varphi_{1})$. We get
$d_{j,m}(x_{0})-c(x_{0}) \in {\mathbb Z}/m$ for
$(j,m) \in E^{p}(\varphi_{1})$ and
$a_{j,m,K_{p}}^{\varphi_{1}}(x_{0}) \neq 0$.
Thus the image of $d_{j,m} - d_{j',m'}$ is contained in
${\mathbb Z}/m + {\mathbb Z}/m'$ for
$(j,m), (j',m') \in E^{p}(\varphi_{1})$;
since $d_{j,m}(x_{1})=d_{j',m'}(x_{1})$
we deduce that $d_{j,m} \equiv d_{j',m'}$.
Denote by $d_{p}$ any function
$d_{j,m}$ for $(j,m) \in E^{p}(\varphi_{1})$.
We obtain
\[ \xi_{\varphi_{2},K_{X}^{\mu_{p}}}^{j}(x_{0},z) \equiv
(z + d_{p}(x_{0})) \circ  \xi_{\varphi_{1},K_{X}^{\mu_{p}}}^{j}(x_{0},z)
\circ  (z-d_{p}(x_{0})) \]
by construction
for all $(j,x_{0}) \in {\mathbb Z} \times (0,\delta_{0}) \dot{K}_{p}$.
We have $|Img (d_{p})| \leq M$ in $(0,\delta_{0}) \dot{K}_{p}$.

Consider $p,q \in \{1,\hdots,l\}$ such that
$\dot{K}_{p} \cap \dot{K}_{q} \neq \emptyset$.
Consider $(j,m) \in E^{p}(\varphi_{1})$ and
$(j',m') \in E^{q}(\varphi_{1})$. We have
$d_{p}(x_{0}) - c(x_{0}) \in {\mathbb Z}/m$ and
$d_{q}(x_{0}) - c(x_{0}) \in {\mathbb Z}/m'$
for all $x_{0} \in (0,\delta_{0})
(\dot{K}_{p} \cap \dot{K}_{q})$ such that
$(a_{j,m,K_{p}}^{\varphi_{1}}
a_{j',m',K_{q}}^{\varphi_{1}})(x_{0}) \neq 0$.
We deduce that $d_{p} - d_{q}$ is a constant function, moreover
$d_{p} - d_{q} \in {\mathbb Q}$. Then we can extend $d_{p}$ to
$(0,\delta_{0})(\dot{K}_{p} \cup \dot{K}_{q})$.
We get that $d_{1}$ is a multi-valuated
function in $B(0,\delta_{0}) \setminus \{0\}$ such that
$d_{1}(e^{2 \pi i}x)-d_{1}(x) \equiv C$ for some $C \in {\mathbb Q}$.
We also have $|Img(d_{1})| \leq M$ in $B(0,\delta_{0}) \setminus \{0\}$
and then $d_{1} \in \vartheta(B(0,\delta_{0}))$ by lemma \ref{lem:nomon}.
Then $\varphi_{1}$ and $\varphi_{2}$ are conjugated by an element
of $\diff{p}{2}$ by theorem \ref{teo:forcen}.
%
%
\end{proof}
We give now a geometrical interpretation of our complete system of
analytic invariants. Roughly speaking, given  $\varphi \in \diff{p1}{2}$
the next theorem claims that the analytic classes of $\varphi_{|x=x_{0}}$
for $x_{0} \in B(0,\delta_{0}) \setminus \{0\}$
characterize the analytic class of $\varphi$ whenever we exclude
singularities of the conjugating mappings at $x_{0}=0$. The result is the
generalization of proposition \ref{pro:mod}.
\begin{teo}
\label{teo:mod}
Let $\varphi_{1}, \varphi_{2} \in \diff{p1}{2}$ with
$Fix \varphi_{1}=Fix \varphi_{2}$. Then
$\varphi_{1} \sim \varphi_{2}$ if and only if
$(\varphi_{1})_{|x=x_{0}}$ and $(\varphi_{2})_{|x=x_{0}}$
are conjugated by a r-mapping $\kappa_{x_{0}}$ for some
$r \in {\mathbb R}^{+}$ and all $x_{0}$ in a pointed neighborhood of $0$.
\end{teo}
\begin{proof}
By proposition \ref{pro:mod} we can suppose that $Fix \varphi_{1}$ is
not of trivial type.

We have $Fix \varphi_{1}=Fix \varphi_{2}$ by hypothesis and
$Res(\varphi_{1}) \equiv Res(\varphi_{2})$
since the residues are analytic invariants.
Let $\alpha_{j}$ be a convergent normal form of $\varphi_{j}$ for
$j \in \{1,2\}$. Then there exists
$\zeta \in \diff{p}{2}$ such that
$\zeta \circ \alpha_{1} = \alpha_{2} \circ \zeta$
by proposition \ref{pro:anconfl}.
By replacing $\varphi_{2}$ with
$\zeta^{\circ (-1)} \circ \varphi_{2} \circ \zeta$
we can suppose that $\varphi_{1}$ and $\varphi_{2}$ have a common normal
form $\alpha_{1}$.
The mapping $\kappa_{x_{0}}$ has to be replaced with
$(\zeta^{\circ (-1)})_{|x=x_{0}} \circ \kappa_{x_{0}}$, it is
still a rR-mapping (maybe for a smaller $r \in {\mathbb R}^{+}$)
by lemma \ref{lem:quitR} for all $x_{0}$ in a pointed neighborhood of $0$.

There exists $m \in {\mathbb N}$
such that $X = (x^{m},y)^{*} \log \alpha_{1}$ belongs to $\Xt$.
Fix a privileged curve
$(y=\gamma_{1}(x)) \in Sing_{V} X$ and an
EV-covering. Let us denote
$c(x_{0})=
(\psi^{X} \circ \kappa_{x_{0}^{m}} - \psi^{X})(x_{0},\gamma_{1}(x_{0}))$.
We are done since proposition \ref{pro:conimeq} and lemma
\ref{lem:bddfldi} assure that the hypothesis of theorem \ref{teo:imgbdd}
is satisfied.
\end{proof}
\begin{teo}
\label{teo:rig}
Let $\varphi_{1}, \varphi_{2} \in \diff{p1}{2}$ satisfying that
$Fix \varphi_{1}=Fix \varphi_{2}$ and
$Res(\varphi_{1}) \equiv Res(\varphi_{2})$. Suppose that
$(\varphi_{1})_{|x=0} \in \diff{1}{}$ is
not analytically trivial. Then $\varphi_{1} \sim \varphi_{2}$
if and only if $m_{\varphi_{1}} \equiv m_{\varphi_{2}}$.
\end{teo}
The analogue of this theorem for the generic case when $N(X)=2$
is the main theorem in \cite{MRR}.  They do not impose any condition
on $(\varphi_{1})_{|x=0}$. The next section provides counterexamples
if $(\varphi_{1})_{|x=0}$ is analytically trivial. They did not notice
that the hypothesis $m_{\varphi_{1}} \equiv m_{\varphi_{2}}$ does not
prevent the degeneration of conjugations in the neighborhood of $x=0$.
\begin{proof}
We can suppose that $Fix \varphi_{1}$ is not of trivial type by
corollary \ref{cor:rig}. Moreover we can suppose that $\varphi_{1}$
and $\varphi_{2}$ have a common convergent normal form.
Consider the notations at the beginning of this section.

%
%

We have $\xi_{\varphi, K_{b}}^{j}(0,z) =
\xi_{\varphi(0,y)}^{\Lambda(j)}(z)$ for all
$\varphi \in \{ \varphi_{1}, \varphi_{2} \}$,
$b \in \{1,\hdots,l\}$ and $j \in {\mathbb Z}$
where $\Lambda \equiv \Lambda(\varphi_{1}) \equiv \Lambda(\varphi_{2})$
(cor. \ref{cor:EV0}).
Since $(\varphi_{1})_{|x=0}$ is not analytically trivial then
there exists $s(0) \in \{-1,1\}$ and
$(j(0),b(0),\beta) \in
D_{s(0)}(\varphi_{1})  \times {\mathbb N} \times {\mathbb C} \setminus \{0\}$
such that
$a_{j(0),b(0),K_{p}}^{\varphi_{1}}(0) = \beta$ for all $p \in \{1,\hdots,l\}$.
Then $m_{\varphi_{1}}(0) = m_{\varphi_{2}}(0)$ implies
that there exists $(j(1),\beta') \in
D_{s(0)}(\varphi_{1}) \times {\mathbb C} \setminus \{0\}$ such that
$a_{j(1),b(0),K_{p}}^{\varphi_{2}}(0) = \beta'$ for all $p \in \{1,\hdots,l\}$.
Since $m_{\varphi_{1}} \equiv m_{\varphi_{2}}$ we have
\[ \left\{
{\begin{array}{c}
a_{j(0),b(0),K_{b(x)}}^{\varphi_{2}}(x) =
a_{j(0),b(0),K_{b(x)}}^{\varphi_{1}}(x) e^{2 \pi i s(0) b(0) c(x)} \\
a_{j(1),b(0),K_{b(x)}}^{\varphi_{2}}(x) =
a_{j(1),b(0),K_{b(x)}}^{\varphi_{1}}(x) e^{2 \pi i s(0) b(0) c(x)} .
\end{array} }\right. \]
The first equation implies $s(0) Img c(x) > K_{1}$
in a pointed neighborhood of $0$ for some
$K_{1} \in {\mathbb R}$. We obtain
$s(0) Img c(x) < K_{2}$ for $x \neq 0$ and some $K_{2} \in {\mathbb R}$
from the second equation.
This implies $|Img \ c(x)| \leq \max (|K_{1}|,|K_{2}|)$ for all
$x \neq 0$ in a neighborhood of $0$. Now
$\varphi_{1} \sim \varphi_{2}$ is a consequence of
theorem \ref{teo:imgbdd}.
\end{proof}
\section{Optimality of the results}
\label{sec:otimo}
We introduce an example which proves that the hypothesis on the
non-analytical triviality of $(\varphi_{1})_{|x=0}$ in theorem
\ref{teo:rig} can not be dropped. It also shows that the uniform
hypothesis in theorem \ref{teo:mod} is essential.
\begin{pro}
\label{pro:couMRR}
Let  $X \in \Xnt$. There exist
$\varphi_{1}, \varphi_{2} \in \diff{p1}{2}$ with normal form
${\rm exp}(X)$ and such that $m_{\varphi_{1}} \equiv m_{\varphi_{2}}$
but $\varphi_{1} \not \sim \varphi_{2}$.
Moreover there exists an analytic injective mapping $\sigma$
conjugating $\varphi_{1}$ and $\varphi_{2}$ and defined
in a domain $|y|< C_{0}/\sqrt[\nu(X)]{|\ln x|}$ for some
$C_{0} \in {\mathbb R}^{+}$.
\end{pro}
In particular we provide a counter-example to the main theorem in
\cite{MRR}.
The domain $|y|< C_{0}/\sqrt[\nu(X)]{|\ln x|}$ is defined in the universal
covering of ${\mathbb C}^{*} \times {\mathbb C}$; its size
decays when $x$ tends to $0$.
Anyway the decay is slower than algebraic.

Let $X=f(x,y) \partial /\partial y \in \Xnt$.
We consider vector fields of the form
\[ X_{v}= \frac{f(x,y)}{1+ f(x,y) v(x,y,t)} \frac{\partial}{\partial y} +
2 \pi i t \frac{\partial}{\partial t} \]
where $v$ is defined in a domain of the form
$B(0,\delta) \times B(0,\epsilon) \times B(0,2)$ in coordinates
$(x,y,t)$. The vector field $X_{v}$ supports a dimension $1$ foliation
$\Omega_{v}$
preserving the hypersurfaces $x=cte$. Moreover since
$X_{v}(t)=2 \pi i t$ then $X_{v}$ is transversal to
every hypersurface $t=cte$ except $t=0$. As a consequence we can
consider the holonomy mapping $hol_{v}(x,y,t_{0},z_{0})$ of the foliation
given by $X_{v}$ along a path
$t \in e^{2 \pi i[0,z_{0}]} t_{0}$, it maps the transversal
$t=t_{0}$ to $t=t_{0} e^{2 \pi i z_{0}}$ for $t_{0} \neq 0$.
The restriction of $hol_{v}(x,y,t,z)$ to $(x,y) \in  Sing X$ is the identity.
Supposed that $v=v(x,y)$ we have
\[ hol_{v}(x,y,t,z) = \left({ {\rm exp}
\left({ z \frac{f(x,y)}{1+f(x,y) v(x,y)} }\right)(x,y),e^{2 \pi i z} t
}\right) . \]
The restriction $(\Omega_{v})_{|x=0}$ is a germ of saddle-node
for $v \in {\mathbb C}\{x,y,t\}$.
The holonomy $hol_{v}(0,y,t_{0},1)$
at a transversal $t=t_{0}$ to the strong integral
curve $y=0$ is analytically trivial if
and only if $(\Omega_{v})_{|x=0}$ is analytically normalizable
\cite{MaRa:ihes}. In particular $(\Omega_{0})_{|x=0}$ is analytically
normalizable. Every foliation in the same formal class than
$(\Omega_{0})_{|x=0}$ is analytically conjugated to some
$(\Omega_{v})_{|x=0}$ with $v \in {\mathbb C}\{y,t\} \cap (y,t)$, 
we just truncate the formal conjugation. Every formal
class contains non-analytically normalizable elements, hence there exists
$v^{0} \in {\mathbb C}\{y,t\} \cap (y,t)$ such that 
\[ (X_{v^{0}})_{|x=0} = \frac{f(0,y)}{1+f(0,y) v^{0}(y,t)}
\frac{\partial}{\partial y} + 2 \pi i t \frac{\partial}{\partial t} \]
is not analytically normalizable. Hence the holonomy
$hol_{v^{0}}(0,y,t_{0},1)$ is not analytically trivial for $t_{0} \neq 0$.
Moreover up to change of coordinates $(x,y,t) \to (x,y,\eta t)$
for some $\eta \in {\mathbb R}^{+}$ there exists
$(\delta_{0}, \epsilon_{0}) \in {\mathbb R}^{+}$ such that
\begin{itemize}
\item $v^{0} \in \vartheta(B(0,\epsilon_{0}) \times B(0,2))$ and
$\sup_{B(0,\delta_{0}) \times B(0,2)} |v^{0}| < 1$.
\item $\sup_{B(0,\delta_{0}) \times B(0,\epsilon_{0})} |f| < C_{0} < 1/16$.
\item $1/2 < \sup_{B(0,\delta_{0}) \times B(0,\epsilon_{0})}
|f \circ {\rm exp}(zX)(x,y)|/|f(x,y)| <2$ for all $z \in B(0,2)$.
\end{itemize}
The constant $C_{0}>0$ will be determined later on.
There exists $k \in {\mathbb N}$ such that $(x^{k},y)^{*} X \in \Xt$.
Denote $Y=(x^{k},y)^{*} X$.
Consider $U=B(0,\delta) \times B(0,\epsilon)$ such that
there exists a EV-covering $K_{1}= K_{Y}^{\mu_{1}}$, $\hdots$,
$K_{l}= K_{Y}^{\mu_{l}}$ fulfilling that $H(x)$ is well-defined for
all $x \in [0,\delta) \dot{K}_{p}$,
$H \in Reg(\epsilon, \mu_{p} X, K_{p})$ and $p \in \{1,\hdots,l\}$.
We can also suppose that there exists $C>0$ such that
\[ |f(x,y)| \leq \frac{C}{(1+|\psi_{H,\kappa}^{X}(x,y)|)^{1+ 1 /\nu(X)}}
\ \ \forall (x,y) \in H^{\kappa} \]
for all $H \in Reg(\epsilon, \mu_{p} X,K_{p})$, $p \in \{1,\hdots,l\}$ and
$\kappa \in \{L,R\}$ by proposition \ref{pro:bddconf}.
Finally we suppose that ${\rm exp}(B(0,4)X)(U)$ is contained in
$B(0,\delta_{0}) \times B(0,\epsilon_{0})$.

Denote $V=B(0,\delta) \times B(0,\epsilon_{0}) \times B(0,2)$.
Let $v \in \vartheta(V)$ such that $\sup_{V} |v| < 2$.
Consider an integral $\psi$ of the time form of $X$. We have
\[ X_{v} \left({ \psi - \frac{1}{2 \pi i} \ln t }\right) =
\frac{1}{1+vf} -1 = - \frac{vf}{1+vf} . \]
We obtain
\begin{equation}
\label{equ:phol}
 \psi \circ hol_{v}(x,y,t,z_{0})
= \psi(x,y) + z_{0} - \int_{0}^{z_{0}} \frac{vf}{1+vf}
\circ hol_{v}(x,y,t,z) dz .
\end{equation}
We claim that $hol_{v}(U \times B(0,2) \setminus \{ 0 \}
\times [0,1]) \subset V$.
Otherwise there exist $(x_{0},y_{0},t_{0})$ in $U \times B(0,2)$ and a
minimum $z_{0} \in [0,1]$ such that
$y \circ hol_{v}(x_{0},y_{0},t_{0},z_{0}) \in \partial B(0,\epsilon_{0})$.
This leads us to
\[ |\psi \circ hol_{v}(x_{0},y_{0},t_{0},z_{0}) - \psi(x_{0},y_{0})| \leq
|z_{0}| + |z_{0}| \frac{2 C_{0}}{1-2C_{0}} \leq \frac{8}{7} |z_{0}| <2 \]
and that contradicts the choice of $U$.
Denote $\Delta_{v}(x,y,t) =  \psi \circ hol_{v}(x,y,t,1) - (\psi + 1)$.
We obtain
\[ | \Delta_{v}(x,y,t) | \leq  \frac{32}{7}  |f(x,y)| < 5 C_{0}
\ \ \forall (x,y,t) \in U \times B(0,2) . \]
We define $\Delta_{v}^{1}(x,y)=\Delta_{v}(x,y,1)$ and
$\Delta_{v}^{2}(x,y)= \Delta_{v}(x,y,x)$.
The function $\Delta_{v}^{1}$ is holomorphic in
$U$. The same property is true for
$\Delta_{v}^{2}$ since it is holomorphic in
$U \setminus [x=0]$ and bounded.

We define $\varphi_{1,v}=hol_{v}(x,y,1,1)$ and
$\varphi_{2,v} = {\rm exp}(zX)(x,y,1+\Delta_{v}^{2}(x,y))$.
Clearly $\varphi_{2,v}(x,y) = hol_{v}(x,y,x,1)$ for $x \neq 0$.
\begin{lem}
${\rm exp}(X)$ is a convergent normal form of
$\varphi_{1,v}$, $\varphi_{2,v}$ for all $v$ in $\vartheta(V)$.
\end{lem}
\begin{proof}
The equation \ref{equ:phol} implies that $\Delta_{v}^{1}$ and
$\Delta_{v}^{2}$ belong to $(f)$. Since we have
\[ y \circ \varphi = y + \sum_{j=1}^{\infty}
(1+ \Delta_{\varphi})^{j} \frac{X^{\circ (j)}(y)}{j!} =
y \circ {\rm exp}(X) + O(f^{2}) \]
for $\varphi \in \{\varphi_{1,v} , \varphi_{2,v}\}$
then $\varphi_{1,v}$ and $\varphi_{2,v}$ have convergent
normal form ${\rm exp}(X)$.
\end{proof}
Fix a privileged $\gamma \in Sing_{V} Y$. We choose $C_{0}>0$ such that
there exists $I >0$ holding
that $\forall s \in \{-1,1\}$ and $\forall j \in D_{s}({\rm exp}(X))$ we have
\[ \xi_{\varphi,K_{p}}^{j} \in C^{0}([0,\delta) \dot{K}_{p} \times
[s Img z < -I]) \cap \vartheta
((0,\delta) \dot{K}_{p} \times [s Img z < -I])
\ \forall  1 \leq p \leq l \]
whenever $\varphi$ has convergent normal form ${\rm exp}(X)$
and $|\Delta_{\varphi}(x,y)| \leq 5 \min (C_{0}, |f(x,y)|)$
for all $(x,y) \in B(0,\delta) \times B(0,\epsilon)$ (remark \ref{rem:smalld}).

By choice $({\varphi}_{1,v^{0}})_{|x=0}$ is not analytically trivial.
Thus there exists $(j(0),p(0))$ in ${\mathbb Z} \times \{1,\hdots,l\}$
and $x_{0} \in (\delta/2 ,\delta) \times \dot{K}_{p(0)}$ such that
$\xi_{{\varphi}_{1,v^{0}},K_{p(0)}}^{j(0)}(x_{0},z) \not \equiv z+
\zeta_{\varphi_{1, v^{0}}}(x_{0})$.
Denote $u=(x/x_{0}) v^{0}(y,t)$, we get
$sup_{V} |u| <2$. We define $\varphi_{1}=\varphi_{1,u}$
and $\varphi_{2}=\varphi_{2,u}$.
\begin{lem}
$\varphi_{1}$ is not analytically trivial.
\end{lem}
\begin{proof}
By construction
$\xi_{{\varphi}_{1},K_{p(0)}}^{j(0)}(x,z)$ is well-defined in
$x \in [0,\delta) \times \dot{K}_{p(0)}$ and
$\xi_{{\varphi}_{1},K_{p(0)}}^{j(0)}(x_{0},z) \not \equiv z +
\zeta_{\varphi_{1}}(x_{0})$.
We deduce that $\varphi_{1}$ is not analytically trivial.
\end{proof}
The next lemma is a consequence of $u(0,y,t) \equiv 0$.
\begin{lem}
$(\varphi_{1})_{|x=0} \equiv (\varphi_{2})_{|x=0} \equiv {\rm exp}(X)_{|x=0}$.
In particular
$(\varphi_{1})_{|x=0}$ and $(\varphi_{2})_{|x=0}$ are analytically trivial.
\end{lem}
Denote by $\sigma(x,y)$ the analytic mapping $hol_{u}(x,y,1,ln x/(2 \pi i))$.
\begin{lem}
The mapping $\sigma(x,y)$ conjugates
$\varphi_{1}$ and $\varphi_{2}$ in a domain of the form
$|y| < C_{0}/\sqrt[\nu(X)]{|\ln x|}$ for some $C_{0} \in {\mathbb R}^{+}$.
Moreover $\sigma$ is not univaluated since
\[ \sigma(e^{2 \pi i}x,y) = hol_{u} \left({
x,y,1,\frac{\ln x}{2 \pi i} + 1 }\right) =
hol_{u}(x,y,x,1) \circ \sigma(x,y) = \varphi_{2} \circ \sigma(x,y) . \]
\end{lem}
\begin{proof}
Consider a domain $W \subset B(0,\delta) \times B(0,\epsilon_{0})$ such that
\[ {\rm exp}
\left({ B \left({ 0, \frac{|\ln x|}{\pi} }\right) X }\right)(x,y) \in
B(0,\delta) \times B(0,\epsilon_{0}) \ \ \forall (x,y) \in W. \]
Since $y \circ hol_{u}(x,y,1,w \ln x/(2 \pi)) \subset
\overline{B}(0,\epsilon_{0})$ for all $w \in [0,1]$ implies
\begin{equation}
\label{equ:repmod}
 \left|{ \psi \circ hol_{u}
\left({ x,y,1, w \frac{\ln x}{2 \pi i} }\right) - \psi(x,y) }\right| \leq
w \frac{|\ln x|}{2 \pi} + \frac{w}{7} \frac{|\ln x|}{2 \pi}
< \frac{|\ln x|}{\pi}
\end{equation}
by equation \ref{equ:phol}
then $hol_{u}(x,y,1,w \ln x/(2 \pi i))$ is well-defined and belongs
to $V$ for all $(x,y,w) \in W \times [0,1]$.
We have $\psi \sim 1/y^{\nu(X)}$ in the first exterior set by
remark \ref{rem:bditfpet}, we can deduce that $W$ contains a domain
of the form $|y| < C_{0}/\sqrt[\nu(X)]{|\ln x|}$ for some
$C_{0} \in {\mathbb R}^{+}$.
\end{proof}
The domain $W_{0}=[|y| < C_{0}/\sqrt[\nu(X)]{|\ln x|}]$ contains the germ
of all the ``algebraic'' domains of the form $|y| < |x|^{b}$ for
$b \in {\mathbb Q}^{+}$, in particular $W_{0}$ contains
$Sing X \setminus \{(0,0)\}$, every intermediate set and every exterior set
except the first one.
\begin{lem}
We have
\[ \xi_{\varphi_{2}, K_{p}}^{j}(x_{0},z) =
\left({ z +  \frac{\ln x_{0}}{2 \pi i} }\right) \circ
\xi_{\varphi_{1},K_{p}}^{j}(x_{0},z) \circ
\left({ z -\frac{\ln x_{0}}{2 \pi i} }\right) \]
for all $(j,p) \in {\mathbb Z} \times \{1,\hdots,l\}$ and
$x_{0} \in (0,\delta) \times \dot{K}_{p}$. Then
we get $m_{\varphi_{1}} \equiv m_{\varphi_{2}}$ and
$\varphi_{1} \not \sim \varphi_{2}$.
\end{lem}
\begin{proof}
Let $(x_{0},y_{0}) \in Sing X \setminus \{(0,0)\}$. We remark that
\[ \lim_{(x,y) \to (x_{0},y_{0})}
\psi \circ hol_{v} \left({ x,y,1,w \frac{\ln x}{2 \pi i} }\right)
- \psi = w \frac{\ln x_{0}}{2 \pi i} \]
for all $w \in [0,1]$ by equation \ref{equ:phol}. Basically the uniform
hypothesis in proposition \ref{pro:conimeq} is used to estimate
$\psi \circ \kappa - \psi$ for a r-conjugation $\kappa$.
Such an estimation is provided here by the inequality \ref{equ:repmod}, hence
we can proceed like in proposition \ref{pro:conimeq} to obtain
\[ \xi_{\varphi_{2}, K_{p}}^{j}(x_{0},z) =
\left({ z +  \frac{\ln x_{0}}{2 \pi i} }\right) \circ
\xi_{\varphi_{1},K_{p}}^{j}(x_{0},z) \circ
\left({ z -\frac{\ln x_{0}}{2 \pi i} }\right) \]
for all $(j,p) \in {\mathbb Z} \times \{1,\hdots,l\}$ and
$x_{0} \in (0,\delta) \times \dot{K}_{p}$. We deduce
$m_{\varphi_{1}} \equiv m_{\varphi_{2}}$ from the previous equation and
$(\varphi_{1})_{|x=0} \equiv (\varphi_{2})_{|x=0}$.

We know that $\varphi_{1}$ and $\varphi_{2}$ are not analytically trivial.
Therefore we have $\varphi_{1} \not \sim \varphi_{2}$;
otherwise $|\ln|x||$  would be bounded in a neighborhood of
$0$ by theorem \ref{teo:imgbdd}.
%
%
\end{proof}
\bibliography{rendu}
\end{document}